\documentclass[10pt]{amsart}
\usepackage{amssymb,amsmath}
\usepackage{hyperref}
\usepackage{bm}
\usepackage[]{graphicx}
\usepackage{amssymb, epsfig}
\newcommand{\gammaa}{\gamma}

\newcommand{\ssy}{\scriptscriptstyle}
\newcommand{\cunit}{\text{\rm i}}
\newtheorem{thm}{Theorem}[section]

\newtheorem{lem}[thm]{Lemma}

\newtheorem{proposition}[thm]{Proposition}
\newtheorem{remm}{Remark}[section]
\newtheorem{remark}[remm]{Remark}
\numberwithin{equation}{section}
\begin{document}
\title[]
{GALERKIN METHODS FOR PARABOLIC AND \\
SCHR\"ODINGER EQUATIONS WITH DYNAMICAL BOUNDARY CONDITIONS \\
AND APPLICATIONS TO UNDERWATER ACOUSTICS}
\author{D.C.~Antonopoulou$^{\dag\P}$,
V.A.~Dougalis$^{\ddag\P}$ and G.E.~Zouraris$^{\S\P}$}
\thanks
{$^{\ddag}$ Department of Mathematics,
University of Athens,
Panepistimiopolis,
GR--157 84 Zographou, Greece.}
\thanks
{$^{\dag}$ Department of Applied Mathematics,
University of Crete,
GR--714 09 Heraklion, Greece.}
\thanks
{$^{\S}$ Department of Mathematics,
University of Crete,
GR--714 09 Heraklion, Greece.}
\thanks
{$^{\P}$ Institute of Applied and Computational Mathematics,
FO.R.T.H., GR--711 10 Heraklion, Greece.}
\maketitle
%
%
%
%
\begin{abstract}
In this paper we consider Galerkin-finite element methods that
approximate the solutions of initial-boundary-value problems in
one space dimension for parabolic and Schr\"odinger evolution
equations with dynamical boundary conditions. Error estimates of
optimal rates of convergence in $L^2$ and $H^1$ are proved for
the accociated semidiscrete and fully discrete
Crank-Nicolson-Galerkin approximations. The problem involving the
Schr\"odinger equation is motivated by considering the standard
`parabolic' (paraxial) approximation to the Helmholtz equation,
used in underwater acoustics to model long-range sound
propagation in the sea, in the specific case of a domain with a
rigid bottom of variable topography. This model is contrasted
with alternative ones that avoid the dynamical bottom boundary
condition and are shown to yield qualitatively better
approximations. In the (real) parabolic case, numerical
approximations are considered for dynamical boundary conditions
of reactive and dissipative type.
\end{abstract}
%
%
%
\section{Introduction}\label{sec1}
Our main goal in this paper is to analyze Galerkin-finite element
methods for initial-boundary-value problems, involving dynamical
boundary conditions, for the linear Schr\"odinger and the heat
equations. In addition, in a specific problem arising in
underwater acoustics and modelled by the Schr\"odinger equation,
we will also consider an alternative boundary condition and
evaluate, analytically and numerically, the two models.
\par
We start with the underwater acoustic application. Consider the
Helmholtz equation (HE) in cylindrical coordinates in the
presence of cylindrical symmetry
\begin{equation}\tag{HE}
\Delta p+k_0^2\,\eta^2(r,z)p=0.
\end{equation}
Here $z\geq 0$ is the depth variable increasing downwards and
$r\geq 0$ is the horizontal distance (range) from a harmonic point
source of frequency $f_0$ placed on the $z$ axis. For simplicity
we shall assume that the medium consists of a single layer of
water of constant density, occupying the region, $0\leq z\leq
\ell(r)$, $r\geq 0$, between the free surface $z=0$ and the
range-dependent bottom $z=\ell(r)$ (see Fig. 1); $\ell=\ell(r)$
will be assumed to be smooth and positive. The function
$p=p(r,z)$ is the acoustic pressure, $k_0=\frac{2\pi f_0}{c_0}$
is a reference wave number, $c_0$ a reference sound speed, and
$\eta(r,z)$ the index of refraction, defined as
$\frac{c_0}{c(r,z)}$, where $c(r,z)$ is the speed of sound in the
water. (HE) is supplemented by the surface `pressure-release'
condition $p(r,0)=0$. In the case of a soft bottom the
homogeneous Dirichlet boundary condition
\begin{equation}\tag{D}
p=0\quad\mbox{\rm at}\quad z=\ell(r)
\end{equation}
is assumed to hold. The case of a rigid bottom is modelled by a
Neumann boundary condition (with ${\dot\ell}=\frac{d\ell}{dr}$)
\begin{equation}\tag{N}
p_z-{\dot\ell}(r)p_r=0\quad\mbox{\rm at}
\quad z=\ell(r).
\end{equation}
Applying the change of variables $p(r,z)=\psi(r,z)
\,\frac{e^{\cunit k_0 r}}{\sqrt{k_0r}}$, assuming that $|2\cunit
k_0\psi_r|\!\!>\!\!>\!\!|\psi_{rr}|$ (narrow-angle paraxial
approximation) and neglecting terms of $O(\frac{1}{r^{2}})$
(far-field approximation) we arrive (cf., e.g., \cite{ref54},
\cite{ref42}, \cite{ref12}) at the standard `Parabolic' Equation
(PE), which is a linear Schr{\"o}dinger equation of the form
\begin{equation}\tag{PE}
\psi_r=\tfrac{\cunit}{2k_0}\,\psi_{zz}
+\cunit\,\tfrac{k_0}{2}\,(\eta^2(r,z)-1)\,\psi,
\end{equation}
where $\psi=\psi(r,z)$ is a complex-valued function of the two
real variables $r$ and $z$. The (PE) has been widely used in
underwater acoustics to model one-way, long-range sound
propagation near the horizontal plane of the source, in
inhomogeneous, weakly range-dependent marine environments. Its
solution will be sought in the domain $0\leq z\leq \ell(r)$,
$r\geq 0$. The (PE) will be supplemented by an initial condition
$\psi(0,z)=\psi_0(z)$, $0\leq z\leq\ell(0)$, modelling the source
at $r=0$, the surface boundary condition $\psi=0$ for $z=0$,
$r\geq 0$, and a bottom boundary condition obtained by
transforming (D) or (N). The Dirichlet boundary condition (D)
remains of the same type ($\psi=0$ at $z=\ell(r)$) while the
Neumann boundary condition (N) is transformed to a condition of
the form
\begin{equation}\tag{PN}
\psi_z-{\dot\ell}(r)\,\psi_r-g_{\ssy B}(r)
\,{\dot\ell}(r)\,\psi=0
\quad\mbox{\rm at}\quad z=\ell(r),
\end{equation}
where $g_{\scriptscriptstyle B}(r)$ is complex-valued and is
usually taken simply as $\cunit\, k_0$.

The theory and numerical analysis of this initial-boundary-value
problem (ibvp) with the Dirichlet bottom boundary condition is
standard, cf.e.g. \cite{ref46}, \cite{ref88}, and will not be
considered any further. On the other hand the analysis is
complicated in the case of the Neumann boundary condition, when
${\dot\ell}(r)$ is not the zero function, due to the presence of
the term $\psi_r$ in (PN). In \cite{ref1} Abrahamsson and Kreiss
proved existence and uniqueness of solutions, in the case of a
strictly monotone bottom, i.e. when ${\dot\ell}(r)$ is of one
sign for $r\geq 0$.
\begin{center}
\resizebox{8cm}{4cm}{\includegraphics{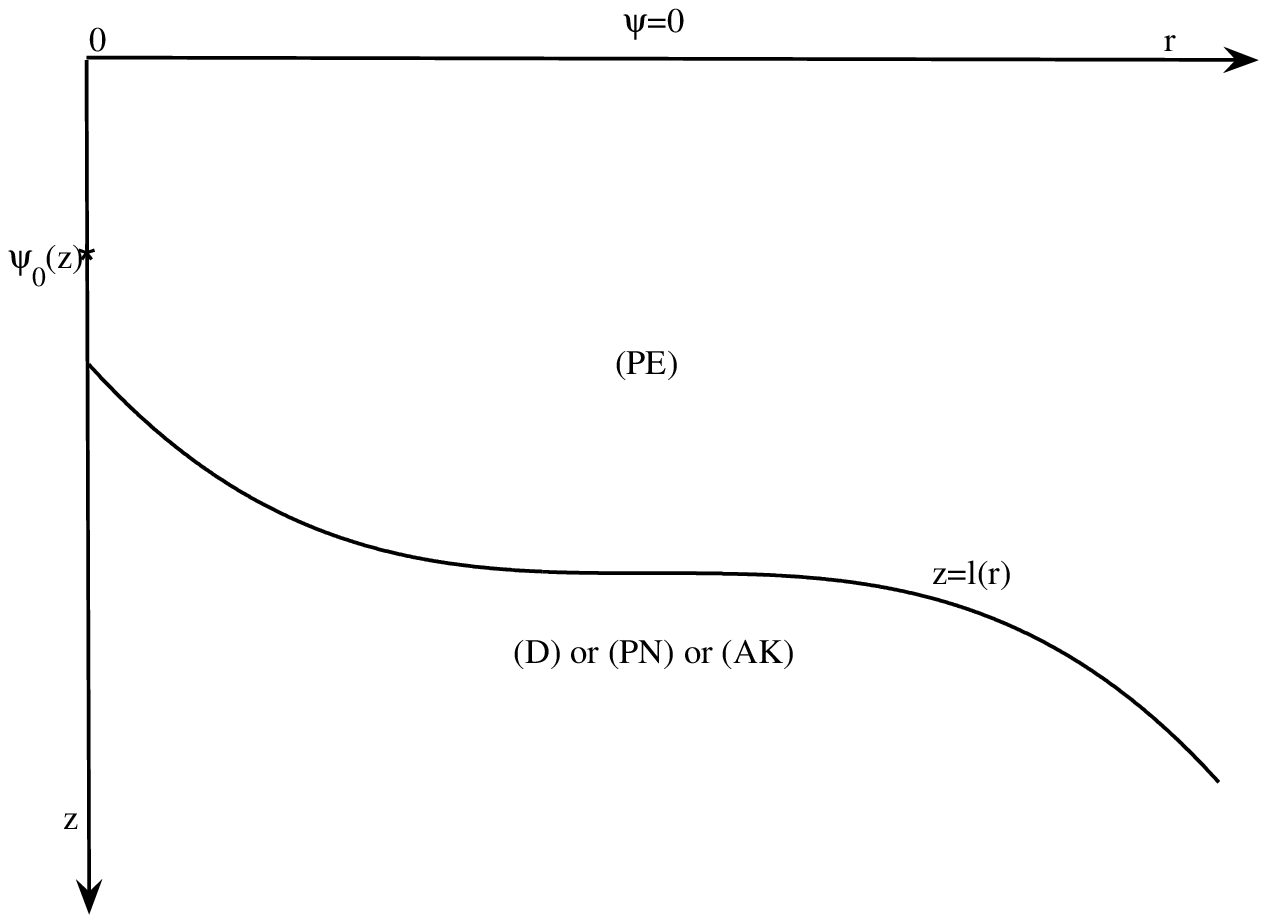}}\\
{\small {Figure 1. The domain of the initial boundary value
problems for the (PE) in the $r$, $z$ variables}}
\end{center}
\vspace{0.5cm}
\par
We shall transform the above ibvp's to equivalent ones posed on a
horizontal strip. With this aim in mind, we first introduce
non-dimensional variables as in \cite{ref8}, defined by
$y:=\frac{z}{L}$, $t:=\frac{r}{L}$, $w:=\frac{\psi}{\psi_{\rm
ref}}$, where we take $L:=\frac{1}{k_0}$ and $\psi_{\rm
ref}:=\max|\psi_0|$. Then, letting
$s(t):=k_0\,\ell\big(\frac{t}{k_0}\big)$,
$g(t)=k_0\,g_{\scriptscriptstyle B}\big(\frac{t}{k_0}\big)$,
$\gammaa(t,y):=\frac{1}{2}\,\big[\eta^2\big(\frac{t}{k_0},
\frac{y}{k_0}\big)-1\,\big]$ we see that the (PE) becomes
\begin{equation}\label{1.1}
w_t=\tfrac{\cunit}{2}\,w_{yy}+\cunit\,\gammaa(t,y)\,w, \quad
0\leq y\leq s(t),\quad t\geq 0.
\end{equation}
We note that the index of refraction $\eta$, and consequently the
function $\gamma$, may be taken to be complex-valued in order to
model attenuation of sound in the water. The initial condition
becomes
\begin{equation}\label{1.2}
w(0,y)=w_0(y):=\tfrac{1}{\psi_{\rm
ref}}\,\psi_0\big(\tfrac{y}{k_0}\big),\quad 0\leq y\leq s(0).
\end{equation}
The surface condition remains the same, i.e.,
\begin{equation}\label{1.3}
w(t,0)=0,\quad t\geq 0,
\end{equation}
while the boundary condition (PN) becomes
\begin{equation}\label{1.4}
w_y(t,s(t))-\dot{s}(t)\,\big[w_t(t,s(t))+g(t)\,w(t,s(t))\big]=0,
\quad t\geq 0.
\end{equation}
We now perform the range-dependent change of depth variable
$x:=\frac{y}{s(t)}$, that maps the domain of the problem onto the
horizontal strip $0\leq x\leq 1$, $t\geq 0$. We also make the
transformation
\begin{equation}\label{1.5}
u(t,x)=\exp(-\cunit\,\delta(t)x^2)\,w(t,s(t)x),
\end{equation}
which defines the new field variable $u(t,x)$ for $0\leq x\leq
1$, $t\geq 0$. In \eqref{1.5}
$\delta(t):=\frac{\dot{s}(t)s(t)}{2}$, $t\geq 0$, where a dot
denotes differentiation with respect to $t$. In terms of the new
variables \eqref{1.1} becomes
\begin{equation}\label{1.6}
u_t=\cunit\,a(t)\,u_{xx}+\cunit\,\beta(t,x)u, \quad 0\leq x\leq
1,\quad t\geq 0,
\end{equation}
where, for $0\leq x\leq 1$, $t\geq 0$,
\begin{equation}\label{1.7}
\gathered
a(t)=\tfrac{1}{2s^2(t)},\quad \beta(t,x)=
\beta_{\scriptscriptstyle R}(t,x)
+{\rm i}\beta_{\scriptscriptstyle I}(t,x),\\
\beta_{\scriptscriptstyle R}(t,x)={\rm Re}[\gamma(t,xs(t))]
-\tfrac{\ddot{s}(t)s(t)x^2}{2},\quad
\beta_{\scriptscriptstyle I}(t,x)={\rm
Im}[\gamma(t,xs(t))]+\tfrac{\dot{s}(t)}{2\,s(t)}\cdot
\endgathered
\end{equation}
The purpose of introducing in \eqref{1.5} the factor
$e^{-\cunit\,\delta(t)x^2}$ with $\delta=\frac{\dot{s}s}{2}$ is
to avoid the presence of a $u_x$ term in the right-hand side of
\eqref{1.6} and, consequently, simplify somewhat the analysis.
Under the transformation \eqref{1.5}, the initial and boundary
conditions \eqref{1.2}-\eqref{1.4} change accordingly.
Specifically, we have
\begin{equation}\label{1.8}
u(0,x)=u_0(x):=e^{-\cunit\delta(0)x^2}w_0(xs(0))
\quad\forall\,x\in[0,1],
\end{equation}
\begin{equation}\label{1.9}
u(t,0)=0,\quad t\geq 0,
\end{equation}
and
\begin{equation}\label{1.10}
u_x(t,1)=s_1(t)\,u_t(t,1)+m(t)\,u(t,1),
\quad t\geq 0,
\end{equation}
where
\begin{equation}\label{1.11}
s_1(t):=\frac{\dot{s}(t)s(t)}{1+\dot{s}(t)^2},
\quad m(t):=g(t)s_1(t)+{\rm
i}(s_1(t)\dot{\delta}(t)-2\delta(t)),\quad t\geq 0.
\end{equation}
\par
The boundary condition \eqref{1.10} is an example of a
\textit{dynamical} boundary condition, because it involves (if
$\dot{s} \neq 0$) the value of $u_t$ at the boundary. As was
already mentioned, the well-posedness of ibvp's of the type
$\{$\eqref{1.6}, \eqref{1.8}, \eqref{1.9}, \eqref{1.10}$\}$, for
$t$ in a finite interval $[0,T]$, was proved in \cite{ref1} under
the assumption that $\dot{s}(t)$ is of one sign for all
$t\,\in\,[0,T]$. One of our main purposes in this paper is to
construct and analyze fully discrete Galerkin-finite element
methods for approximating the solution of the above ibvp.
\par
We consider the ibvp consisting of \eqref{1.6}-\eqref{1.11}. We
assume that the bottom is \textit{upsloping}, i.e. that
$\dot{s}(t)\leq 0$, and that the problem has a unique solution,
smooth enough for the purposes of the error estimation. In
paragraphs 2.1 and 2.2 we discretize the problem in $x$ by the
standard Galerkin method and prove optimal-order $L^2$ and $H^1$
estimates for the error of the resulting semidiscretization. This
is achieved by using appropriate properties of the $L^2$ and the
elliptic projections onto the finite element subspace and a
relevant $H^1$ superconvergence result. (The difficulty of the
problem lies in the presence of the $u_t$ term in \eqref{1.10};
the condition $\dot{s}(t)\leq 0$, which implies that $s_1(t)\leq
0$, is needed to obtain a basic energy inequality for the error
of the semidiscretization). Subsequently, in paragraph 2.2, we
discretize the semidiscrete problem in the $t$ variable using a
Crank-Nicolson type method with a variable step-length. Again,
under the assumption that $\dot{s}(t)\leq 0$ for $0\leq t\le T$,
we prove $L^2$ and $H^1$ error estimates which are of optimal
order in $x$ and $t$.
\par
In order to overcome the analytical and numerical difficulties
caused by dynamical boundary conditions of the form \eqref{1.10}
Abrahamsson and Kreiss proposed in \cite{ref2} an alternative
rigid bottom boundary condition, which, in the case of (PE), is
of the form
$$\psi_z-\cunit\,k_0\,{\dot\ell}(r)\,\psi=0\quad\mbox{\rm at}\quad
z=\ell(r).\leqno{{\rm (AK)}}$$
This condition may be viewed as a `paraxialization' of (PN). When
the nondimensionalization $z\rightarrow y$, $r\rightarrow t$,
$\psi\rightarrow w$ is performed, (AK) becomes
\begin{equation}\label{1.12}
w_y(t,s(t))-\cunit\,\dot{s}(t)\,w(t,s(t))=0,\quad t\geq 0.
\end{equation}
Finally, after changing the depth variable by $x=y/s(t)$ and the
dependent variable by \eqref{1.5}, it is not hard to see that
\eqref{1.12} becomes simply
\begin{equation}\label{1.13}
u_x(t,1)=0,\quad t\geq 0.
\end{equation}
The proof of the well-posedness of the ibvp consisting of
\eqref{1.6}-\eqref{1.9} and \eqref{1.13} is standard, cf.
\cite{ref46}. Its numerical analysis too is straightforward;
under no restriction on the sign of $\dot{s}(t)$ we prove in
paragraph 2.3 optimal-order $L^2$ and $H^1$ error estimates for
the standard semidiscrete Galerkin scheme and its Crank-Nicolson
full discretization.
\par
In Section~\ref{sec55} we present results of various numerical
experiments that we performed for problems on variable domains
with the Neumann and Abrahamsson-Kreiss bottom boundary
conditions, using the fully discrete finite element methods
analyzed in Section~\ref{sec3}. As predicted by the theoretical
stability and convergence analysis, the finite element scheme is
stable and second-order accurate when Neumann boundary conditions
are considered in domains with upsloping bottoms. (It also
appears to be convergent in small scale problems with downsloping
bottoms and also in more realistic examples if the downsloping
bottom has very small slope.) The scheme with the
Abrahamsson-Kreiss condition behaved well, as predicted by the
theory, in all examples of bottoms of arbitrary shape that we ran.
\par
When we compared the results of the schemes using both boundary
conditions in the case of the upsloping and downsloping rigid
bottom ASA wedge (a standard test problem for long range sound
propagation in underwater acoustics, \cite{JF}), we found that in
the upsloping case there was very good agreement between the two
schemes. In the downsloping case, the scheme implementing the
Neumann boundary condition was not convergent. This is in
agreement with the results of Abrahamsson and Kreiss,
\cite{ref1}, \cite{ref2}, who pointed out that for some
downsloping bottom profiles one may observe instabilities in the
case of the Neumann boundary condition. On the other hand, the
scheme with the Abrahamsson-Kreiss condition was convergent and
its results agreed well with those furnished by the finite
difference code IFD, \cite{LB}, \cite{LBP}, \cite{ref42},
implemented with the rigid bottom boundary condition option. The
IFD scheme uses a discretized version of the Neumann boundary
condition (PN), wherein the $\psi_r$ term is replaced by the
right-hand side of the (PE). We prove \textit{a priori} $L^2$
estimates for the resulting ibvp.
\par
A final point of interest emerging from the numerical experiments
is that, for some downsloping bottom profiles $s(t)$ with an
inflection point at some $t=t^*$, we observed violent growth of
the $L^2$-norm of the numerical solution of the problem with the
Neumann boundary condition for $t>t^*$. This growth (blow-up?) of
the solution seems to be a feature of the problem and not an
artifact of the numerical scheme.
\par
Error estimates for a \textit{finite difference} scheme of
second-order of accuracy in $x$ and $t$ for some of the ibvp's
considered here were proved in \cite{ref8}. In the case of the
Neumann boundary condition \eqref{1.10} these error estimates
were shown to hold not only when $\dot{s}(t)\leq 0$ but also in
the strictly downsloping case $\dot{s}(t)>0$, $t\in [0,T]$, as a
result of the validity of a certain discrete $H^1$ estimate; this
estimate mimics an analogous $H^1$ estimate for the continuous
problem, which holds provided $\dot{s}(t)\leq 0$ or $\dot{s}(t)>0$
when $t\in [0,T]$. In \cite{ref53} Sturm considered the
Abrahamsson-Kreiss condition for the (PE) in three dimensions
over a variable bottom in the more general case of a multilayered
fluid medium with homothetic layers. When restricted to single
layer problems in the presence of azimuthal symmetry, the scheme
of \cite{ref53} is similar to the one analyzed here in the case
of the Abrahamsson-Kreiss bottom boundary condition. We have
considerably modified the analysis of \cite{ref53} and obtain
optimal-order estimates, since, by using the transformation
\eqref{1.5}, we essentially avoid an elliptic projection with
time-dependent terms. We finally mention that a uniform range
step version of the scheme of this paper and also
three-dimensional extensions thereof were analyzed in
\cite{refthes}.
\par
The problem addressed in the present paper, namely sound
propagation modelled by the (PE) in a single layer of water over
a rigid bottom, is, of course, an idealized model problem in
underwater acoustics. More realistic environments consist, for
example, of a layer of water above several layers of fluid
sediments of different density, speed of sound and attenuation
overlying a rigid or soft bottom. If the layers are separated by
interfaces of weakly range-dependent topography and low
backscatter is expected, long-range sound propagation may again
be modelled by the (PE) in each layer with transmission
conditions (continuity of $\psi$ and of
$\frac{1}{\rho}\frac{\partial \psi}{\partial n}$, where $\rho$ is
the density and $n$ the normal direction to the interface)
imposed across the layer interfaces. Hence, the issue arises of
how to treat the dynamical interface condition, now involving
$\psi_r$  on both sides of an interface, and the ensuing problems
are analogous to those encountered in the case of the dynamical
bottom boundary condition. The analysis is more complicated now,
as it appears that possible non-homotheticity of the layers has
to be balanced by the jump across the interface in the imaginary
part of the analog of the function $\gamma$, cf. \eqref{1.1}, in
order to ensure the well-posedeness of the problem, \cite{DZ}.
For a recent review of several issues regarding the interface
problem for the (PE), we refer the reader to \cite{DKSZ};
references to underwater acoustics computations with the (PE) in
the presence of interfaces with change-of-variable techniques
include e.g. \cite{ref8}, \cite{ref53} and \cite{DKSZ}. Here we
just wish to point out that range-dependent topography has often
been approximated in practice by `staircase' (piecewise
horizontal) bottoms and interfaces. This raises the issue of what
boundary / interface conditions to pose on the vertical part of
the steps of the staircase. Moreover, it is well documented that
staircase approximations lead to nonphysical energy losses or
gains, cf. e.g. \cite{JF}, \cite{PJF}. To alleviate this problem
of energy non-conservation, change-of-variable techniques may be
used as in the present paper. They may also be extended to
interface, \cite{DKSZ}, or 3D-problems, \cite{ref53}.
\par
We turn now to one-dimensional (real) parabolic problems with
dynamical boundary conditions. We consider the following model
problem: For $0<T<\infty$ we seek a real-valued function
$u=u(t,x)$ defined for $(t,x)\in [0,T]\times [0,1]$ and satisfying
\begin{equation}\label{1.14}
\begin{split}
&u_t=a(t)u_{xx}+\beta(t,x)u+f(t,x)\quad\forall (t,x)\in\,[0,T]\times[0,1],\\
&u(t,0)=0\quad\forall t\in\,[0,T],\\
&a(t)u_x(t,1)=\varepsilon(t)u_t(t,1)+\delta(t)u(t,1)+g(t)\quad\forall t\in\,[0,T],\\
&u(0,x)=u_0(x)\quad\forall x\in\,[0,1],
\end{split}
\end{equation}
where $a(t)\geq a_*>0$ for $t\,\in\,[0,T]$ and $\beta$, $f$,
$\varepsilon$, $\delta$, $g$, $u_0$ are smooth, real-valued
functions. Such problems occur in heat conduction, (\cite{C},
Section 4.3.5), and in other areas; see \cite{Esc1993} for a
fuller list of references. Our aim is to construct fully discrete
Galerkin-finite element approximations for the ibvp \eqref{1.14}
and prove error estimates, with techniques analogous to those used
in the case of the Schr\"odinger equation. We consider two
different cases depending on the sign of the function
$\varepsilon$ in the dynamical boundary condition.
\par
We treat first the \textit{dissipative} case, characterized by
the hypothesis that $\varepsilon(t)\leq 0$ for all
$t\,\in\,[0,T]$, and in which the ibvp \eqref{1.14} is well
posed, cf. e.g. \cite{Esc1993}. In paragraph 4.1, applying the
standard Galerkin method to this case, we prove optimal-order
$L^2$ and $H^1$ estimates for the error of the resulting
semidiscretization and for the Crank-Nicolson-Galerkin fully
discrete scheme. Matters are more complicated in the
\textit{reactive} case, wherein $\varepsilon(t)>0$ for
$t\,\in\,[0,T]$. In this case the problem is well posed in one
space dimension as in the case at hand, but in general is not
well posed in higher dimensions, \cite{VazVit2008}, \cite{BBR}.
To construct a Galerkin-finite element method in this case, we
replace the term $u_t$ in the dynamical boundary condition using
the p.d.e. in \eqref{1.14}, thus obtaining a boundary condition
involving $u_{xx}(t,1)$. The resulting ibvp is discretized in
space by means of a $H^1$-type Galerkin method that uses finite
element spaces consisting of piecewise polynomial functions in
$H^2$ of degree at least three. In paragraph 4.2 we analyze this
method and prove optimal-order $H^1$ error estimates for the
semidiscrete approximation and the fully discrete one when the
Crank-Nicolson scheme is used in time-stepping. The case where
$\varepsilon(t)$ changes sign in $[0,T]$ is under investigation;
for a discussion cf. \cite{Ba-Re}.
\section{Numerical Schemes and Error Estimates for the (PE)}\label{sec3}
\subsection{Preliminaries}\label{subsec40}
Let $D:=(0,1)$. We will denote by $L^2(D)$ the space of the
Lebesgue measurable complex-valued functions which are square
integrable on $D$, and by $\|\cdot\|$ the standard norm of
$L^2(D)$, i.e., $\|f\|:=\{\int_{\ssy
D}|f(x)|^2\,dx\}^{\frac{1}{2}}$ for $f\in L^2(D)$.
The inner product in $L^2(D)$ that induces the norm $\|\cdot\|$
will be denoted by $(\cdot,\cdot)$, i.e. $(f_1,f_2):=\int_{\ssy D}
f_1(x)\,\overline{f_2(x)}\,dx$ for $f_1$, $f_2\in L^2(D)$.
Also, we will denote by $L^{\infty}(D)$ the space of the Lebesgue
measurable functions which are bounded a.e. on $D$, and by
$|\cdot|_{\infty}$ the associated norm, i.e.,
$|f|_{\infty}:=\text{\rm ess}\sup_{\ssy D}|f|$ for $f\in
L^{\infty}(D)$.
For $s\in{\mathbb N}_0$, we denote by $H^s(D)$ the Sobolev space
of complex-valued functions having generalized derivatives up to
order $s$ in $L^2(D)$, and by $\|\cdot\|_s$ its usual norm, i.e.
$\|f\|_s:=\bigl\{\sum_{\ell=0}^s
\|\partial_x^{\ell}f\|^2\bigr\}^{\frac{1}{2}}$ for $f\in H^s(D)$.
In addition, we set $|v|_1:=\|v'\|$ for $v\in H^1(D)$.
Also, ${\mathbb H}^1(D)$ will denote the subspace of $H^1(D)$
consisting of functions which vanish at $x=0$ in the sense of
trace$;$ we set ${\mathbb H}^s(D)=H^s(D)\cap{\mathbb H}^1(D)$ for
$s\ge2$.
In addition, for $s\in{\mathbb N}_0$, we denote by
$W^{s,\infty}(D)$ the Sobolev space of complex-valued functions
having generalized derivatives up to order $s$ in $L^{\infty}(D)$,
and by $|\cdot|_{s,\infty}$ its usual norm, i.e.
$|f|_{s,\infty}:=\max_{0\leq{\ell}\leq{s}}
|\partial_x^{\ell}f|_{\infty}$ for $f\in W^{s,\infty}(D)$. In
what follows, $C$ will denote a generic constant independent of
the discretization parameters and having in general different
values at any two different places.
\par
For later use, we recall the well-known Poincar{\'e}-Friedrichs
inequality
\begin{equation}\label{PoincareF}
\|v\|\leq\,C_{\ssy P\!F}\,|v|_1\quad\forall\,v\in{\mathbb H}^1(D),
\end{equation}
the Sobolev-type inequality
\begin{equation}\label{SobolevIneq}
|v|_{\infty}\leq\,|v|_1\quad\forall\,v\in {\mathbb H}^1(D)
\end{equation}
and the trace inequality
\begin{equation}\label{TraceIneq}
|v(1)|^2\leq\,2\,\|v\|\,|v|_1 \quad\forall\,v\in{\mathbb H}^1(D).
\end{equation}
\par
Let $r\in{\mathbb N}$ and $S_h$ be a finite dimensional subspace
of ${\mathbb H}^1(D)$ consisting of complex-valued functions that
are polynomials of degree less or equal to $r$ in each interval
of a non-uniform partition of $D$ with maximum length
$h\in(0,h_{\star}]$. It is well-known, \cite{ref14}, that the
following approximation property holds:
\begin{equation}\label{2.35}
\displaystyle{\inf_{\chi\in S_h}}\big\{
\|v-\chi\|+h\|v-\chi\|_1\big\}\leq
C\,h^{s+1}\,\|v\|_{s+1},\quad\forall\,v\in{\mathbb H}^{s+1}(D),
\,\,\,\forall\,h\in(0,h_{\star}],\quad s=0,\dots,r.
\end{equation}
Also, we assume that the following inverse inequality holds
\begin{equation}\label{2.36}
|\phi|_1\leq\,C\,h^{-1} \,\|\phi\|\quad\forall\,\phi\in
S_h,\quad\forall\,h\in(0,h_{\star}],
\end{equation}
which is true when, for example, the partition of $D$ is
quasi-uniform, \cite{ref14}.
In addition, we define the $L^2-$projection operator
$P_h:L^2(D)\rightarrow S_h$ by
\begin{equation*}
(P_hv,\phi)=(v,\phi)\quad\forall\,\phi\in S_h,\quad\forall\,v\in
L^2(D),
\end{equation*}
and the elliptic projection operator $R_h:H^1(D)\rightarrow S_h$
by
\begin{equation}\label{R_def}
{\mathcal B}(R_hv,\phi)={\mathcal B}(v,\phi)\quad\forall\,\phi\in
S_h, \quad\forall\,v\in H^1(D),
\end{equation}
where $\mathcal B$ is the sesquilinear form defined for $u$, $w$
$\in\, H^1(D)$ by $\mathcal B(u,\,w):=(u',\,w')$. It follows,
\cite{ref14}, \cite{ref55}, that
\begin{equation}\label{4.39}
\|R_hv-v\|+h\,\|R_hv-v\|_1\leq\,C\,h^{s+1}
\,\|v\|_{s+1},\quad\forall\,v\in{\mathbb H}^{s+1}(D),
\,\,\,\forall\,h\in(0,h_{\star}],\quad s=0,\dots,r.
\end{equation}
\par
Finally, for $v\in L^2(D)$, we define the discrete negative norm
\begin{equation*}
\|v\|_{-1,h}:=\sup\Big\{\tfrac{|(v,\phi)|}{|\phi|_1}:\quad\phi\in
S_h\quad{\rm
and}\quad\phi\not=0\Big\},\quad\forall\,h\in(0,h_{\star}].
\end{equation*}
%
%
%
%
%
%
\begin{lem}\label{2.3.1}
The elliptic projection operator $R_h$ has the following property:
\begin{equation}\label{4.40}
R_hv(1)=v(1),\quad\forall\,v\in \mathbb{H}^1(D).
\end{equation}
\end{lem}
%
%
%
%
%
%
%
%
%
%
%
%
\begin{proof}
Let $v\in \mathbb{H}^1(D)$ and $\omega$ be the element of $S_h$
given by $\omega(x)=x$ for $x\in{\overline D}$. Then \eqref{R_def}
gives $R_hv(1)-v(1)={\mathcal B}(R_hv-v,\omega)=0$, which is the
desired result. $\Box$
\end{proof}
%
%
%
%
%
%
\begin{lem}\label{2.3.2}
Let $\omega\in C^1({\overline D})$. Then
\begin{equation}\label{TheBound1}
|P_h(\omega\phi)|_1\,\leq\,C\,|\omega|_{1,\infty}
\,|\phi|_1\quad\forall\,\phi\in
S_h,\quad\forall\,h\in(0,h_{\star}].
\end{equation}
\end{lem}
%
%
%
%
%
\begin{proof}
Let $h\in(0,h_{\star}]$ and $\phi\in S_h$.
Since $|P_h(\omega\phi)|_1\leq|P_h(\omega\phi-R_h(\omega\phi))|_1
+|R_h(\omega\phi)|_1$, using \eqref{2.36} and \eqref{R_def} we
arrive at
$|P_h(\omega\phi)|_1\leq\,C\,h^{-1}\,\|\omega\phi-R_h(\omega\phi)\|
+|\omega\phi|_1$. Next, we use the estimate \eqref{4.39} for
$s=0$ to obtain
$|P_h(\omega\phi)|_1\leq\,C\,\big[\,|\omega|_{\infty}\,|\phi|_1
+|\omega'|_{\infty}\|\phi\|\,\big]$. Thus, the bound
\eqref{TheBound1} follows by combining the latter inequality and
\eqref{PoincareF}. $\Box$
\end{proof}
%
%
%
%
%
%
\subsection{The Neumann (dynamical) boundary condition}
In this subsection, we shall consider the (PE) with the Neumann
boundary condition, i.e. the ibvp \eqref{1.6}, \eqref{1.8},
\eqref{1.9}, \eqref{1.10}. We shall write this problem in a
slightly more general form, as follows. For $T>0$ given, we seek
a function $u:[0,T]\times \overline{D}\rightarrow{\mathbb C}$
satisfying
\begin{equation}\notag{({\mathcal N})}
\begin{split}
&u_t=\cunit\,a(t)\,u_{xx}+\cunit\,\beta(t,x)\,u+f(t,x)
\quad\forall\,(t,x)\in\,[0,T]\times\overline{D},\\
&u(t,0)=0\quad\forall\,t\in\,[0,T],\\
&u_x(t,1)=\mu(t)\,\left[\,S(t)\,u_t(t,1)+G(t)\,u(t,1)\,\right]
\quad\forall\,t\in\,[0,T],\\
&u(0,x)=u_0(x)\quad\forall\,x\in\,\overline{D}.\\
\end{split}
\end{equation}
We shall assume that $a:[0,T]\rightarrow{\mathbb R}\backslash\{0\}$,
$\beta$, $f:[0,T]\times\overline{D}\rightarrow{\mathbb C}$,
$u_0:\overline{D}\rightarrow \mathbb{C}$, $\mu$,
$S:[0,T]\rightarrow{\mathbb R}$ and $G:[0,T]\rightarrow{\mathbb C}$
are given functions.
We shall assume that the solution $u$ of ($\mathcal{N}$) exists
uniquely, and that the data and the solution of $(\mathcal{N})$
are smooth enough for the purposes of the error estimates that
will follow.
(In some numerical experiments of Section~\ref{sec55} we shall
revert to the specific physical data in \eqref{1.9},
\eqref{1.10}, \eqref{1.11}, and take the functions $a(t)$,
$\beta(t,x)$ as in \eqref{1.8}, $\mu(t)=\frac{\dot{s}(t)}{s(t)}$,
$S(t)=\frac{s^2(t)}{1+(\dot{s}(t))^2}$, $G(t)=g(t)\,S(t)+\cunit\,
[\,S(t)\,\dot{\delta}(t)-s^2(t)\,]$, where
$\delta=\frac{s\dot{s}}{2}\cdot$)
\subsubsection{Semidiscrete approximation}\label{N_Semi}
The weak formulation of $(\mathcal{N})$, obtained by taking the
$L^2(D)$ inner product of the p.d.e. in $(\mathcal{N})$ with a
function in ${\mathbb H}^1(D)$, integrating by parts and using the
boundary conditions, motivates defining $u_h:[0,T]\rightarrow
S_h$, the semidiscrete approximation of $u$, by the equation
\begin{equation}\label{2.34}
\begin{split}
(\partial_t u_h(t,\cdot),\phi)=&\,\cunit \,a(t)\,\mu(t)\,\big[\,
S(t)\,\partial_t u_h(t,1)+
G(t)\,u_h(t,1)\,\big]\,\overline{\phi(1)}\\
&-\cunit\,a(t)\,{\mathcal B}(u_h(t,\cdot),\phi)
+\cunit\,(\beta(t,\cdot)\,u_h(t,\cdot),\phi)
+(f(t,\cdot),\phi)\quad\forall\,\phi\in S_h,
\quad\forall\,t\in[0,T],\\
\end{split}
\end{equation}
and
\begin{equation}\label{semi_init}
u_h(0,\cdot)=R_hu_0(\cdot).
\end{equation}
%
%
%
\begin{proposition}\label{2.3.3**}
The problem \eqref{2.34}-\eqref{semi_init}
admits a unique solution $u_h\in C^1([0,T];S_h)$.
\end{proposition}
%
%
%
\begin{proof}
Let $\dim(S_h)=J$ and $\{\phi_j\}_{j=1}^{\ssy J}$ be a basis of
$S_h$ consisting of real-valued functions. Hence, we have
$R_hu_0=\sum_{j=1}^{\ssy J}\gamma^0_j\,\phi_j$ and
$u_h(t,x)=\sum_{j=1}^{\ssy J}\gamma_j(t)\,\phi_j(x)$, where
$\gamma_j:[0,T]\rightarrow{\mathbb C}$ for $j=1,\dots,J$. Then
\eqref{2.34}-\eqref{semi_init} is equivalent to the following
o.d.e. initial-value problem: Find ${\widetilde G}\in
C^1([0,T];{\mathbb C}^{\ssy J})$ such that ${\widetilde
G}(0)={\widetilde G}^{0}$ and
${\widetilde A}(t)\,{\widetilde G}'(t)={\widetilde
B}(t)\,{\widetilde G}(t)+{\widetilde F}(t)$, $\forall\,t\in[0,T],$
where ${\widetilde G}(t):=(\gamma_1(t),\dots,\gamma_{\ssy
J}(t))^T$, ${\widetilde G}^{0}:=(\gamma^0_{1},\dots,\gamma^0_{\ssy
J})^T$, ${\widetilde A}:[0,T]\rightarrow{\mathbb C}^{\ssy J\times
J}$ with ${\widetilde A}_{\ell j}(t):=(\phi_{\ell},\phi_j)
-\cunit\,a(t)\,S(t)\,\mu(t)\,\phi_{\ell}(1)\,\phi_j(1)$,
${\widetilde B}:[0,T]\rightarrow{\mathbb C}^{\ssy J\times J}$ with
${\widetilde B}_{\ell
j}(t):=\cunit\,a(t)\,G(t)\,\mu(t)\,\phi_{\ell}(1)\,\phi_j(1)
-\cunit\,a(t)\,{\mathcal B}(\phi_{\ell},\phi_j)
+\cunit\,(\beta(t,\cdot)\,\phi_{\ell},\phi_j)$, and
$\widetilde{F}:[0,T]\rightarrow{\mathbb C}^{\ssy J}$ with
$\widetilde{F}(t):=((f(t,\cdot),\phi_1),\dots,(f(t,\cdot),\phi_{\ssy
J}))^T$.
Since ${\widetilde A}$, ${\widetilde B}$, ${\widetilde F}$ are
continuous maps, to ensure existence and uniqueness of the
solution ${\widetilde G}$, it is sufficient to show that
${\widetilde A}(t)$ is nonsingular for $t\in[0,T]$. Indeed,
letting $t\in[0,T]$ and $x\in\text{\rm Ker}({\widetilde A}(t))$,
we have $\text{\rm Re}(\overline{x}^T{\widetilde A}(t)x)=0$, from
which we conclude that $\big\|\sum_{j=1}^{\ssy
J}x_j\phi_j\big\|^2=0$ and hence $x=0$. $\Box$
\end{proof}
%
%
%
\par
Let us first present a $H^1$ superconvergence error estimate for
the semidiscrete approximation $u_h$.
%
%
\begin{proposition}\label{H1SUP}
Let $u$ be the solution of $(\mathcal{N})$ and $u_h$ its
semidiscrete approximation defined by
\eqref{2.34}-\eqref{semi_init}. Assume that $\mu(t)\leq0$ and
$S(t)>0$ for $t\in[0,T]$. Then
\begin{equation}\label{H1superconv}
\|u_h(t,\cdot)-R_hu(t,\cdot)\|_1\leq\,C\,h^{r+1}
\,\left(\int_0^t\Gamma_{\ssy{\mathcal N}}(\tau)\;d\tau
\right)^{\frac{1}{2}}
\quad\forall\,t\in[0,T],\quad\forall\,h\in(0,h_{\star}],
\end{equation}
where
$\Gamma_{\ssy{\mathcal N}}(\tau):=\|u(\tau,\cdot)\|^2_{r+1}
+\|\partial_tu(\tau,\cdot)\|^2_{r+1}
+\sum_{\ell=0}^2\int_0^{\tau}
\|\partial_t^{\ell}u(s,\cdot)\|^2_{r+1}\,ds$.
\end{proposition}
%
%
%
%
%
%
%
\begin{proof}
Let $h\in(0,h_{\star}]$, $\theta_h:=u_h-R_hu$ and
$\xi(t):=\tfrac{1}{a(t)}$. Using \eqref{R_def} and \eqref{4.40}
we obtain
\begin{equation}\label{2.42magna}
\begin{split}
(\partial_t\theta_h(t,\cdot),\phi)&=\cunit\,a(t)\,\mu(t)\,\big[
S(t)\,\partial_t\theta_h(t,1)+G(t)\,\theta_h(t,1)\,\big]
\,\overline{\phi(1)}\\
&\quad-\cunit\,a(t)\,{\mathcal B}(\theta_h(t,\cdot),\phi)
+\cunit\,(P_h(\beta(t,\cdot)\,\theta_h(t,\cdot)),\phi)\\
&\quad+\big(\Psi_{\star}(t,\cdot),\phi\big)\quad\forall\,\phi\in
S_h,\quad\forall\,t\in[0,T],\\
\end{split}
\end{equation}
where $\Psi_{\star}:=[\partial_tu-R_h(\partial_tu)]
-\cunit\,\beta\,(u-R_hu)$.
Set $\phi=\partial_t\theta_h$ in \eqref{2.42magna} and then take
imaginary parts to obtain
\begin{equation}\label{2.43magna}
\begin{split}
\tfrac{d}{dt}|\theta_h(t,\cdot)|_1^2&\leq\,|\mu(t)|\big[-2\,
S^{\star}\,|\partial_t\theta_h(t,1)|^2+2\,|G(t)|\,|\theta_h(t,1)|
\,|\partial_t\theta_h(t,1)|\,\big]\\
&\quad +2\,|\xi(t)|\,\,\|\partial_t\theta_h(t,\cdot)\|_{-1,h}\,\,
|P_h(\beta(t,\cdot)\,\theta_h(t,\cdot))|_1\\
&\quad+2\,\xi(t)\,\text{\rm Im}(\Psi_{\star}(t,\cdot),
\partial_t\theta_h(t,\cdot))\quad\forall\,t\in[0,T],\\
\end{split}
\end{equation}
where $S^{\star}:=\inf_{[0,T]}S>0$.
In order to bound properly the quantity
$\|\partial_t\theta_h\|_{-1,h}$, first use \eqref{4.39} to obtain
\begin{equation}\label{auto_1}
\|\Psi_{\star}(t,\cdot)\|\leq\,C\,h^{r+1}
\,\big[\,\|u(t,\cdot)\|_{r+1}+
\|\partial_tu(t,\cdot)\|_{r+1}\,\big]\quad\forall\,t\in[0,T].
\end{equation}
Then, use of \eqref{SobolevIneq} and \eqref{auto_1} in
\eqref{2.42magna} gives
\begin{equation*}
\begin{split}
\big|(\partial_t\theta_h(t,\cdot),\phi)\big|\leq&\,|a(t)|\,\Big[
S(t)\,|\mu(t)|\,|\partial_t\theta_h(t,1)|
+\big(\,|G(t)|\,|\mu(t)|+1\,\big)
\,|\theta_h(t,\cdot)|_1\,\Big]\,|\phi|_1\\
&\quad
+|\beta(t,\cdot)|_{\infty}\,\|\theta_h(t,\cdot)\|\,\|\phi\|\\
&\quad+C\,h^{r+1}\big(\,\|\partial_tu(t,\cdot)\|_{r+1}
+\|u(t,\cdot)\|_{r+1}\,\big)\,\|\phi\|
\quad\forall\,\phi\in S_h,\quad\forall\,t\in[0,T],\\
\end{split}
\end{equation*}
which, along with \eqref{PoincareF}, yields that
\begin{equation}\label{auto_2}
\begin{split}
2\,|\xi(t)|\,\|\partial_t\theta_h(t,\cdot)\|_{-1,h}&
\leq\,2\,S(t)\,|\mu(t)|\,|\partial_t\theta_h(t,1)|\\
&\quad+C\,\Big[\,|\theta_h(t,\cdot)|_1+h^{r+1}\big(\,
\|\partial_tu(t,\cdot)\|_{r+1} +\|u(t,\cdot)\|_{r+1}\,\big)\,\Big]
\quad\forall\,t\in[0,T].\\
\end{split}
\end{equation}
Thus, combining \eqref{2.43magna}, \eqref{auto_2},
\eqref{PoincareF}, \eqref{SobolevIneq}, and \eqref{TheBound1}, we
arrive at
\begin{equation*}
\begin{split}
\tfrac{d}{dt}|\theta_h|_1^2\leq\,C\,\,\big[\,|\theta_h|_1^2
+h^{2(r+1)}\,(\, \|\partial_tu\|_{r+1}^2 +\|u\|_{r+1}^2\,)\,\big]
+2\,\xi\,\,\text{\rm Im}(\Psi_{\star},
\partial_t\theta_h)\quad\text{\rm on}\ \ [0,T].\\
\end{split}
\end{equation*}
Since $\theta_h(0,\cdot)=0$, integrating with respect to $t$ in
the inequality above yields
\begin{equation*}
\begin{split}
|\theta_h(t,\cdot)|_1^2\leq&\,C\,\left[\,
\int_0^t\,|\theta_h(s,\cdot)|_1^2\,ds
+h^{2(r+1)}\int_0^t\big(\,\|\partial_tu(s,\cdot)\|^2_{r+1}
+\|u(s,\cdot)\|_{r+1}^2\,\big)\,ds\,\right]\\
&+{\rm
Im}\Big{\{}\,2\,\xi(t)\,(\Psi_{\star}(t,\cdot),\theta_h(t,\cdot))
-2\,\int_0^t\xi'(s)\,(\Psi_{\star}(s,\cdot),\theta_h(s,\cdot))\,ds\\
&-2\,\int_0^t\xi(s)\,(
\partial_t\Psi_{\star}(s,\cdot),\theta_h(s,\cdot))\,ds\Big{\}}
\quad\forall\,t\in[0,T].\\
\end{split}
\end{equation*}
Using in the above the Cauchy-Schwarz inequality,
\eqref{PoincareF}, and \eqref{auto_1}, we obtain
\begin{equation}\label{auto_horse}
|\theta_h(t,\cdot)|_1^2\leq\,C\,
\int_0^t\,|\theta_h(s,\cdot)|_1^2\,ds
+{C}\,h^{2(r+1)}\,\Gamma_{\ssy{\mathcal N}}(t)
\quad\forall\,t\in[0,T].
\end{equation}
The estimate \eqref{H1superconv} follows from \eqref{auto_horse}
using Gr{\"o}nwall's lemma and \eqref{PoincareF}. $\Box$
\end{proof}
%
%
%
%
%
%
%
%
A simple consequence of this superconvergence estimate and the
approximation property \eqref{4.39} of the elliptic projection is
the following convergence result:
%
%
%
%
%
\begin{thm}\label{TheoremConv1}
Let $u$ be the solution of $(\mathcal{N})$ and $u_h$ its
semidiscrete approximation defined by
\eqref{2.34}-\eqref{semi_init}. Assume that $\mu(t)\leq0$ and
$S(t)>0$ for $t\in[0,T]$. Then
\begin{equation}\label{L2H1conv}
\|u_h(t,\cdot)-u(t,\cdot)\| +h\,\|u_h(t,\cdot)-u(t,\cdot)\|_1
\leq\,C\,h^{r+1}\,\left(\|u(t,\cdot)\|_{r+1}^2+
\int_0^t\Gamma_{\ssy{\mathcal N}}(\tau)\;d\tau
\,\right)^{\frac{1}{2}}\quad\forall\,t\in[0,T],
\end{equation}
where $\Gamma_{\ssy{\mathcal N}}$ is the function defined
in the statement of Proposition~\ref{H1SUP}. $\Box$
\end{thm}
%
%
%
%
%
%
%
%
\par
Therefore, taking into account the relation of $a$, $\mu$ and $S$
to the function $s(t)$ describing the bottom topography, we
conclude that the error estimate of Theorem~\ref{TheoremConv1}
holds in the case of domains with upsloping bottom profiles,
i.e., when $\dot{s}(t)\leq0$ for $t\in[0,T]$.
%
%
%
%
%
%
%
\begin{remark}
The $H^1$ superconvergence estimate
\eqref{H1superconv}, \eqref{SobolevIneq}, and a standard
$L^{\infty}$ estimate for the error of the elliptic projection
$($\cite{Wheeler1973}$)$ yield as usual an optimal-order
estimate of the error $|u-u_h|_{\infty}$ on $[0,T]$
$($cf. \cite{ref55}$)$.
\end{remark}
%
%
%
%
%
%
%
%
%
%
%
\subsubsection{Crank-Nicolson fully discrete approximations}\label{N_fully}
%
%
Let $N\in{\mathbb N}$ and $(t^n)_{n=0}^{\ssy N}$ be the nodes of
the partition of $[0,T]$ where, $t^0=0$, $t^{\ssy N}=T$ and
$t^n<t^{n+1}$ for $n=0,\dots,N-1$. Define $k_n:=t^n-t^{n-1}$ for
$n=1,\dots,N$, $t^{n+\frac{1}{2}}:=\frac{t^n+t^{n+1}}{2}$ for
$n=0,\dots,N-1$, and $k:=\max_{1\leq{n}\leq{\ssy N}}k_n$. We set
$u^n:=u(t^n,\cdot)$ for $n=0,\dots,N$, where $u$ is the solution
of $(\mathcal{N})$. Finally, for sequences $(V^m)_{m=0}^{\ssy M}$,
we define $\partial V^m:=\tfrac{1}{k_n}(V^m-V^{m-1})$ and
${\mathcal A}V^m=\tfrac{1}{2}(V^m+V^{m-1})$ for $m=1,\dots,M$.
\par
For $n=0,\dots,N$, the Crank-Nicolson method yields an
approximation $U_h^n\in S_h$ of $u(t^n,\cdot)$ as follows:
\par
{\tt Step 1}: Set
\begin{equation}\label{FDmethod1}
U^0_h:=R_hu_0.
\end{equation}
\par
{\tt Step 2}: For $n=1,\dots,N$, find $U^n_h\in S_h$ such that
\begin{equation}\label{FDmethod2}
\begin{split}
(\partial U_h^n,\chi)=&\,\cunit\,a^{n-\frac{1}{2}}\,
\mu^{n-\frac{1}{2}}\, \big[\,S^{n- \frac{1}{2}}\,\partial
U_h^n(1)+G^{n-\frac{1}{2}}\,{\mathcal A}U_h^n(1)\,\big]
\,\overline{\chi(1)}\\
&-\cunit\,a^{n-\frac{1}{2}}\,{\mathcal B}\big({\mathcal
A}U_h^{n},\chi\big) +\cunit\,\big(\beta^{n-\frac{1}{2}}
\,{\mathcal A}U_h^{n},\chi\big)
+\big(f^{n-\frac{1}{2}},\chi\big)\quad\forall\,\chi\in S_h,\\
\end{split}
\end{equation}
where $S^{n-\frac{1}{2}}:=S(t^{n-\frac{1}{2}})$,
$\mu^{n-\frac{1}{2}}:=\mu(t^{n-\frac{1}{2}})$,
$a^{n-\frac{1}{2}}:=a(t^{n-\frac{1}{2}})$,
$G^{n-\frac{1}{2}}:=G(t^{n-\frac{1}{2}})$,
$f^{n-\frac{1}{2}}:=f(t^{n-\frac{1}{2}},\cdot)$ and
$\beta^{n-\frac{1}{2}}:=\beta(t^{n-\frac{1}{2}},\cdot)$.
%
%
\par
We first examine the problem of existence and uniqueness of the
fully discrete approximation $U_h^n$.
%
%
%
%
%
%
%
%
\begin{proposition}\label{FDExistence}
Let $n\in\{1,\dots,N\}$ and suppose that $U_h^{n-1}\in S_h$ is
well defined. If $S^{n-\frac{1}{2}}>0$ and
$\mu^{n-\frac{1}{2}}\leq0$, then, there exists a constant $C_n$
such that if $k_n<C_n$, then $U_h^n$ is well defined by
\eqref{FDmethod2}.
\end{proposition}
%
%
%
%
%
\begin{proof}
Since \eqref{FDmethod2} is equivalent to a linear system of
algebraic equations with unknowns the coefficients of $U_h^{n}$
with respect to a basis of $S_h$, existence and uniqueness of
$U_h^{n}$ will follow if we show that if there is a $V\in S_h$
such that
\begin{equation}\label{uniq_up1}
\begin{split}
\tfrac{1}{k_n}\,(V,\phi)=&\,\cunit\,a^{n-\frac{1}{2}}
\,\mu^{n-\frac{1}{2}}\,\Big[\,
S^{n-\frac{1}{2}}\,\tfrac{1}{k_n}\,V(1)
+G^{n-\frac{1}{2}}\,\tfrac{1}{2}\,V(1)\,\Big]
\,\overline{\phi(1)}\\
&-\cunit\,\tfrac{a^{n-\frac{1}{2}}}{2}\,{\mathcal
B}\big(V,\phi\big)
+\tfrac{\cunit}{2}\,(P_h\big(\beta^{n-\frac{1}{2}}
\,V),\phi\big)\quad\forall\,\phi\in S_h,\\
\end{split}
\end{equation}
then $V=0$. Set $\phi=\tfrac{1}{k_n}V$ in \eqref{uniq_up1}, and
then take imaginary parts and use the arithmetic-geometric mean
inequality and \eqref{TheBound1} to obtain
\begin{equation}\label{uniq_up2}
\begin{split}
|V|_1^2=&\,\mu^{n-\frac{1}{2}}\,k_n\,\left[\,2\,S^{n-\frac{1}{2}}
\,\Big|\tfrac{V(1)}{k_n}\Big|^2 +\tfrac{1}{k_n}\,\text{\rm
Re}(G^{n-\frac{1}{2}})\,|V(1)|^2\,\right]
+\tfrac{k_n}{a^{n-\frac{1}{2}}}\,\text{\rm
Re}(P_h(\beta^{n-\frac{1}{2}}\,V),\tfrac{V}{k_n})\\
\leq&\,|\mu^{n-\frac{1}{2}}|\,k_n\,\left[\,-S^{n-\frac{1}{2}}
\,\Big|\tfrac{V(1)}{k_n}\Big|^2
+\tfrac{|G^{n-\frac{1}{2}}|^2}{4S^{n-\frac{1}{2}}}\,|V(1)|^2\right]
+\tfrac{C\,|\beta^{n-\frac{1}{2}}|_{1,\infty}}
{|a^{n-\frac{1}{2}}|}
\,k_n\,|V|_1\,\|\tfrac{V}{k_n}\|_{-1,h}.\\
\end{split}
\end{equation}
For $\phi\in S_h$, we use \eqref{uniq_up1}, \eqref{SobolevIneq}
and \eqref{PoincareF} to obtain
\begin{equation*}
\begin{split}
\big|(\tfrac{V}{k_n},\phi)\big|\leq&\,|a^{n-\frac{1}{2}}|
\,|\mu^{n-\frac{1}{2}}|\,\Big[\,
S^{n-\frac{1}{2}}\,\Big|\tfrac{V(1)}{k_n}\Big|
+|G^{n-\frac{1}{2}}|\,\tfrac{1}{2}\,|V(1)|\,\Big]
\,|\phi|_1\\
&+\tfrac{1}{2}\,\big[\,|a^{n-\frac{1}{2}}|
+C\,|\beta^{n-\frac{1}{2}}|_{\infty}\,\big]
\,|V|_1\,|\phi|_1,\\
\end{split}
\end{equation*}
which yields
\begin{equation}\label{minusEstim}
\|\tfrac{V}{k_n}\|_{-1,h}\leq\,|a^{n-\frac{1}{2}}|
\,|\mu^{n-\frac{1}{2}}|\,
S^{n-\frac{1}{2}}\,\Big|\tfrac{V(1)}{k_n}\Big|
+C_{{\ssy E}}\,|V|_1,
\end{equation}
where $C_{{\ssy E}}:=\tfrac{1}{2}\,\big[\,|a^{n-\frac{1}{2}}|
+C\,|\beta^{n-\frac{1}{2}}|_{\infty}
+|a^{n-\frac{1}{2}}|\,|\mu^{n-\frac{1}{2}}|
\,|G^{n-\frac{1}{2}}|\,\big]$.
Using \eqref{uniq_up2} and \eqref{minusEstim}, \eqref{SobolevIneq}
gives
\begin{equation}\label{uniq_up3}
|V|_1^2\,\left\{1-k_n\,\left[
\tfrac{|\mu^{n-\frac{1}{2}}|\,|G^{n-\frac{1}{2}}|^2}
{4S^{n-\frac{1}{2}}}+\tfrac{C\,C_{{\ssy
E}}\,|\beta^{n-\frac{1}{2}}|_{1,\infty}}
{|a^{n-\frac{1}{2}}|}
+\tfrac{C^2}{4}\,|\mu^{n-\frac{1}{2}}|
\,S^{n-\frac{1}{2}}\,|\beta^{n-\frac{1}{2}}|_{1,\infty}^2
\right]\right\}\leq0,
\end{equation}
which ends the proof. $\Box$
\end{proof}
%
%
%
%
%
%
%
\par
In particular, if we suppose that $\beta$ is in
$C([0,T],\,W^{1,\infty}(D))$, that $a,\;\mu,\;S,\;G$ are
continuous functions on $[0,T]$, and that $S(t)>0$ and
$\mu(t)\leq\,0$ for $t\in\,[0,T]$, (i.e. the upsloping case),
then the existence and uniqueness of the fully discrete
approximation $U_h^n$ follows if $k_n\leq\,C$, where $C$ is a
constant independent of $n$. This follows from
Proposition~\ref{FDExistence} and the fact that the quantity
multiplying $k_n$ in \eqref{uniq_up3} may be uniformly bounded
with respect to $n$.
\par
In the case of a general bottom topography we have:
%
%
%
\begin{proposition}
Let $n\in\{1,\dots,N\}$ and suppose that $U_h^{n-1}\in S_h$ is
well defined. Then, there exist constants $C_{n,1}$ and $C_{n,2}$
such that if $\tfrac{k_n}{h}<C_{n,1}$ and $k_n<C_{n,2}$, then
$U_h^n$ is well defined by \eqref{FDmethod2}.
\end{proposition}
%
%
%
%
%
%
%
\begin{proof}
Let $\dim S_h=J$ and $\{\phi_j\}_{j=1}^{\ssy J}$ be a basis of
$S_h$ consisting of real-valued functions. It is easily seen that
existence and uniqueness of $U_h^n$ is equivalent to the
invertibility of a matrix $\widetilde{M}\in{\mathbb C}^{\ssy
J\times J}$ defined by $\widetilde{M}_{\ell
j}:=\mathcal{M}(\phi_j,\phi_{\ell})$ for $j$, $\ell=1,\dots,J$,
where $\mathcal{M}:S_h\times S_h\rightarrow{\mathbb C}$ is given
by
$\mathcal{M}(\chi,\phi):=(\chi,\phi)
-\cunit\,a^{n-\frac{1}{2}}\,\mu^{n-\frac{1}{2}}
\,S^{n-\frac{1}{2}}\,\chi(1)\,{\overline{\phi(1)}}+
\tfrac{k_n}{2}\,\big[\, -\cunit\,(\beta^{n-\frac{1}{2}}\chi,\phi)
-\cunit\,\,\mu^{n-\frac{1}{2}}\,a^{n-\frac{1}{2}}
\,G^{n-\frac{1}{2}}\,\chi(1)\,{\overline{\phi(1)}}
+\cunit\,a^{n-\frac{1}{2}}\,{\mathcal B}(\chi,\phi) \,\big]$ for
$\chi$, $\phi\in S_h$. If $x\in \text{\rm Ker}(\widetilde{M})$ we
have ${\rm Re}\big[\mathcal{M}(\phi_{\star},\phi_{\star})\big]=0$
with $\phi_{\star}:=\sum_{j=1}^{\ssy J}x_j\,\phi_j$.
Then, using \eqref{TraceIneq} and \eqref{2.36}, we get
\begin{equation*}
\begin{split}
\|\phi_{\star}\|^2\leq&\,\tfrac{k_n}{2}\,
\Big[\,|\beta^{n-\frac{1}{2}}|_{\infty}
\,\|\phi_{\star}\|^2+2|\mu^{n-\frac{1}{2}}|
\,|a^{n-\frac{1}{2}}|\,|G^{n-\frac{1}{2}}|
\,\|\phi_{\star}\|\,|\phi_{\star}|_1\,\Big]\\
\leq&\,\tfrac{k_n}{2}\,\Big[
\,|\beta^{n-\frac{1}{2}}|_{\infty}
\,+\tfrac{C}{h}\,|\mu^{n-\frac{1}{2}}|
\,|a^{n-\frac{1}{2}}|\,|G^{n-\frac{1}{2}}|
\,\Big]\,\|\phi_{\star}\|^2,
\end{split}
\end{equation*}
which, under our hypotheses, yields $x=0$ and ends the proof.
$\Box$
\end{proof}
%
%
%
%
%
%
%
%
%
\par
Hence, if $\beta\in C([0,T],L^{\infty}(D))$ and $a$, $\mu$, $G$
are continuous on $[0,T]$ (i.e. in the case of general bottom
topography), the existence and uniqueness of $U_h$ follows if we
take $k_n\leq C_1$ and $\tfrac{k_n}{h}\leq C_2$, for some
constants $C_1$ and $C_2$ independent of $n$.
\par
We next establish the consistency of our fully discrete scheme in
the $t$ variable.
%
%
%
%
%
%
%
\begin{proposition}\label{2.8}
Let $u$ be the solution of $(\mathcal{N})$. For $n=1,\dots,N$,
define $\sigma^n:\overline{D}\rightarrow{\mathbb C}$ by
\begin{equation}\label{ConsistencyEquations}
\tfrac{u^{n}-u^{n-1}}{k_n}=\cunit\,a^{n-\frac{1}{2}}\,
u_{xx}(t^{n-\frac{1}{2}},\cdot)
+\cunit\,\beta^{n-\frac{1}{2}}\,{\mathcal A}u^n
+f^{n-\frac{1}{2}}+\sigma^n.
\end{equation}
Then,
\begin{gather}
\|\sigma^n\|\leq\,C\,(k_n)^2\,\,B^n_1(u),
\quad n=1,\dots,N,\label{Synepeia1}\\
\|\sigma^{n+1}-\sigma^n\|\leq\,C\,\big[\,(k_n)^2
+|k_{n+1}-k_n|\,\big]\,(k_n+k_{n+1})\,\,B^n_2(u),\quad
n=1,\dots,N-1, \label{Synepeia2}
\end{gather}
where
$B^n_1(u):=\sum_{\ell=2}^3\max_{\ssy[t^{n-1},t^n]}
\|\partial_t^{\ell}u\|$ and
$B^n_2(u):=\sum_{\ell=2}^4\max_{\ssy[t^{n-1},t^{n+1}]}
\|\partial_t^{\ell}u\|$.
\end{proposition}
%
%
%
%
%
%
%
%
\begin{proof}
It follows easily by using the partial differential equation and
Taylor's formula. $\Box$
\end{proof}
%
%
%
%
%
%
%
\par
We prove now that the following $H^1$ superconvergence estimate
holds in the fully discrete case.
%
%
%
%
%
%
\begin{proposition}\label{MagnaProposition}
Let $u$ be the solution of $(\mathcal{N})$ and
$(U_h^n)_{n=0}^{\ssy N}$ be the fully discrete approximations
that the method \eqref{FDmethod1}-\eqref{FDmethod2} produces.
Assume that $\mu(t)\leq0$ and $S(t)>0$ for $t\in[0,T]$. In
addition, assume that there exists a constant $C\ge 0$ such that
\begin{equation}\label{MeshCondition}
|k_{n+1}-k_n|\leq\,C\,\max\big\{k_n^2,k_{n+1}^2\big\}, \quad
n=1,\dots,N-1.
\end{equation}
Then, there exists a constant $C_1$ such that, if
$\displaystyle{\max_{1\leq{n}\leq{\ssy N}}}
(k_n\,C_1)\leq\tfrac{1}{3}$, there exists a constant $C>0$
such that
\begin{equation}\label{MagnaEst}
\max_{1\leq{n}\leq{\ssy N}}\|U_h^n-R_hu^n\|_1\leq\,C
\,(k^2+h^{r+1})\,\,\Xi_{\ssy{\mathcal N}}(u)
\quad\forall\,h\in(0,h_{\star}],
\end{equation}
where
$\Xi_{\ssy{\mathcal N}}(u):=\sum_{\ell=0}^2\max_{\ssy[0,T]}
\|\partial_t^{\ell}u\|_{r+1}+\max_{\ssy[0,T]}\|\partial_t^3u\|_1
+\max_{\ssy[0,T]}\|\partial_t^4u\|
+\max_{\ssy t\in[0,T]}|\partial_t^3u(t,1)|$.
\end{proposition}
%
%
%
%
%
%
%
\begin{proof}
Let $h\in(0,h_{\star}]$, $\theta_h^n:=U_h^n-R_hu^n$ for
$n=0,\dots,N$, $\xi:=\tfrac{1}{a}$ and
$\xi^{n-\frac{1}{2}}:=\xi(t^{n-\frac{1}{2}})$ for $n=1,\dots,N$.
We use \eqref{FDmethod2}, \eqref{ConsistencyEquations},
\eqref{R_def}, and \eqref{4.40}, to obtain
\begin{equation}\label{TOGOO1}
\begin{split}
(\partial\theta_h^n,\chi)=&\,
\cunit\,a^{n-\frac{1}{2}}\,\mu^{n-\frac{1}{2}}\,\Big[S^{n-\frac{1}{2}}\,
\partial\theta_h^n(1)
+G^{n-\frac{1}{2}}\,{\mathcal A}\theta_h^{n}(1)
-{\mathcal E}_3^n\,\Big]\,\overline{\chi(1)}\\
&-\cunit\,a^{n-\frac{1}{2}}\,{\mathcal B}({\mathcal
A}\theta_h^n,\chi)
+\cunit\,(P_h(\beta^{n-\frac{1}{2}}\,{\mathcal A}\theta_h^n),\chi)\\
&+({\mathcal
E}_1^n-\sigma^n,\chi)+\cunit\,a^{n-\frac{1}{2}}\,{\mathcal
B}({\mathcal E}_2^n,\chi)\quad\forall\,\chi\in S_h,
\quad n=1,\dots,N,\\
\end{split}
\end{equation}
where
\begin{equation*}
\begin{split}
{\mathcal E}_1^n&:=\partial u^n-R_h(\partial
u^n)-\cunit\,P_h[\,\beta^{n-\frac{1}{2}}({\mathcal
A}u^n-R_h({\mathcal A}u^n))\,\big],\\
{\mathcal E}_2^n&:=u(t^{n-\frac{1}{2}})-{\mathcal A}u^n,\\
{\mathcal E}_3^n&:=S^{n-\frac{1}{2}}\,\big[
\,\partial_tu(t^{n-\frac{1}{2}},1)-\partial u^n(1)\big]
+G^{n-\frac{1}{2}}\,\big[
\,u(t^{n-\frac{1}{2}},1)-{\mathcal A}u^n(1)\,\big].\\
\end{split}
\end{equation*}
Using Taylor's formula and \eqref{4.39}, we deduce the following
estimates:
\begin{equation}\label{FragmaE1}
\begin{split}
\|{\mathcal
E}_1^n\|&\leq\,C\,\tfrac{h^{r+1}}{k_n}\,\Big\|\int_{t^{n-1}}^{t^n}
\partial_tu(s,\cdot)\,ds\Big\|_{r+1}
+C\,|\beta^{n-\frac{1}{2}}|_{\infty}\,h^{r+1}\,\|{\mathcal
A}u^n\|_{r+1}\\
&\leq\,C\,h^{r+1}\,\left(\,\max_{\ssy[t^{n-1},t^n]}\|u\|_{r+1}
+\max_{\ssy[t^{n-1},t^n]}\|\partial_tu\|_{r+1}\,\right),\\
\end{split}
\end{equation}
\begin{equation}\label{FragmaE3}
|{\mathcal
E}_2^n|_1\leq\,C\,k_n^2\,\max_{\ssy[t^{n-1},t^n]}|\partial_t^2u|_1,
\end{equation}
and
\begin{equation}\label{FragmaE4}
|{\mathcal
E}_3^n|\leq\,C\,k_n^2\,\left[
\,\max_{t\in[t^{n-1},t^n]}|\partial_t^2u(t,1)|
+\max_{t\in[t^{n-1},t^n]}|\partial_t^3u(t,1)|\,\right]
\end{equation}
for $n=1,\dots,N$.
Set $\chi=\partial\theta_h^n$ in \eqref{TOGOO1}, and then take
imaginary parts to obtain
\begin{equation}\label{AstroEst1}
\begin{split}
|\theta_h^n|_1^2\leq&\,|\theta_h^{n-1}|_1^2
+2\,k_n\,|\xi^{n-\frac{1}{2}}|\,
\,|P_h(\beta^{n-\frac{1}{2}}\,{\mathcal A}\theta_h^n)|_1
\,\,\|\partial\theta_h^n\|_{-1,h}\\
&+k_n\,|\mu^{n-\frac{1}{2}}|\,\big[\,-2\,S_{\star}\,
|\partial\theta_h^n(1)|^2 +2\,|G^{n-\frac{1}{2}}|
\,|{\mathcal A}\theta_h^{n}(1)|\,|\partial\theta_h^n(1)|
+2\,|{\mathcal E}_3^n|\,|\partial\theta_h^n(1)|\,\big]\\
&+2\,k_n\,{\rm Re}[{\mathcal B}({\mathcal
E}_2^n,\partial\theta_h^n)]+2\,k_n\,\xi^{n-\frac{1}{2}}
\,{\rm Im}({\mathcal E}_1^n-\sigma^n,\partial\theta_h^n),
\quad n=1,\dots,N,\\
\end{split}
\end{equation}
where $S_{\star}:=\inf_{[0,T]}S$.
\par
Now let us estimate $\|\partial\theta_h^n\|_{-1,h}$. For
$\varphi\in S_h$, \eqref{TOGOO1}-\eqref{FragmaE4},
\eqref{Synepeia1}, \eqref{SobolevIneq}, and \eqref{PoincareF} give
\begin{equation*}
\begin{split}
|(\partial\theta_h^n,\varphi)|\leq&\,
|a^{n-\frac{1}{2}}|\,|\mu^{n-\frac{1}{2}}|\,S^{n-\frac{1}{2}}
\,|\partial\theta_h^n(1)|\,|\varphi|_1\\
&+C\,|{\mathcal A}\theta_h^n|_1\,|\varphi|_1
+C\,(h^{r+1}+k_n^2)\,\,|\varphi|_1\,\,\Xi_1(u),
\quad n=1,\dots,N,\\
\end{split}
\end{equation*}
where
$\Xi_1(u):=\max_{\ssy[0,T]}\|u\|_{r+1}
+\max_{\ssy[0,T]}\|\partial_tu\|_{r+1}
+\max_{\ssy[0,T]}\|\partial_t^2u\|_1
+\max_{\ssy[0,T]}\|\partial_t^3u\|
+\max_{\ssy t\in[0,T]}|\partial_t^3u(t,1)|$.
Hence, we conclude that
\begin{equation}\label{TOGOOO1}
\begin{split}
2\,k_n\,|\xi^{n-\frac{1}{2}}|\,\|\partial\theta_h^n\|_{-1,h}\leq&\,
2\,k_n\,|\mu^{n-\frac{1}{2}}|\,S^{n-\frac{1}{2}}\,
|\partial\theta_h^n(1)|\\
&+C\,k_n\,|{\mathcal A}\theta_h^n|_1
+C\,k_n\,(h^{r+1}+k_n^2)\,\,\Xi_1(u),
\quad n=1,\dots,N.\\
\end{split}
\end{equation}
\par
Now, combining \eqref{TOGOOO1} and \eqref{AstroEst1} we have
\begin{equation*}
\begin{split}
|\theta_h^n|_1^2\leq&\,|\theta_h^{n-1}|_1^2 +C\,k_n\,|{\mathcal
A}\theta_h^n|_1^2
+C\,k_n\,\big[\,(k_n)^4+(h^{r+1}+k_n^2)\,|{\mathcal
A}\theta_h^n|_1\,\big]\,\Xi_1(u)\\
&+2\,k_n\,{\rm Re}[{\mathcal B}({\mathcal
E}_2^n,\partial\theta_h^n)]+2\,k_n\,\xi^{n-\frac{1}{2}}\,{\rm
Im}({\mathcal E}_1^n-\sigma^n,\partial\theta_h^n),\quad
n=1,\dots,N,
\end{split}
\end{equation*}
from which there follows that for some constant $C_1\geq 0$
\begin{equation}\label{AstroEst3}
\begin{split}
(1-C_1\,k_n)\,|\theta_h^n|_1^2\leq&\,(1+C_1\,k_n)
\,|\theta_h^{n-1}|_1^2
+C_2\,k_n\,(h^{r+1}+k_n^2)^2\,(\Xi_1(u))^2\\
&+2\,k_n\,{\rm Re}[{\mathcal B}({\mathcal
E}_2^n,\partial\theta_h^n)]+2\,k_n\,\xi^{n-\frac{1}{2}}\,{\rm
Im}({\mathcal E}_1^n-\sigma^n,\partial\theta_h^n),
\quad n=1,\dots,N.\\
\end{split}
\end{equation}
To continue, we assume that
$\displaystyle{\max_{1\leq{n}\leq{\ssy N}}}(C_1\,k_n)
\leq\,\tfrac{1}{3}$,
which allows us to conclude that
$\tfrac{1+C_1\,k_n}{1-C_1\,k_n}\leq\,e^{3C_1k_n}$ for
$n=1,\dots,N$. Hence, \eqref{AstroEst3} yields
\begin{equation*}
\begin{split}
|\theta_h^n|_1^2\leq&\,e^{3C_1k_n}\,|\theta_h^{n-1}|_1^2
+\tfrac{C_2\,k_n}{1-C_1\,k_n}\,(h^{r+1}+k_n^2)^2\,(\Xi_1(u))^2\\
&+\tfrac{2\,k_n}{1-C_1\,k_n}\,\Big[\,{\rm Re}[{\mathcal
B}({\mathcal E}_2^n,\partial\theta_h^n)]+\xi^{n-\frac{1}{2}}\,{\rm
Im}({\mathcal E}_1^n-\sigma^n,\partial\theta_h^n)\,\Big],
\quad n=1,\dots,N.\\
\end{split}
\end{equation*}
Next, we define
$\lambda_j^n:=\tfrac{\exp\left(3C_1\sum_{\ell=j+1}^{n}
k_{\ell}\right)}{1-C_1\,k_j}$
and use a simple induction argument to arrive at
\begin{equation*}
\begin{split}
|\theta_h^n|_1^2\leq&\,C_2\,(\Xi_1(u))^2\,\sum_{j=1}^{n}k_j
\,\lambda_j^n\,(h^{r+1}+k_j^2)^2\\
&+2\,\sum_{j=1}^{n}k_j \,\lambda_j^n\,\Big[\,{\rm Re}[{\mathcal
B}({\mathcal E}_2^j,\partial\theta_h^j)]
+\xi^{j-\frac{1}{2}}\,{\rm Im}({\mathcal
E}_1^j-\sigma^j,\partial\theta_h^j)\,\Big],\quad n=1,\dots,N,\\
\end{split}
\end{equation*}
which yields
\begin{equation}\label{AstroEst6}
|\theta_h^n|_1^2\leq\,C\,(h^{r+1}+k^2)^2\,(\Xi_1(u))^2+T_{\ssy
A}^n+T_{\ssy B}^n,\quad n=1,\dots,N,
\end{equation}
where
\begin{equation*}
\begin{split}
T_{\ssy A}^n&:=2\,\sum_{j=1}^{n} \lambda_j^n\,{\rm
Re}\big[\,{\mathcal B}({\mathcal
E}_2^j,\theta_h^j-\theta_h^{j-1})\,\big],\\
T_{\ssy B}^n&:=2\,\sum_{j=1}^{n} \lambda_j^n\,
\xi^{j-\frac{1}{2}}\,{\rm Im}({\mathcal E}_1^j
-\sigma^j,\theta_h^j-\theta_h^{j-1}).\\
\end{split}
\end{equation*}
First we observe that
\begin{equation}\label{AstroMagma1}
\begin{split}
T_{\ssy A}^n&=\tfrac{2}{1-C_1\,k_n}\,\,{\rm Re}\big[\,{\mathcal
B}({\mathcal E}_2^n,\theta_h^n)\,\big]\\
&\quad+2\,\sum_{j=1}^{n-1} \lambda_j^n\,\,{\rm
Re}[{\mathcal B}({\mathcal E}_2^j-{\mathcal E}_2^{j+1},\theta_h^j)]\\
&\quad+2\sum_{j=1}^{n-1}\,\exp\Bigg(3\,C_1\sum_{\ell=j+2}^{n}k_{\ell}\Bigg)
\,\Big[\tfrac{\exp\left(3\,C_1\,k_{j+1}\right)-1+C_1\,k_j}{1-C_1\,k_j}
-\tfrac{C_1\,k_{j+1}}{1-C_1\,k_{j+1}}\Big]\,\,\,{\rm
Re}\big[\,{\mathcal B}({\mathcal E}^{j+1}_2,\theta_h^j)\,\big],\\
\end{split}
\end{equation}
for $n=1,\dots,N$.
Since
\begin{equation*}
|{\mathcal E}_2^j-{\mathcal
E}_2^{j+1}|_1\leq\,C\,(k_j+k_{j+1})\,\big[\,(k_j)^2+|k_{j+1}-k_j|\,\big]\,
\Xi_2(u),\quad j=1,\dots,N-1,
\end{equation*}
with
$\Xi_2(u):=\max_{\ssy[0,T]}|\partial_t^2u|_1
+\max_{\ssy[0,T]}|\partial_t^3u|_1$,
we see that \eqref{AstroMagma1}, \eqref{FragmaE3}, and
\eqref{MeshCondition} yield
\begin{equation}\label{TEstim1}
|T_{\ssy A}^n|\leq\,C\,k^2\,\Xi_2(u)
\,\max_{1\leq{m}\leq{n}}|\theta_h^m|_1,
\quad n=1,\dots,N.
\end{equation}
In addition, we have
\begin{equation}\label{AstroMagma2}
\begin{split}
T_{\ssy B}^n=&\tfrac{2}{1-C_1\,k_n}\, \xi^{n-\frac{1}{2}}
\,{\rm Im}({\mathcal E}_1^n-\sigma^n,\theta_h^n)\\
&+2\,\sum_{j=1}^{n-1} \lambda_j^n\, \xi^{j-\frac{1}{2}}
\,{\rm Im}({\mathcal E}_1^j-\sigma^j
-{\mathcal E}_1^{j+1}+\sigma^{j+1},\theta_h^j)\\
&+\,2\,\sum_{j=1}^{n-1}\xi^{j-\frac{1}{2}}\,\,
\exp\Bigg(3\,C_1\sum_{\ell=j+2}^{n}k_{\ell}\Bigg)\,\Big[
\tfrac{\exp\left(3\,C_1\,k_{j+1}\right)-1+C_1\,k_j}{1-C_1\,k_j}
-\tfrac{C_1\,k_{j+1}}{1-C_1\,k_{j+1}}\Big]\,\,{\rm Im}({\mathcal
E}_1^{j+1}-\sigma^{j+1},\theta_h^j)\\
&+\,2\,\sum_{j=1}^{n-1}(\xi^{j-\frac{1}{2}}-\xi^{j+\frac{1}{2}})\,
\lambda_{j+1}^n
\,\,{\rm Im}({\mathcal E}_1^{j+1}-\sigma^{j+1},\theta_h^j),
\quad n=1,\dots,N.\\
\end{split}
\end{equation}
Observing that
\begin{equation*}
|{\mathcal E}_1^j-{\mathcal E}_1^{j+1}|_1
\leq\,C\,(k_j+k_{j+1})\,h^{r+1}\,\Xi_3(u),\quad
j=1,\dots,N-1,
\end{equation*}
with
$\Xi_3(u):=\max_{\ssy[0,T]}\|\partial_tu\|_{r+1}
+\max_{\ssy[0,T]}\|\partial_t^2u\|_{r+1}$,
we see that \eqref{AstroMagma2}, \eqref{FragmaE1},
\eqref{Synepeia1}-\eqref{MeshCondition} and, \eqref{PoincareF}
yield
\begin{equation}\label{TEstim2}
|T_{\ssy B}^n|\leq\,C\,(k^2+h^{r+1})\,\Xi_4(u)
\,\max_{1\leq{m}\leq{n}}|\theta_h^m|_1, \quad n=1,\dots,N,
\end{equation}
where
$\Xi_4(u):=\sum_{\ell=0}^2\max_{\ssy[0,T]}
\|\partial_t^{\ell}u\|_{r+1}
+\sum_{\ell=3}^4\max_{\ssy[0,T]}
\|\partial_t^{\ell}u\|$.
Now, from \eqref{AstroEst6}, \eqref{TEstim1}, and \eqref{TEstim2}
there follows that
\begin{equation*}
|\theta_h^n|_1^2\leq\,C\,(h^{r+1}+k^2)^2\,(\Xi_1(u))^2
+C\,(k^2+h^{r+1})\,\big(\,\Xi_2(u)+\Xi_4(u)\,\big)
\,\max_{1\leq{m}\leq{n}}|\theta_h^m|_1,
\quad n=1,\dots,N,
\end{equation*}
which easily yields
\begin{equation}\label{AstroEst7}
\max_{0,\leq{n}\leq{\ssy
N}}|\theta_h^n|_1^2\leq\,C\,(h^{r+1}+k^2)^2
\,(\Xi_1(u)+\Xi_2(u)+\Xi_4(u))^2.
\end{equation}
The desired estimate \eqref{MagnaEst} is then a simple consequence
of \eqref{AstroEst7} and \eqref{PoincareF}. $\Box$
\end{proof}
%
%
%
%
%
%
%
\par
Now we are ready to prove error estimates in the $L^2$ and $H^1$
norms.
%
%
%
%
%
%
%
%
\begin{thm}\label{MagnaTheorema}
Let $u$ be the solution of $({\mathcal N})$ and
$(U_h^n)_{n=0}^{\ssy N}$ be the fully discrete approximations
that the method \eqref{FDmethod1}-\eqref{FDmethod2} produces.
Assume that $\mu(t)\leq0$, $S(t)>0$, for $t\in[0,T]$, that
\eqref{MeshCondition} holds and $\max_{1\leq{n}\leq{\ssy
N}}(C_1\,k_n)\leq\tfrac{1}{3}$, where $C_1$ is the constant
specified in Proposition~\ref{MagnaProposition}. Then
\begin{equation*}
\max_{0\leq{n}\leq{\ssy N}}\|U_h^n-u^n\|_{\ell}\leq
\,C\,\big(\,k^2+h^{r+1-\ell}\,\big)
\,\,\Xi_{\ssy{\mathcal N}}(u),\quad\forall\,h\in(0,h_{\star}],
\end{equation*}
for $\ell=0,1$, where $\Xi_{\ssy{\mathcal N}}(u)$ was specified
in Proposition~\ref{MagnaProposition}.
\end{thm}
%
%
%
%
%
%
%
%
%
\begin{proof}
It is a simple consequence of \eqref{MagnaEst} and \eqref{4.39}.
$\Box$
\end{proof}
%
%
%
%
%
\par
We conclude that in the case of upsloping bottoms, the fully
discrete Crank-Nicolson-Galerkin method
\eqref{FDmethod1}-\eqref{FDmethod2} yields fully discrete
approximations $U_h^n$ that converge to the solution $u$ of
$(\mathcal{N})$ at optimal rates in the $L^2$ and $H^1$ norms.
%
%
%
%
\subsection{The Abrahamsson-Kreiss boundary condition}\label{sec4}
We consider now the (PE) with the Abrahamsson-Kreiss bottom
boundary condition, i.e. the ibvp \eqref{1.6}, \eqref{1.8},
\eqref{1.9}, \eqref{1.12}, which we rewrite here, in slightly
more general form, for the convenience of the reader.
For $T>0$ given, seek a function $u:[0,T]\times
\overline{D}\rightarrow \mathbb{C}$ satisfying
\begin{equation}\notag{({\mathcal A}{\mathcal K})}
\begin{split}
&u_t=\cunit\,a(t)\,u_{xx}
+\cunit\,\beta(t,x)\,u+f(t,x)\quad\forall
(t,x)\in\,[0,T]\times\overline{D},\\
&u(t,0)=0\quad\forall t\in\,[0,T],\\
&u_x(t,1)=0\quad\forall t\in\,[0,T],\\
&u(0,x)=u_0(x)\quad\forall x\in\,\overline{D}.\\
\end{split}
\end{equation}
We assume again that $a:[0,T]\rightarrow{\mathbb
R}\backslash\{0\}$, $\beta$, $f:[0,T]\times
\overline{D}\rightarrow \mathbb{C}$, $u_0:\overline{D}\rightarrow
\mathbb{C}$ are given functions. We shall assume that the
solution of $(\mathcal{AK})$ exists uniquely and that the data
and the solution of $(\mathcal{AK})$ are smooth enough for the
purposes of the error estimation. We note that $(\mathcal{AK})$
may be considered as a special case of $(\mathcal{N})$ obtained
by setting $\mu$ equal to zero in $(\mathcal{N})$. (This does not
imply of course that we assume that $\dot{s}$ is zero. We recall
that in the Abrahamsson-Kreiss formulation the effect of variable
bottom enters explicitly in the definition of $a$ and $\beta$,
cf. \eqref{1.7}, and in the change-of-variable formula
\eqref{1.5}.) All the error estimates for $(\mathcal{AK})$ that
follow may then be considered as special cases of the analogous
estimates in the two preceding paragraphs but with some important
simplifications. For the convenience of the reader we shall
restate the results but not prove them in detail; we shall just
point out some differences between them and the analogous
estimates for the problem $(\mathcal{N})$. It will be seen that
the finite element approximations of $(\mathcal{AK})$ exist and
satisfy optimal-order error estimates under no further
assumptions (except smoothness) on the shape of the bottom.
\subsubsection{Semidiscrete approximation}
Using the finite element subspace $S_h$ and the notation
established in paragraph~\ref{subsec40}, we define the
semidiscrete approximation $u_h$ of the solution of
$(\mathcal{AK})$ as the map $u_h:[0,T]\rightarrow S_h$ satisfying
\begin{equation}\label{2.43}
(\partial_t u_h(t,\cdot),\phi)=-\cunit\,a(t)\,{\mathcal
B}(u_h(t,\cdot),\phi) +\cunit\,(\beta(t,\cdot)\,u_h(t,\cdot),\phi)
+(f(t,\cdot),\phi)\quad\forall\,\phi\in S_h,
\quad\forall\,t\in[0,T],\\
\end{equation}
and
\begin{equation}\label{semi_init2}
u_h(0,\cdot)=u_h^0,
\end{equation}
where $u_h^0\in S_h$ is an approximation of $u_0$, which may be taken,
for example, as $P_h u_0$ or $R_h u_0$.
%
%
%
%
\begin{proposition}\label{2.11}
The problem \eqref{2.43}-\eqref{semi_init2} admits a unique
solution in $C^1([0,T],S_h)$. If $f\equiv0$ and
$\beta_{I}\equiv0$, then the solution preserves the $L^2(D)$ norm,
i.e.,
\begin{equation}\label{2.45}
\|u_h(t,\cdot)\|=\|u_h^0\|\quad\forall\,t\in [0,T].
\end{equation}
\end{proposition}
%
%
%
%
%
%
%
\begin{proof}
The first part follows from Proposition~\ref{2.3.3**} for $\mu=0$.
The conservation of the $L^2(D)$ norm follows by taking $\phi=u_h$
in \eqref{2.43} and then real parts. $\Box$
\end{proof}
%
%
%
%
%
%
%
%
%
%
\begin{thm}\label{2.12}
Let $u$ be the solution of $(\mathcal{AK})$ and
$u_h$ its semidiscrete approximation defined by
\eqref{2.43}-\eqref{semi_init2}.
Then
\begin{equation}\label{2.46}
\|u(t,\cdot)-u_h(t,\cdot)\|_{\ell}\leq\,C\,\left[\,
\|u_h^0-R_hu_0\|_{\ell} +h^{r+1-\ell}\,\left(\|u(t,\cdot)\|_{r+1}^2
+\int_0^t\Gamma_{\ssy{{\mathcal A\mathcal K},\ell}}(\tau)
\,d\tau\right)^{\frac{1}{2}}
\right]\quad\forall\,t\in[0,T],
\end{equation}
for $\ell=0,1$ and $h\in(0,h_{\star}]$, where
$\Gamma_{\ssy{{\mathcal A\mathcal K},0}}(\tau):=\sum_{m=0}^1
\|\partial_t^mu(\tau,\cdot)\|_{r+1}^2$ and
$\Gamma_{\ssy{{\mathcal A\mathcal K},1}}(\tau)
:=\Gamma_{\ssy{\mathcal N}}(\tau)$, where $\Gamma_{\ssy{\mathcal
N}}$ is the function defined in the statement of
Proposition~\ref{H1SUP}.
\end{thm}
%
%
%
%
%
%
%
\begin{proof}
Let $h\in(0,h_{\star}]$.
Defining as usual $\theta_h:=u_h-R_h u$ we obtain
\begin{equation}\label{2.48}
(\partial_t\theta_h(t,\cdot),\phi)=-\cunit\,a(t)\,{\mathcal
B}(\theta_h(t,\cdot),\phi)
+\cunit\,(\beta(t,\cdot)\,\theta_h(t,\cdot),\phi)
+\big(\Psi_{\star}(t,\cdot),\phi\big)\quad\forall\,\phi\in
S_h,\quad\forall\,t\in[0,T],\\
\end{equation}
where $\Psi_{\star}$ was defined in the course of the proof of
Proposition \ref{H1SUP}. Taking $\phi=\theta_h$ in \eqref{2.48}
and then real parts, we may prove \eqref{2.46} with $\ell=0$
in a straightforward manner. The proof of \eqref{2.46} with $\ell=1$
follows the steps of the proof of Proposition~\ref{H1SUP}
if we take $\phi=\partial_t\theta_h$ in \eqref{2.48} and then
imaginary parts. $\Box$
\end{proof}
%
%
%
%
\begin{remark}
Hence, if $u_h^0$ is taken equal to $P_hu^0$ or $R_hu^0$,
Theorem~\ref{2.12} yields optimal-order estimates of the
error $u-u_h$ in the $L^2$ or $H^1$ norm, respectively.
Also, we note that to obtain the estimate \eqref{2.46} with
$\ell=0$, we do not need the inverse inequality \eqref{2.36}.
\end{remark}
%
%
%
\subsubsection{Crank-Nicolson fully discrete approximations}
We now proceed to the full discretization of $(\mathcal{AK})$ by
discretizing the initial-value problem
\eqref{2.43}-\eqref{semi_init2} in $t$ using the Crank-Nicolson
scheme. With notation introduced in paragraph~\ref{N_fully}, we
define for $n=0,\dots,N$ approximations $U_h^n\in S_h$ of
$u(t^n,\cdot)$, the solution of $(\mathcal{AK})$, as follows:
\par
{\tt Step 1}: Set
\begin{equation}\label{2.49}
U^0_h:=u_h^0.
\end{equation}
\par
{\tt Step 2}: For $n=1,\dots,N$, find $U^n_h\in S_h$ such that
\begin{equation}\label{2.50}
\begin{split}
(\partial U_h^n,\chi)=-\cunit\,a^{n-\frac{1}{2}}\,{\mathcal
B}\big({\mathcal A}U_h^{n},\chi\big)
+\cunit\,\big(\beta^{n-\frac{1}{2}}
\,{\mathcal A}U_h^{n},\chi\big)
+\big(f^{n-\frac{1}{2}},\chi\big)\quad\forall\,\chi\in S_h.
\end{split}
\end{equation}
%
%
%
%
%
\begin{proposition}\label{2.13}
Let $n\in\{1,\dots,N\}$ and suppose $U_h^{n-1}$ is well defined.
Then, there exists a constant $C$ independent of $n$ such that if
$k_n\leq C$, $U_h^n$ is well defined by \eqref{2.50}. Moreover, if
$f\equiv0$ and $\beta_I\equiv0$, then
\begin{equation}\label{2.51}
\|U_h^n\|=\|u_h^0\|,\quad n=0,\dots,N.
\end{equation}
\end{proposition}
%
%
%
%
%
%
\begin{proof}
Since \eqref{2.50} is equivalent to a $\dim(S_h)\times\dim(S_h)$
linear system of algebraic equations, existence and uniqueness of
$U_h^n$ will follow if we show that if there is a $V\in S_h$ such
that
\begin{equation}\label{2.52}
\tfrac{1}{k_n}(V,\phi)=
-\tfrac{\cunit}{2}\,a^{n-\frac{1}{2}}\,{\mathcal B}(V,\phi)
+ \tfrac{\cunit}{2}\,(\beta^{n-\frac{1}{2}}V,\phi)
\quad\forall \phi\in S_h,
\end{equation}
then $V=0$. This fact follows easily for $k_n$ sufficiently small, if we
put $\phi=V$ in \eqref{2.52} and take real parts. The conservation
property \eqref{2.51} follows from \eqref{2.50} if we select
$\chi={\mathcal A}U_h^{n}$ and take real parts. $\Box$
\end{proof}
%
%
%
%
%
%
%
%
%
%
\begin{thm}\label{2.14}
Let $u$ be the solution of $(\mathcal{AK})$ and
$(U_h^n)_{n=0}^{\ssy N}$ be the fully discrete approximations
produced by \eqref{2.49}-\eqref{2.50}. Then, if $\max_{1\leq
n\leq{\ssy N}}k_n$ is sufficiently small, we have
\begin{equation}\label{AK_FD_L2}
\max_{0\leq{n}\leq{\ssy N}}\|U_h^n-u^n\|
\leq\,C\,\left[\,\|u_h^0-R_hu_0\|
+(k^2+h^{r+1})\,\Xi_{\ssy{{\mathcal A\mathcal K},0}}(u)
\,\right]\quad\forall\,h\in(0,h_{\star}],
\end{equation}
where $\Xi_{\ssy{{\mathcal A\mathcal K},0}}(u)
:=\sum_{m=0}^1\max_{\ssy [0,T]}\|\partial_t^mu\|_{r+1}
+\sum_{m=2}^3\max_{\ssy[0,T]}\|\partial_t^mu\|
+\max_{\ssy[0,T]}\|\partial_t^2\partial_x^2u\|$. Also, if
\eqref{2.36} and \eqref{MeshCondition} hold, and $\max_{1\leq
n\leq{\ssy N}}k_n$ is sufficiently small, then
\begin{equation}\label{AK_FD_H1}
\max_{0\leq{n}\leq{\ssy N}}\|U_h^n-u^n\|_1
\leq\,C\,\left[\,\|u_h^0-R_hu_0\|_1
+(k^2+h^r)\,\Xi_{\ssy{{\mathcal A\mathcal K},1}}(u)\,\right]
\quad\forall\,h\in(0,h_{\star}],
\end{equation}
where
$\Xi_{\ssy{{\mathcal A\mathcal K},1}}(u):=
\Xi_{\ssy{{\mathcal A\mathcal K},0}}(u)
+\max_{\ssy[0,T]}\|\partial_t^2u\|_{r+1}
+\max_{\ssy[0,T]}\|\partial_t^4u\|
+\max_{\ssy[0,T]}\|\partial_t^3\partial_x^2u\|$.
\end{thm}
%
%
%
%
%
%
%
%
\begin{proof}
First, we modify the consistency argument of Proposition~\ref{2.8}
defining, for $n=1,\dots,N$, $\sigma^n:{\overline D}
\rightarrow{\mathbb C}$ by
$\tfrac{u^{n}-u^{n-1}}{k_n}=\cunit\,a^{n-\frac{1}{2}}\,
{\mathcal A}(u_{xx}(t^n,\cdot))
+\cunit\,\beta^{n-\frac{1}{2}}\,{\mathcal A}u^n
+f^{n-\frac{1}{2}}+\sigma^n$.
Then, we set $\theta_h^n:=U_h^n-R_hu^n$ for $n=0,\dots,N$, to
obtain \eqref{TOGOO1} simplified by setting
$\mu^{n-\frac{1}{2}}=0$ and ${\mathcal E}_2^n=0$. To obtain
\eqref{AK_FD_L2} we put $\chi=\theta_h^n$ and then take real
parts. To obtain \eqref{AK_FD_H1} we proceed along the lines of
the proof of Proposition~\ref{MagnaProposition} appropriately
simplified. $\Box$
\end{proof}
%
%
%
%
%
%
%
%
%
%
%
%
%
%
\section{Numerical experiments}\label{sec55}
In this section we present the results of some numerical
experiments that we performed using the fully discrete
Galerkin-finite element methods, defined and analyzed in the
previous section, to solve the ibvp for the (PE) in domains of
variable bottom topography with Neumann and Abrahamsson-Kreiss
boundary conditions. We also make, in paragraph 3.3, a
theoretical excursion with the aim of explaining some
experimental observations made in paragraph 3.2. Recall that in
the case of the Neumann boundary condition, i.e. for the problem
($\mathcal{N}$), our convergence results were rigorously
established in the case of upsloping bottoms, that is when
$\dot{s}(t)\leq 0$ for all $t\in[0,T]$. One of our goals in this
section is to study numerically the behavior of the Neumann
boundary condition in the presence of downsloping bottoms and
compare the solution of $(\mathcal{N})$ with that of
($\mathcal{AK}$), for which rigorous convergence results hold for
any smooth $s(t)$. In the numerical experiments the finite element
subspace $S_h$ consisted of continuous, piecewise linear
functions defined on a uniform mesh, while the temporal
discretization was effected with uniform time step. All
computations were performed using double precision fortran 77.
\subsection{Order of convergence}
To test numerically the order of convergence of the fully
discrete Crank-Nicolson-Galerkin finite element method (henceforth
referred to as (FE)) in the case of the ibvp ($\mathcal{N}$), we
took $T=1$ and considered three cases of bottom profiles, namely:
\begin{equation*}
\left.\begin{array}{l}
\mbox{Case 1: }s(t)=-0.3t+0.7\mbox{ (upsloping)}.\\
\mbox{Case 2: }s(t)=0.4t+0.3\mbox{ (downsloping)}.\\
\mbox{Case 3: }s(t)=0.2\cos(4\pi t)+0.2\sin(4\pi t)+0.7\mbox{
(oscillatory)}.
\end{array}\right.
\end{equation*}
In ($\mathcal{N}$) we took $a=1/(2s^2)$, $\beta(t,x)=xt+{\rm
i}(3x+t^2)$, $u_0(x)=-x(x-1)^3$. The bottom boundary condition
had the form
$u_x(t,1)=\mu(t){\big[}S(t)u_t(t,1)+G(t)u(t,1){\big]}+f_1(t)$,
where $\mu(t)=\frac{\dot{s}(t)}{s(t)}$,
$S(t)=\frac{s^2(t)}{1+(\dot{s}(t))^2}$, $G(t)={\rm
i}(S(t)\dot{\delta}(t)-s^2(t))$, $\delta=s\dot{s}/2$. The
nonhomogeneous terms $f$ and $f_1$ were chosen so that the exact
solution of the problem was given by $u(t,x)=-x(x-1)^3+\sin(t)x$.
To compare the exact with the numerical solution we calculated
the $l_2$ error at the nodes $x_j$ at $T=1$ (taking $k=h$). Table
4.1 shows the rates of convergence of the numerical solution in
the three cases. The rate is clearly two in the upsloping case
(as predicted by the theory), approaches two in the downsloping
and seems not to have stabilized in the oscillatory case. On the
other hand, as predicted by the convergence theory, (FE) when
applied to ($\mathcal{AK}$) gave clear second-order convergence.
\begin{center}
\begin{tabular}{|c|c|c|c|}\hline
   \mbox{$h$} & \mbox{Case 1} &
  \mbox{Case 2} & \mbox{Case 3} \\ \hline
  1/100 & 1.998 & 1.638 & 1.766 \\ \hline
  1/200 & 1.999 & 1.659 & 1.085 \\ \hline
  1/400 & 1.999 & 2.001 & 1.556 \\ \hline
  1/800 & 2.000 & 2.012 & 2.615\\ \hline
\end{tabular}\\
\vspace{0.5cm}
 Table 4.1. Orders of convergence of (FE) for
($\mathcal{N}$) in $l_2$ in three cases of bottom topography.
\end{center}
\subsection{Comparison of ($\mathcal{N}$) and ($\mathcal{AK}$):
The upsloping and downsloping wedge.}
We first consider the ASA upsloping wedge underwater acoustic
test problem, see \cite{JF}, with rigid bottom given in the
original variables $r$, $z$ by the function $l(r)=200-0.05r\;{\rm
m}$ for $0\leq r\leq 3339\;{\rm m}$. The source, of frequency
$f_0=25\;{\rm Hz}$, was placed at $z_s=100{\rm m}$ and modelled by
the initial value
$\psi_0(z)=\sqrt{\frac{k_0}{2}}\{\exp(-(z-z_s)^2\frac{k_0^2}{4})-
\exp(-(z+z_s)^2\frac{k_0^2}{4})\}$, $0\leq z\leq l(0)$. The water
was assumed to have constant sound speed equal to
$c=c_0=1500\;{\rm m/sec}$ and no attenuation. In (PN) $g_B(r)$
was taken equal to ${\rm i}k_0$. The problem was transformed by
the change of variables \eqref{1.6} to an equivalent one on the
horizontal strip $0\leq x\leq 1$, $0\leq t\leq T$, and it was
solved numerically by (FE) in both the ($\mathcal{N}$) and
($\mathcal{AK}$) formulations with $h=1/1000$, $k=T/1000$,
$T=3339$. (In the figures that follow we present the numerical
results after transforming them back to the original $r$, $z$
variables. Specifically, we present graphs of the numerically
computed field $\psi$, represented as is customary in underwater
acoustics, by the transmission loss function
$TL=-20\log_{10}(|\psi(z,r)|)+10\log_{10}r$ dB depicted as a
function of $r$ at certain depths $z$.) For this upsloping
example we show in Figure 2 the transmission loss curves as
functions of $r\,\in\,[0,2200\;{\rm m}]$ at a depth of $z=90\;{\rm
m}$ for both the ($\mathcal{N}$) and ($\mathcal{AK}$) models,
which evidently agree very well.
\begin{center}
\resizebox{12cm}{6.4cm}{\includegraphics{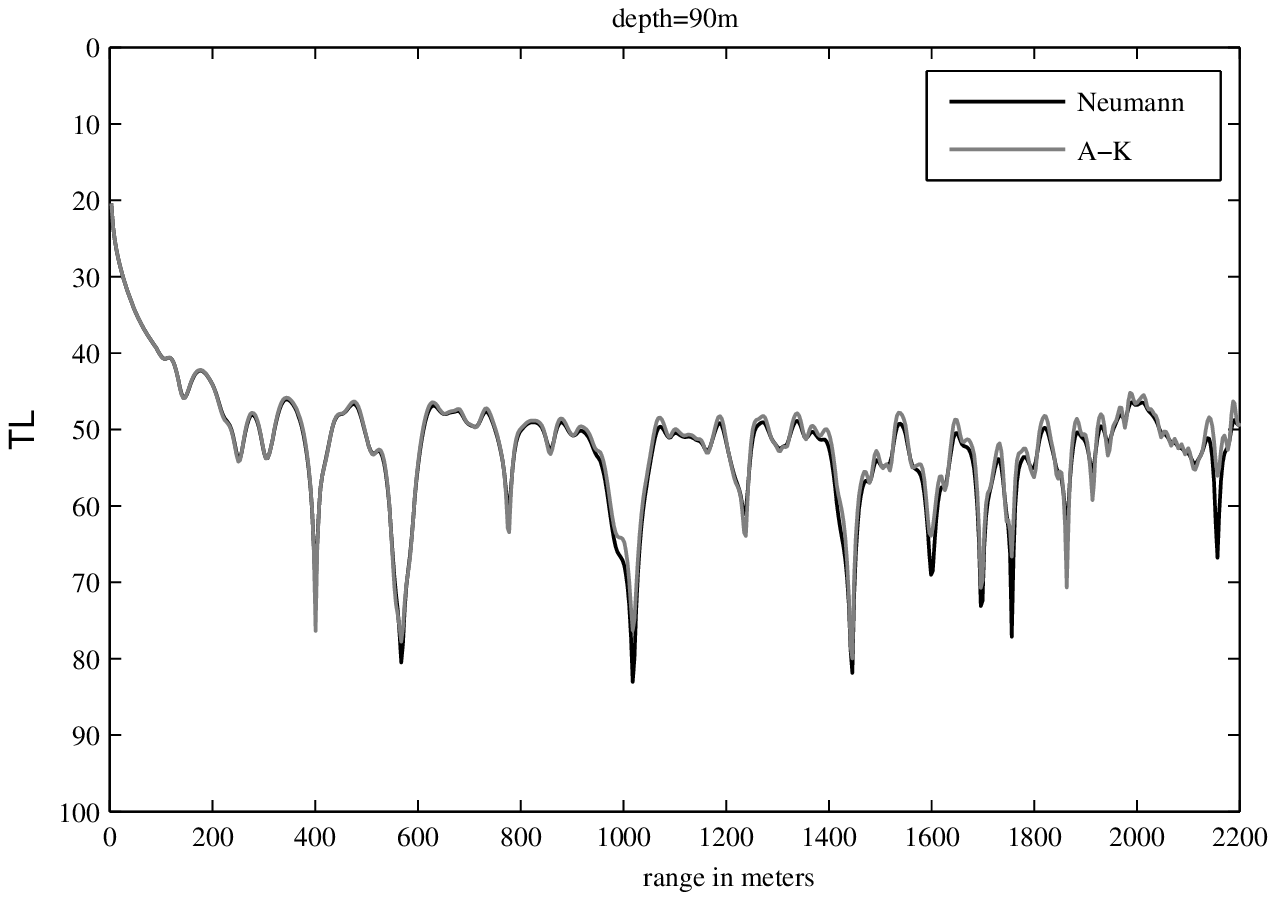}}\\
{\small Figure 2. Upsloping ASA wedge; TL as a function of $r$ at
depth $z=90\;{\rm m}$, comparison of ($\mathcal{N}$) and
($\mathcal{AK}$).}
\end{center}
We then considered the analogous downsloping wedge given by
$l(r)=33.05+0.05r$ for $0\leq r\leq 3339\;{\rm m}$. The source, of
frequency $25\;{\rm Hz}$, was placed at $z_s=25{\rm m}$ and
modelled as in the upsloping case. In this case, we found that the
(FE) numerical solution of the problem ($\mathcal{N}$) apparently
exhibited numerical instabilities and did not seem to converge as
the discretization parameters became smaller. For example, in
Figure 3 we superimpose the TL curves at depth $z=25{\rm m}$
corresponding to the ($\mathcal{N}$) model solved by (FE) with
$h=1/100$, $k=T/100$ and $h=1/1000$, $k=T/1000$, $T=3339$, with
the analogous results obtained by ($\mathcal{AK}$) solved by (FE)
with smaller $h$ and $k$. The ($\mathcal{AK}$) model, when
discretized by (FE), yields reasonable results that converge to
the solution shown in Figure 4 with dotted line. To make sure that
the numerical method used for ($\mathcal{N}$) was not the culprit,
we repeated the numerical experiment using a Crank-Nicolson
\textit{finite difference} discretization for ($\mathcal{N}$), and
found results identical to those of the (FE). We tentatively
conclude, therefore, that in this realistic downsloping bottom
case, the model ($\mathcal{N}$) allows the growth of
instabilities, in agreement with the remarks of Abrahamsson and
Kreiss in \cite{ref1} and \cite{ref2}.
\begin{center}
\resizebox{12cm}{6.4cm}{\includegraphics{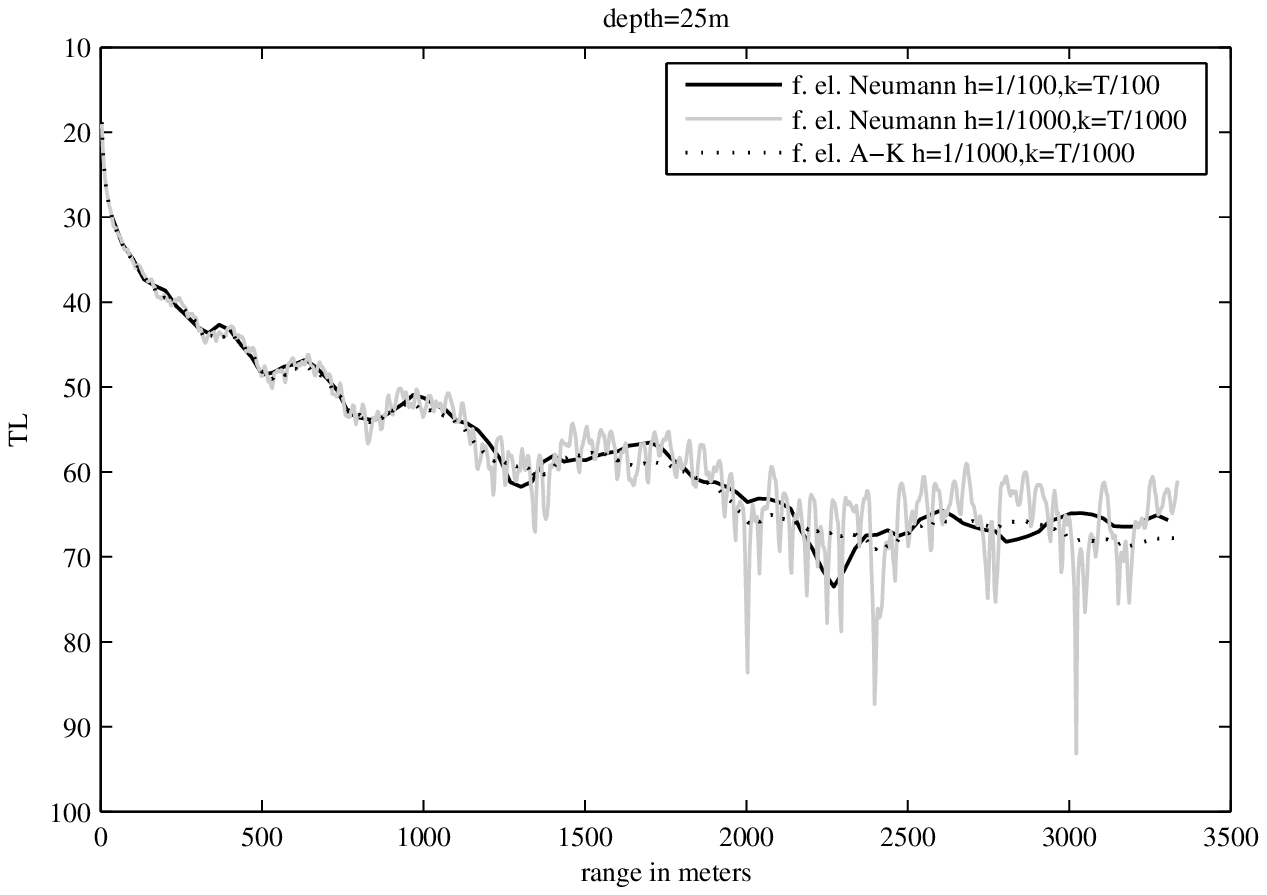}}\\
{\small Figure 3. Downsloping ASA wedge; TL as a function of $r$
at depth $z=25{\rm m}$. (FE) solutions for the ($\mathcal{N}$) and
($\mathcal{AK}$) models.}
\end{center}
To check the validity of the ($\mathcal{AK}$) solution of this
problem we compared the results of Figure 3 with those of yet
another numerical method, the Crank-Nicolson type finite
difference code IFD for the (PE), \cite{LB}, \cite{LBP},
\cite{ref42}, which has been widely used in underwater acoustic
numerical simulations.
\begin{center}
\resizebox{12cm}{5.8cm}{\includegraphics{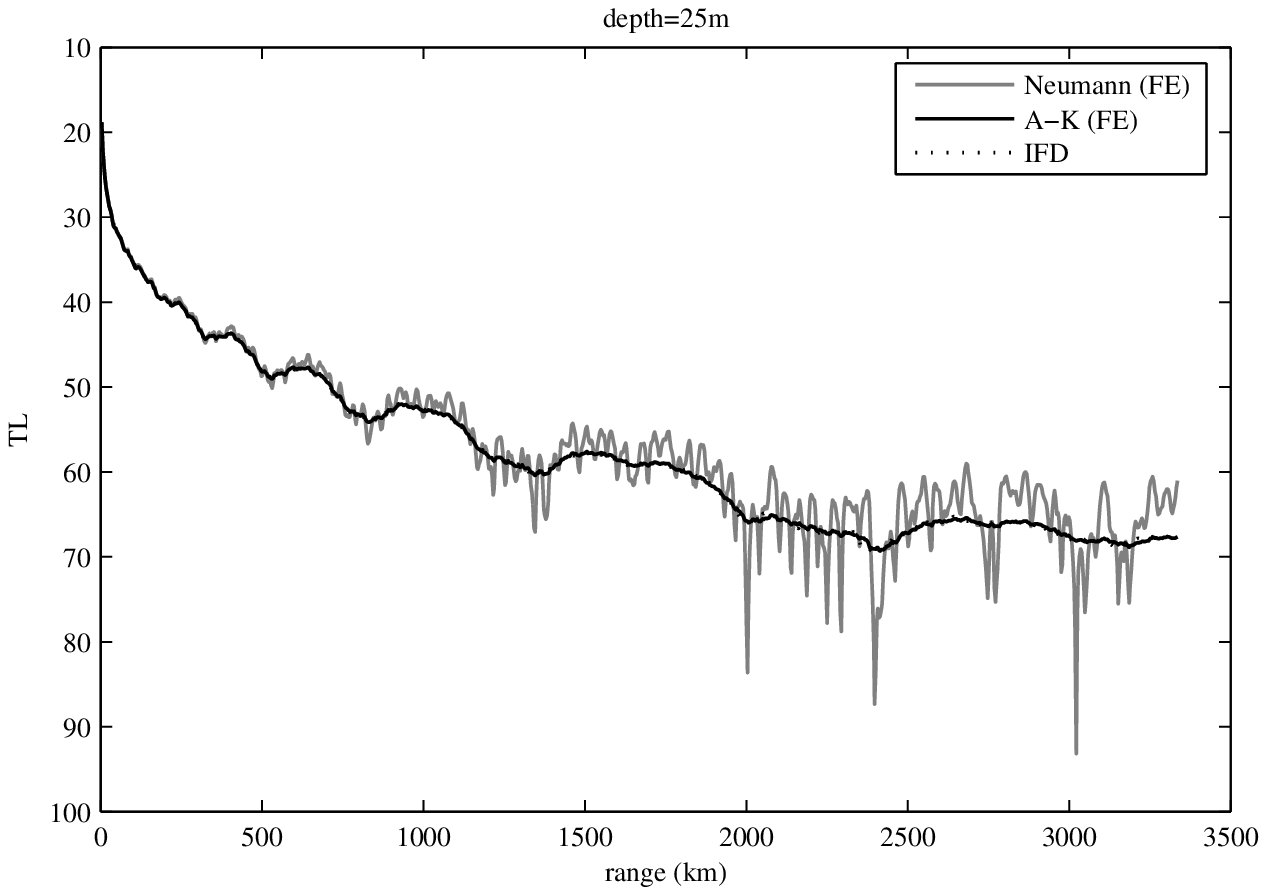}}\\
(a)\\
\resizebox{12cm}{5.8cm}{\includegraphics{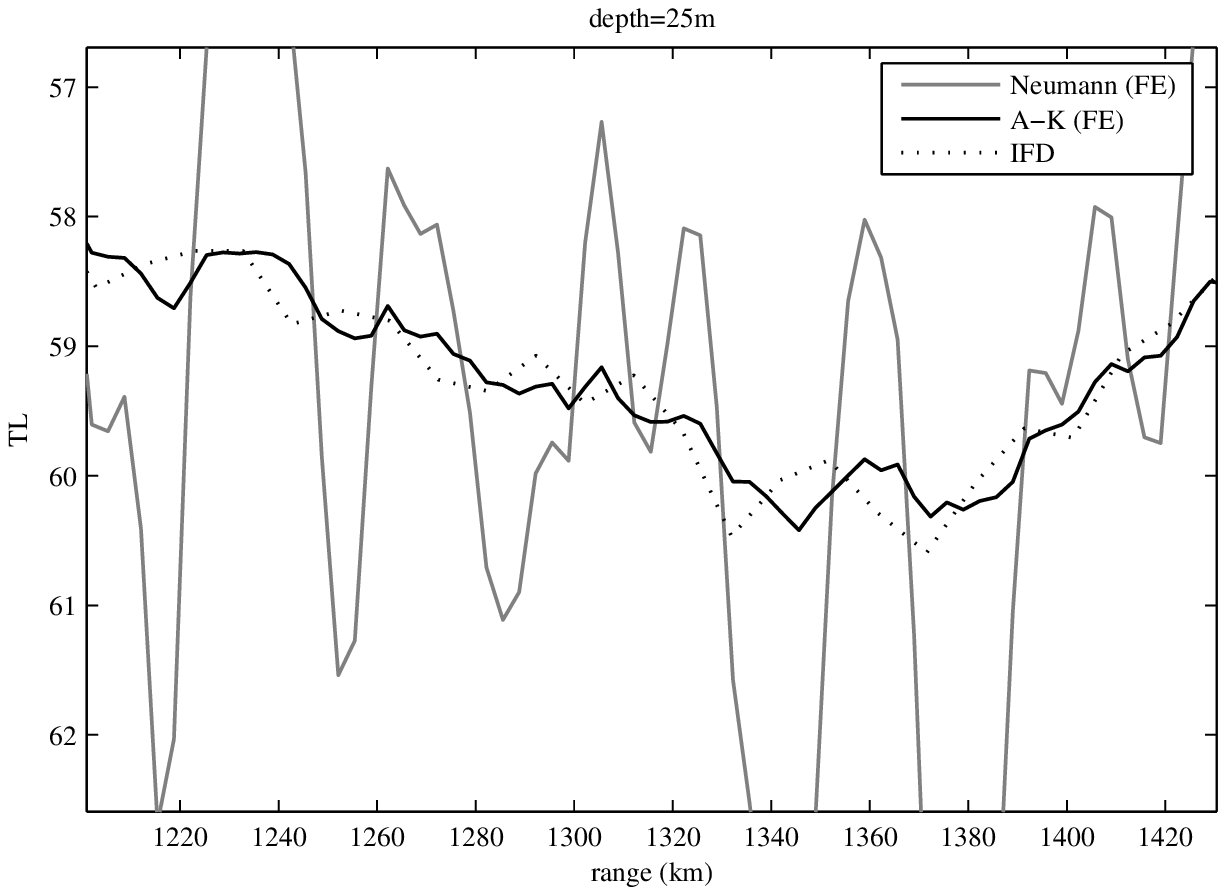}}\\
(b)\\
\resizebox{12cm}{5.8cm}{\includegraphics{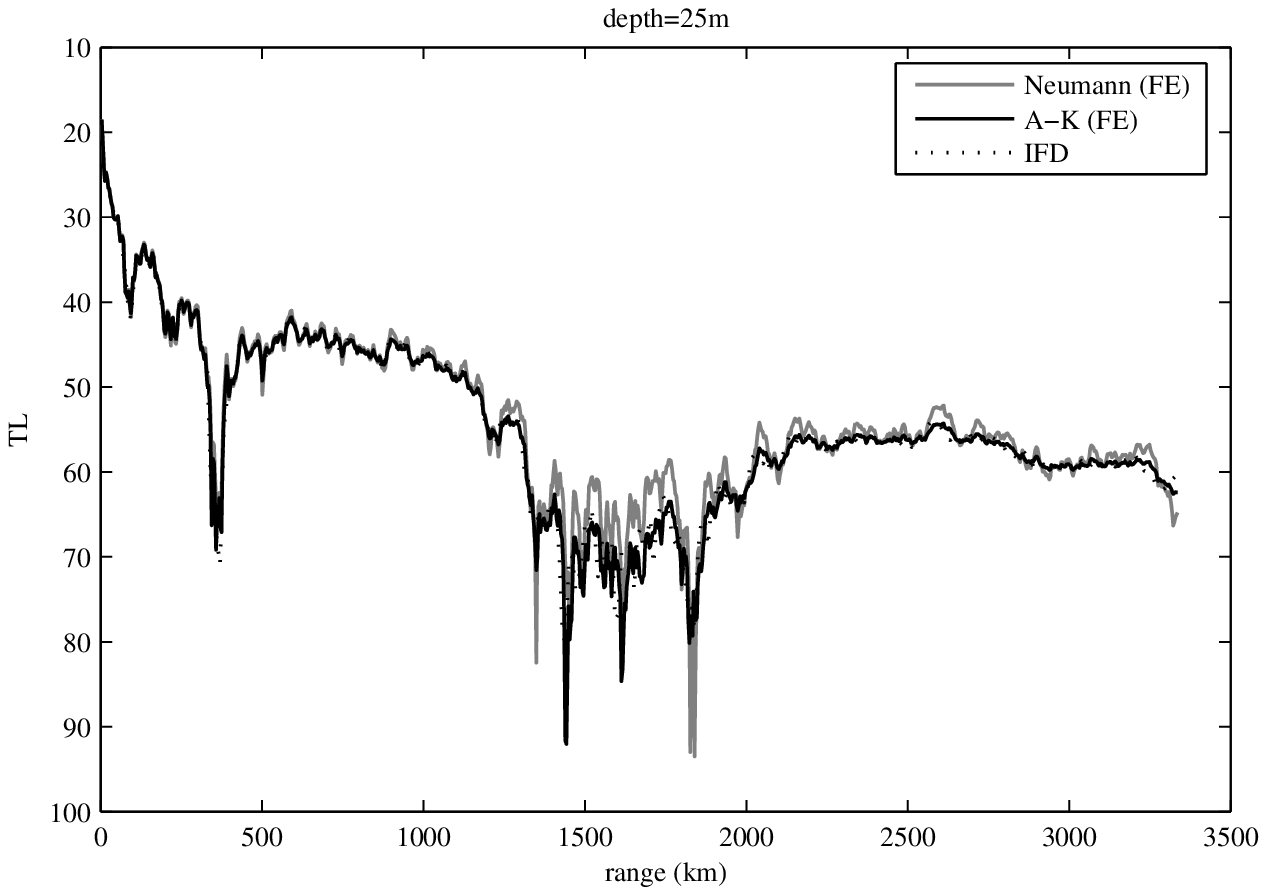}}\\
(c)\\
{\small Figure 4. Downsloping ASA wedge; TL as a function of $r$
at a depth $z=25{\rm m}$. Comparison of ($\mathcal{N}$) and
($\mathcal{AK}$), discretized by (FE), and IFD with rigid bottom
b.c. (a): $f_0=25\;{\rm Hz}$, (b): Magnification of (a) for
$r\in[1210,1430]$, (c): $f_0=80\;{\rm Hz}$.}
\end{center}
We chose the option of the rigid bottom boundary condition in IFD
and solved the problem using $\Delta z=3.31\;{\rm m}$, $\Delta
r=0.17\;{\rm m}$, values by which the IFD solution had converged.
(The IFD code solves the problem in the original $r$, $z$
wedge-shaped domain). Figure 4(a) shows the superimposed TL
curves obtained at $z=25{\rm m}$ by the ($\mathcal{N}$) and
($\mathcal{AK}$) models solved by (FE) with $h=1/1000$,
$k=T/1000$, $T=3339$ (as in Figure 3) and for the IFD with the
rigid bottom boundary condition. The results of ($\mathcal{AK}$)
and IFD agree well. In fact, they differ by about half a dB as
inspection of a typical window of Figure 4(a), shown in Figure
4(b), reveals. (It is worthwhile to note that at a higher
frequency $f_0=80\;{\rm Hz}$ the results of (FE)-($\mathcal{N}$)
approach those of (FE)-($\mathcal{AK}$) and IFD, see Figure 4(c)).
\par
To explain this result we looked closely at how IFD implements
the rigid bottom boundary condition and found that it does not
actually discretize (PN); instead, it uses a different boundary
condition obtained by replacing the $\psi_r$ term in (PN) by
$\frac{\rm i}{2k_0}\psi_{zz}+\frac{{\rm i}k_0}{2}(\eta^2-1)\psi$
using the (PE), and then discretizing the $\psi_{zz}$ term at the
bottom with one-sided finite differences from the interior of the
domain. In the next paragraph we offer an explanation why this
rigid bottom boundary condition yields a stable problem for any
monotone bottom profile.

Our tentative conclusion, then, from this experiment is that in
the case of realistic, downsloping environments, ($\mathcal{AK}$)
and the rigid bottom boundary condition model implemented by IFD
apparently yield correct results, while the Neumann bottom
boundary condition used in ($\mathcal{N}$), which retains the term
$\psi_r$ at the bottom, allows the growth of instabilities.
\subsection{Using the p.d.e. in the dynamical boundary condition}
Let $w=w(t,y)$ be defined for $0\leq y\leq s(t)$, $0\leq{t}\leq
T$, and satisfy \eqref{1.1}--\eqref{1.4}. Replace the term
$w_t(t,s(t))$ in \eqref{1.4} by its value given by the p.d.e. in
\eqref{1.1} to obtain
\begin{equation}\label{3.1}
w_y(t,s(t))-\dot{s}(t)\,\Big\{\,\tfrac{\cunit}{2} \,w_{yy}(t,s(t))
+\big[{\rm i}\,\gamma(t,s(t))+g(t)\big]
\,w(t,s(t))\,\Big\}=0\quad\forall\,t\in[0,T].
\end{equation}
In the IFD code, the rigid bottom boundary condition used is a
finite difference discretization of \eqref{3.1}.
\par
To avoid the presence of the second derivative $w_{yy}(t,s(t))$
in the boundary condition \eqref{3.1}, we differentiate
\eqref{1.1} with respect to $y$ and put
${\widetilde{p}}(t,y)=w_y(t,y)$. (Note that
$w(t,y)=\int_0^y{\widetilde{p}}(t,\xi)\,d\xi$ since $w(t,0)=0$.)
Then, the ibvp \eqref{1.1}--\eqref{1.3}, \eqref{3.1} becomes
\begin{equation}\label{3.2}
\begin{split}
&{\widetilde{p}}_t= \tfrac{\cunit}{2}\,{\widetilde{p}}_{yy}
+\cunit\,\gamma(t,y)\,{\widetilde{p}} +\cunit\,\gamma_y(t,y)\,w
\quad\,\forall\,y\in[0,s(t)],\ \ \forall\,t\in[0,T],\\
&{\widetilde{p}}_y(t,0)=0\quad\forall\,t\in[0,T],\\
&{\widetilde{p}}(t,s(t))-\dot{s}(t)\left\{
\tfrac{\cunit}{2}\,{\widetilde{p}}_y(t,s(t))
+\left[\,\cunit\,\gamma(t,s(t))+g(t)\,\right]
\,w(t,s(t))\right\}=0\quad\forall\,t\in[0,T],\\
&{\widetilde{p}}(0,y)={\widetilde{p}}_0(y) :=w_0'(y)
\quad\forall\,y\in[0,s(0)].\\
\end{split}
\end{equation}
(Note that using the \eqref{1.1} at $y=0$ and the surface
boundary condition $w(t,0)=0$, we obtain that
${\widetilde{p}}_y(t,0)=w_{yy}(t,0)=0$.)
\par
In what follows, we shall obtain an \textit{a priori} $L^2$ bound
for the solution of \eqref{3.2} and then propose a finite element
method for solving it. With this aim in mind, we perform as usual
the range-dependent change of depth variable $x:=\frac{y}{s(t)}$
that maps the domain of the problem onto the horizontal strip
$\{(t,x)\,:\,t\,\in\,[0,T],\,x\,\in\overline{D}\}$, where
$D=(0,1)$. Consider the transformation
\begin{equation}\label{3.3}
{\widetilde{p}}(t,y)=\tfrac{1}{s(t)}\,\exp(-\zeta(t,x))
\,\left(\,p(t,x)-\zeta_x(t,x) \,\int_0^x p(t,\xi)\,d\xi\,\right),
\end{equation}
where the function $\zeta$ will be specified below. Note that the
function $\theta$, defined by
$\theta(t,x):=\int_0^xp(t,\xi)\,d\xi$ for $(t,x)\in [0,T]\times
\overline{D}$, satisfies the first-order o.d.e.
\begin{equation*}
\theta_x(t,x)-\zeta_x(t,x)\,\theta(t,x)=s(t) \,\exp(\zeta(t,x))
\,{\widetilde{p}}(t,x\,s(t)).
\end{equation*}
Solving this differential equation with initial condition
$\theta(t,0)=0$ yields
\begin{equation*}
\theta(t,x)
=s(t)\,\exp(\zeta(t,x))\int_0^x{\widetilde{p}}(t,\xi\,s(t))\,d\xi,
\end{equation*}
from which we may derive the inverse of the transformation
\eqref{3.3} in the form
\begin{equation*}\label{APP_A4}
p(t,x)=s(t)\exp(\zeta(t,x)) \,\left({\widetilde{p}}(t,x s(t))
+\zeta_x(t,x)\int_0^x{\widetilde{p}}(t,\xi\,s(t))\,d\xi\right).
\end{equation*}
After some calculations we also obtain that
\begin{equation}\label{3.4}
\theta(t,x)=\exp(\zeta(t,x))\,w(t,xs(t)),
\quad(t,x)\in[0,T]\times \overline{D}.
\end{equation}
\par
Following the ideas of \cite{ref8}, and after analogous
computations (see, in particular, (2.7) and (2.8) of
\cite{ref8}), we may deduce that $p$ solves a well posed ibvp, in
the case of strictly monotone bottoms, i.e. when $\dot{s}(t)$ is
either positive or negative for all $t\in[0,T]$.
To see this, define first $\zeta$, as in \cite{ref8}, by the
formula
\begin{equation}\label{3.5}
\zeta(t,x)=\tfrac{\cunit}{2}\,(\sigma(t)-1)\,\dot{s}(t)s(t)\,x^2
\quad\forall\,(t,x)\in [0,T]\times \overline{D},
\end{equation}
where
$\sigma(t):=\frac{2(1+\dot{s}(t)^2)}{\dot{s}(t)^2}+\varepsilon$,
if $\dot{s}(t)>0$, where $\varepsilon$ is a positive constant, and
$\sigma(t):=1$, or equivalently $\zeta=0$, if $\dot{s}(t)< 0$.
Then, in the transformed domain, and expressed in terms of the new
field variables $p$ and $\theta$, the ibvp \eqref{3.2} becomes
\begin{equation}\label{3.6}
\begin{split}
&p_t=\tfrac{\cunit}{A(t)}\,p_{xx}
+B(t,x)\,p_x+\left[\,B_x(t,x)+G(t,x)\,\right]\,p
+G_x(t,x)\,\theta\quad
\forall\,(t,x)\in[0,T]\times{\overline{D}},\\
&p_x(t,0)=0\quad\forall\,t\in[0,T],\\
&\cunit\,\tfrac{1}{A(t)}\,p_x(t,1)=
\tfrac{1-R_1(t)\,B(t,1)}{R_1(t)}\,p(t,1)
-\tfrac{R_1(t)\,G(t,1)+R_2(t)}{R_1(t)}\,\theta(t,1)
\quad\forall\,t\in[0,T],\\
&p(0,x)=p_0(x)
\quad\forall\, x\in{\overline{D}},\\
\end{split}
\end{equation}
where
$p_0(x)=s(0)\exp(\zeta(0,x))\left[ w'_{0}(xs(0))
+\zeta_x(0,x)\int_0^x w'_{0}(\xi s(0))\;d\xi\right]$,
$A(t)=2\,s^2(t)$,
$R_1(t)=\tfrac{\dot{s}(t)s(t)}{1+(\dot{s}(t))^2}$,
$B(t,x)=x\,\frac{\dot{s}(t)}{s(t)}
-\frac{\cunit}{s^2(t)}\,\zeta_x(t,x)$,
$G(t,x)=\zeta_t(t,x)-x\,\tfrac{\dot{s}(t)}{s(t)}\,\zeta_x(t,x)
+\cunit\,\gamma(t,xs(t))+\tfrac{\cunit}{2s^2(t)}
\,\left[\,(\zeta_x(t,x))^2-\zeta_{xx}(t,x)\,\right]$,
$R_2(t)=[\,g(t)-\zeta_t(t,1)\,]\,R_1(t)+\zeta_x(t,1)$.
(Recall that $\theta(t,x)=\int_0^xp(t,\xi)\;d\xi$. In addition,
note that \eqref{3.5} yields that $B$ is real-valued and is given
by $B(t,x)=x\tfrac{\dot{s}(t)}{s(t)}\,\sigma(t)$, so that
$B(t,0)=0$ and
$B_x(t,x)=B(t,1)=\tfrac{\dot{s}(t)}{s(t)}\,\sigma(t)$. It is
easily checked that $1-R_1(t)\,B(t,1)\not=0$ for $t\in[0,T]$.)
We may now prove the following result.
%
%
\begin{thm}\label{Qewrhma1}
If the bottom is strictly monotone, the ibvp \eqref{3.6} is
$L^2$-stable.
\end{thm}
%
%
\begin{proof}
Multiply the p.d.e. in \eqref{3.6} by $\overline{p(t,x)}$,
integrate with respect to $x$ in $[0,1]$, use integration by
parts, and take real parts to obtain
\begin{equation*}
\begin{split}
\tfrac{1}{2}\tfrac{d}{dt}\|p(t,\cdot)\|^2&=
\,\tfrac{2-R_1(t)\,B(t,1)}{2R_1(t)}\,|p(t,1)|^2 -{\rm
Re}\left[\,\tfrac{G(t,1)\,R_1(t)+R_2(t)}{R_1(t)}
\,\theta(t,1)\,\overline{p(t,1)}\,\right]\\
&\quad+{\rm Re}\big(\,G_x(t,\cdot)\,\theta(t,\cdot),p(t,\cdot)\,\big)\\
&\quad+\tfrac{1}{2}\,\big(\,B_x(t,\cdot)p(t,\cdot),p(t,\cdot)\,\big)
+{\rm Re}\big(\,G(t,\cdot)\,p(t,\cdot),p(t,\cdot)\,\big)
\quad\forall\,t\in[0,T].\\
\end{split}
\end{equation*}
Using the Cauchy-Schwarz inequality, the arithmetic-geometric mean
inequality, and noting that $|\theta(t,1)|\leq \|p(t,\cdot)\|$,
$\|\theta(t,\cdot)\|\leq \|p(t,\cdot)\|$, we see from the above
that for any $\xi>0$ there exists a constant $C_{\xi}>0$ such that
\begin{equation*}
\tfrac{d}{dt}\|p(t,\cdot)\|^2\leq
\left(\,\tfrac{1}{R_1(t)}-\tfrac{1}{2}\,B(t,1)
+\xi\,\right)\,|p(t,1)|^2 +C_{\xi}\,\|p(t,\cdot)\|^2
\quad\forall\,t\in[0,T].
\end{equation*}
Since $\frac{1}{R_1(t)}-\frac{1}{2}B(1,t)<0$ for $t\in [0,T]$, we
may chose $\xi$ sufficiently small to make the first term in the
right-hand side of the above negative. Hence, by Gr{\"o}nwall's
lemma, we conclude that $\|p(t,\cdot)\|\leq\,C\|p_0\|$ for
$t\in[0,T]$, which ends the proof. $\Box$
\end{proof}
\par
Now, we can define a semiscrete approximation
$p_h:[0,T]\rightarrow S_h$ of the solution $p$ of problem
\eqref{3.6} by
\begin{equation*}
p_h(0,x)=p_h^0(x)\quad\forall\,x\in[0,1],
\end{equation*}
and
\begin{equation*}
\begin{split}
(\partial_t p_h,\phi)=&-\tfrac{\cunit}{A(t)}\,{\mathcal
B}(p_h,\phi) +\Big[\,\tfrac{1-R_1(1)B(t,1)}{R_1(t)}\,p_h(t,1)
-\tfrac{R_1(t)G(1,t)+R_2(t)}{R_1(t)}\,\theta_h(t,1)\Big]
\,\overline{\phi(1)}\\
& +\big(\,B(t,\cdot)\,\partial_xp_h,\phi\,\big)
+\big(\,[B_x(t,\cdot)+G(t,\cdot)]p_h,\phi\,\big)\\
&+\big(\,G_x(t,\cdot)\theta_h,\phi\,\big) \quad\forall\,\phi\in
S_h,\quad\forall\,t\in[0,T],
\end{split}
\end{equation*}
where $\theta_h(t,x):=\int_0^x p_h(t,\xi)\,d\xi$ and $p^0_h\in
S_h$ is a given reasonable approximation of $p_0$.
Consequently, using \eqref{3.4}, we see that
$\exp(-\zeta(t,x))\theta_h(t,x)$ is an approximation of the
solution $w(t,xs(t))$ of the ibvp \eqref{1.1}-\eqref{1.4}.
Also, it follows, as in Theorem~\ref{Qewrhma1}, that there exists
a positive constant $C$ such that
$\|p_h(t,\cdot)\|\leq\,C\,\|p_h^0\|$ for $t\in[0,T]$.
\subsection{Growth of solutions of $(\mathcal{N})$
for various bottom shapes.}
The final set of numerical experiments that we report concern the
behavior of the size of the solutions of $(\mathcal{N})$, as $t$
grows, in the presence of bottom profiles of various shapes.
Recall that in \cite{ref1} it was shown that $(\mathcal{N})$ is
well posed if $s$ is strictly monotone, i.e. if $\dot{s}(t)>0$ or
$\dot{s}(t)<0$ for $0\leq t\leq T$. In addition, downsloping
bottom profiles were identified for which the solution of
$(\mathcal{N})$ grew exponentially with $t$. (The fact that
problems may arise in case $\dot{s}$ changes sign may be
expected, in view of the analogous difficulties encounterd in the
(real) parabolic case, cf.e.g. \cite{Ba-Re}.)
\par
The ibvp $(\mathcal{N})$ was solved numerically with the (FE)
method up to $T=1$, with $\beta=f=g=0$, $u_0(x)=-x(x-1)^3$,
$0\leq x\leq 1$, with mesh parameters $h=k=1/500$, in the case of
the eight bottom profiles $s(t)$, $0\leq t\leq 1$, labeled (a) to
(h) and shown in the left-hand icons of the pairs in Figure 5.
(In all cases depth increases downwards.) The right-hand icon
shows the corresponding, numerically computed, $L^2$-norm of the
solution of $(\mathcal{N})$ $\|u(t,\cdot)\|$ for $0\leq t\leq 1$.
(Note that $\|u(0,\cdot)\|=\frac{1}{6\sqrt{7}}\cong 0.062994$.)
The bottom profiles are given for $0\leq t\leq 1$ by the
expressions:
(a) $s(t)=e^t$, (b) $s(t)=e^{-t}$,
(c) $s(t)=1+(t-0.5)^2$, (d) $s(t)=1-|t-0.5|^3$,
(e) $s(t)=1-(t-0.5)^3$,
(f) $s(t)=2-|2t-1|$,
(g) $s(t)=1+(t-0.5)^3$ and (h) $s(t)=1+t^3$.
%
%
\begin{center}
\resizebox{8cm}{7.0cm}{\includegraphics{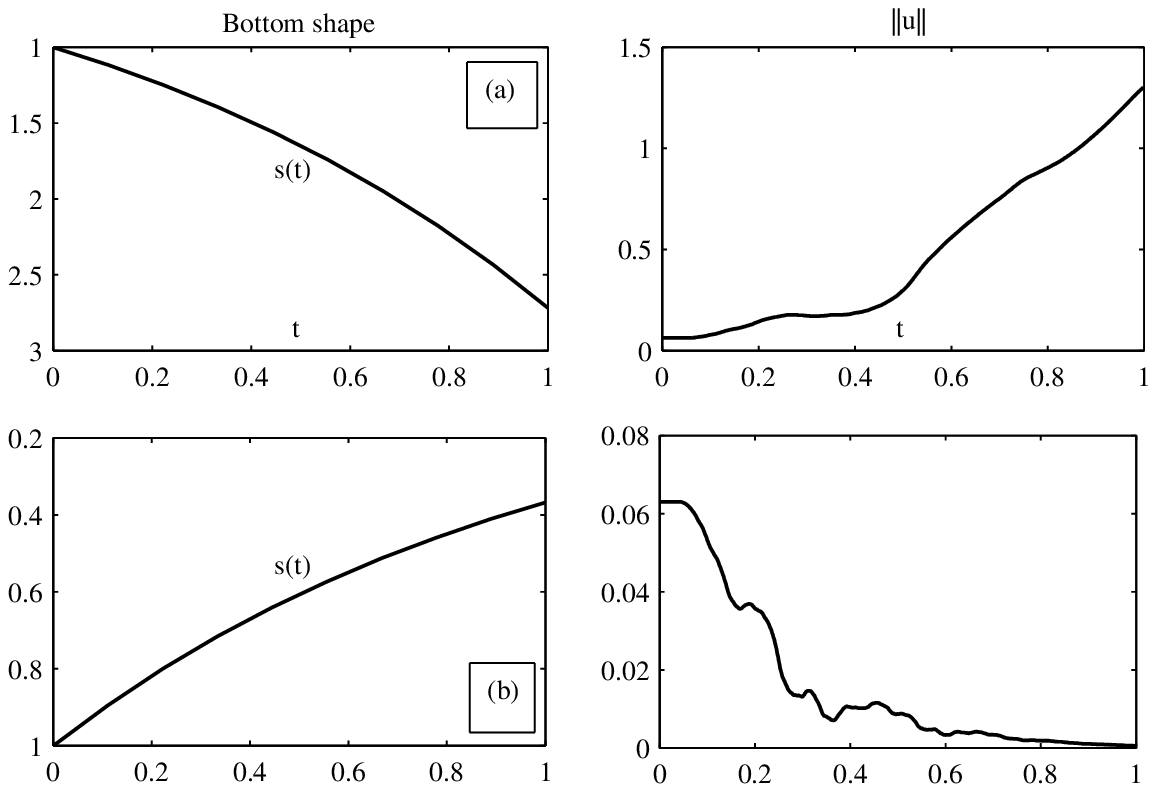}}
\resizebox{8cm}{7.4cm}{\includegraphics{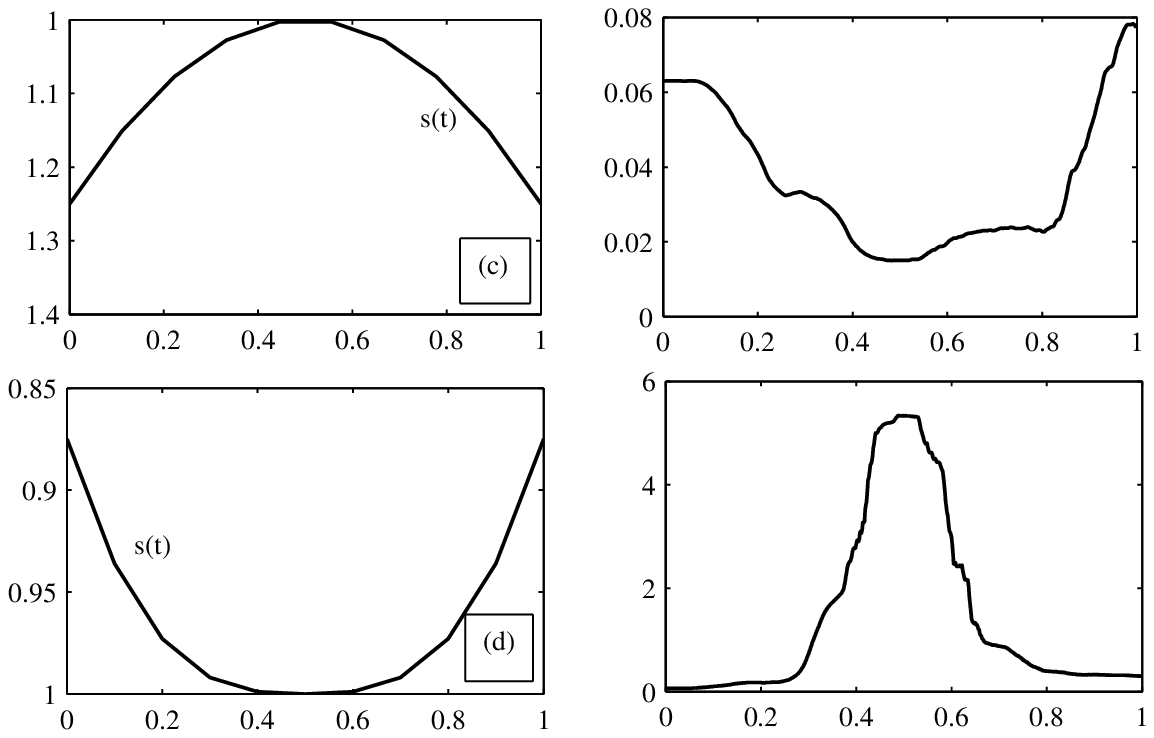}}\\
\resizebox{8cm}{7.2cm}{\includegraphics{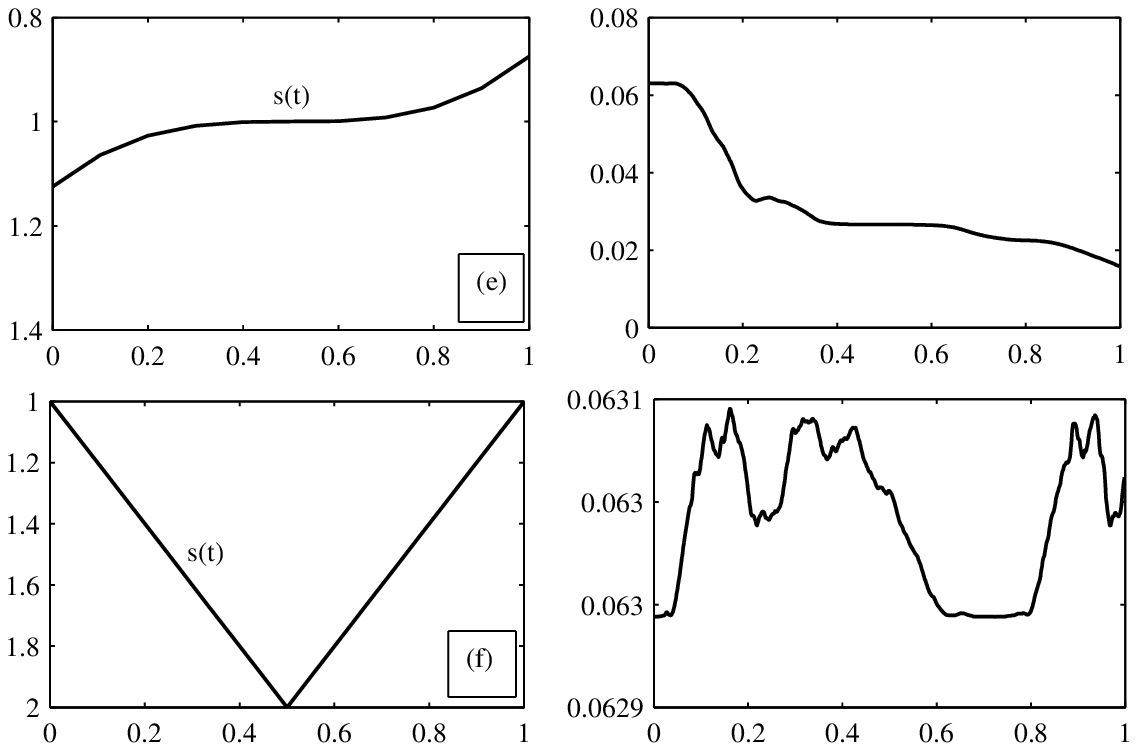}}
\resizebox{8cm}{7.4cm}{\includegraphics{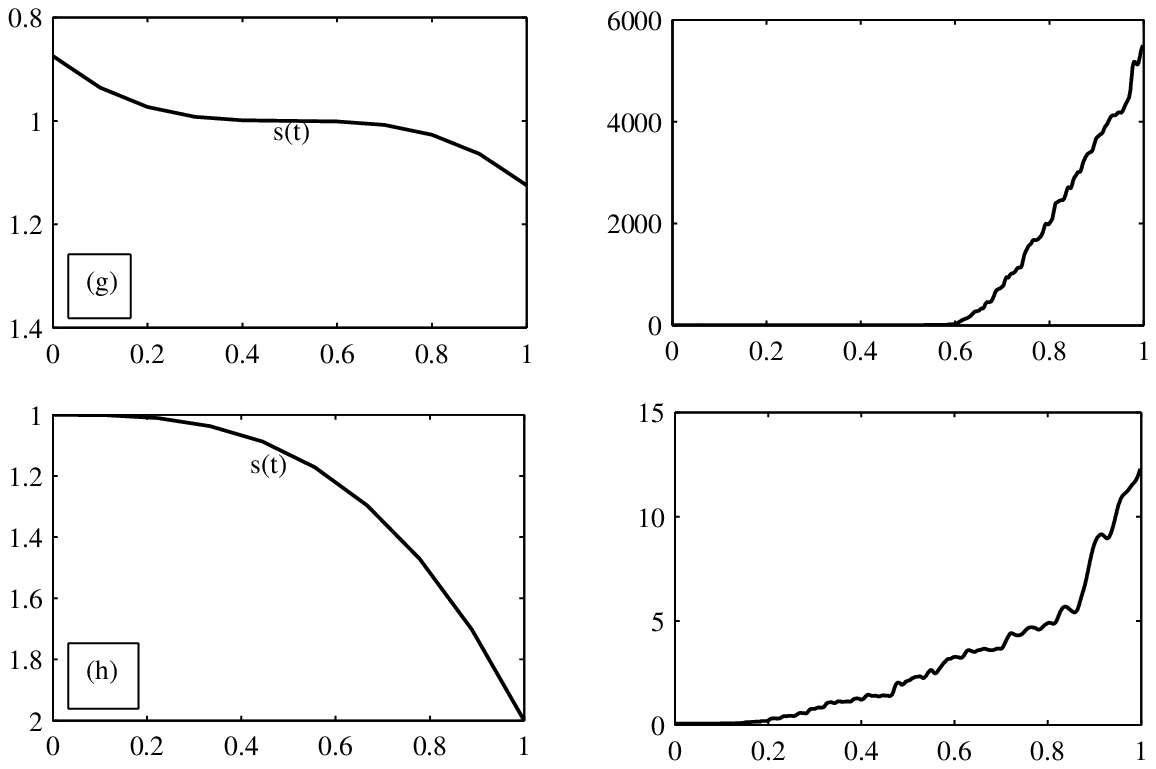}}\\
{\small Figure 5. Behavior of the $L^2$-norm of the numerical
solution of $(\mathcal{N})$ as a function of $t$ for various
bottom profiles $s(t)$.}
\end{center}
%
%
%
Only (a) and (b) correspond to strictly monotone profiles for
which the theory of \cite{ref1} properly applies. In the cases
(c), (d), (f) there is a change in monotonicity, in (e) and (g)
we have that $\dot{s}(t)=\ddot{s}(t)=0$ at $t=1/2$, while in (h)
there holds that $\dot{s}(0)=\ddot{s}(0)=0$. (In the case (f) a
$t$-mesh node was placed at $t=0.5$, where $\dot{s}$ fails to
exist.)
\par
We observe that the solution maintains a small $L^2$-norm in
upsloping, like (b), or eventually upsloping bottoms, as in the
cases of the trenches (d) and (f). There is a considerable growth
of $\|u\|$ in the examples wherein the bottom profile is
eventually downsloping, see (a), (c), (g), and (h), in agreement
with the observations in \cite{ref1}, \cite{ref2}. We note that
in the case (g), an apparent singularity develops at $t=1/2$,
where the bottom curvature changes sign (with horizontal tangent)
and the bottom becomes downsloping. This apparently causes the
$L^2$-norm to grow violently for $t>1/2$. A relatively weaker,
but sizeable growth is also observed in (h), where the bottom is
such that $\dot{s}=\ddot{s}=0$ at $t=0$ and is monotonically
downsloping for $t>0$. One cannot be of course certain about the
existence of a singularity at $t=1/2$ in the case (g), given that
the (FE) code does not at present possess an adaptive refinement
capability in $x$ and $t$. However, when the experiment was
repeated with $k=h=1/800$, it was confirmed that the onset of
rapid growth occurred at about $t=1/2$; for this mesh size,
$\|u\|$ became of order $O(10^4)$ at $t=1$.

%
\section{A parabolic problem with a dynamical boundary condition}
%
%
Here we consider the model one-dimensional (real) parabolic
problem \eqref{1.14} with a dynamical boundary condition
analogous to that of $({\mathcal N})$, which we re-write here for
ease in reading: We seek a real-valued function
$u\,:\,[0,T]\times[0,1]\rightarrow \mathbb{R}$, such that
\begin{equation}\label{TNT1}
\begin{split}
&u_t=a(t)\,u_{xx}+\beta(t,x)\,u+f(t,x)\quad\forall
\,(t,x)\in[0,T]\times[0,1],\\
&u(t,0)=0\quad\forall\,t\in[0,T],\\
&a(t)\,u_x(t,1)=\varepsilon(t)\,u_t(t,1)+\delta(t)\,u(t,1)+g(t)
\quad\forall\,t\in[0,T],\\
&u(x,0)=u_0(x)\quad\forall\,x\in[0,1],\\
\end{split}
\end{equation}
where $\beta:[0,T]\times[0,1]\rightarrow{\mathbb R}$,
$f:[0,T]\times[0,1]\rightarrow{\mathbb R}$,
$\delta:[0,T]\rightarrow{\mathbb R}$, $g:[0,T]\rightarrow{\mathbb
R}$, $a:[0,T]\rightarrow(0,+\infty)$ with
$a_{\star}:=\inf_{[0,T]}a>0$,
$\varepsilon:[0,T]\rightarrow{\mathbb R}$,
$u_0:[0,1]\rightarrow{\mathbb R}$, are given smooth functions. We
shall construct and analyze Galerkin-finite element
approximations for the solution of \eqref{TNT1}, considering two
different cases depending on the sign of $\varepsilon$.
\subsection{The dissipative case}
The dissipative case is characterized by the assumption
$\varepsilon(t)\leq0$ for $t\in[0,T]$; the problem is well posed,
see e.g. \cite{Esc1993}. We assume that its solution is smooth
enough for the purposes of the error estimates to follow. We
adopt the notation and the assumptions of Section~\ref{subsec40},
restricting ourselves to the real case, and avoiding the inverse
inequality \eqref{2.36}.
%
%
\subsubsection{Semidiscrete approximation}
%
%
Find $u_h:[0,T]\rightarrow S_h$, a space-discrete approximation
of $u$, requiring
\begin{equation}\label{TNT2}
\begin{split}
(\partial_t
u_h(t,\cdot),\chi)=&\,\big[\,\varepsilon(t)
\,\partial_tu_h(t,1)+
\delta(t)\,u_h(t,1)+g(t)\,\big]\,\chi(1)\\
&-a(t)\,{\mathcal B}(u_h(t,\cdot),\chi)
+(\beta(t,\cdot)\,u_h(t,\cdot),\chi)
+(f(t,\cdot),\chi)\quad\forall\,\chi\in S_h,
\quad\forall\,t\in[0,T],\\
\end{split}
\end{equation}
and
\begin{equation}\label{TNT3}
u_h(0,\cdot)=R_hu_0(\cdot).
\end{equation}
%
%
%
%
\begin{proposition}
If $\varepsilon\leq0$, then the problem
\eqref{TNT2}-\eqref{TNT3} admits a unique
solution $u_h\in C^1([0,T];S_h)$.
\end{proposition}
%
%
%
\begin{proof}
The result follows if we argue along the lines of the proof of
Proposition~\ref{2.3.3**}. $\Box$
\end{proof}
%
%
%
%
\begin{proposition}\label{P_TNT1}
Let $u$ be the solution of \eqref{TNT1} and $u_h$ its
semidiscrete approximation defined by \eqref{TNT2}-\eqref{TNT3}.
If $\varepsilon(t)\leq0$ for $t\in[0,T]$, then
\begin{equation}\label{TNT41}
\|u_h(t,\cdot)-R_hu(t,\cdot)\|_1^2\leq\,C\,h^{2(r+1)}
\,\left(\int_0^t\Gamma_{\ssy D}(\tau)\;d\tau\right)
\quad\forall\,t\in[0,T],\ \ \forall\,h\in(0,h_{\star}],
\end{equation}
where
$\Gamma_{\ssy D}(\tau):=\|u(\tau,\cdot)\|^2_{r+1}
+\|\partial_tu(\tau,\cdot)\|^2_{r+1}$.
\end{proposition}
%
%
%
%
%
%
%
\begin{proof}
Let $h\in(0,h_{\star}]$ and $\theta_h:=u_h-R_hu$. Using
\eqref{TNT2}, the p.d.e. in \eqref{TNT1}, \eqref{R_def} and
\eqref{4.40} we obtain
\begin{equation}\label{TNT5}
\begin{split}
(\partial_t\theta_h(t,\cdot),\chi)&=\big[\,
\varepsilon(t)\,\partial_t\theta_h(t,1)
+\delta(t)\,\theta_h(t,1)\,\big]\,\chi(1)
-a(t)\,{\mathcal B}(\theta_h(t,\cdot),\chi)\\
&\quad+(\beta(t,\cdot)\,\theta_h(t,\cdot),\chi)
+\big(\Phi_{\star}(t,\cdot),\chi\big)
\quad\forall\,\chi\in S_h,\quad\forall\,t\in[0,T],\\
\end{split}
\end{equation}
where $\Phi_{\star}:=[\partial_tu-R_h(\partial_tu)]
-\beta\,(u-R_hu)$. First we observe that from \eqref{4.39} it
follows that
\begin{equation}\label{TNT55}
\|\Phi_{\star}(t,\cdot)\|\leq\,C\,h^{r+1}
\,\big[\,\|u(t,\cdot)\|_{r+1}+
\|\partial_tu(t,\cdot)\|_{r+1}\,\big]\quad\forall\,t\in[0,T].
\end{equation}
\par
Setting $\chi=\theta_h$ in \eqref{TNT5} and
using \eqref{TraceIneq}, we obtain
\begin{equation*}\label{TNT6}
\begin{split}
\tfrac{d}{dt}\|\theta_h(t,\cdot)\|^2\leq&\,
\tfrac{|\varepsilon(t)|}{\epsilon} \,|\partial_t\theta_h(t,1)|^2
-2\,a_{\star}\,|\theta_h(t,\cdot)|_1^2\\
&+2\,\|\Phi_{\star}(t,\cdot)\|\,\|\theta_h(t,\cdot)\|
+2\,\max_{[0,1]}|\beta(t,\cdot)|\,\|\theta_h(t,\cdot)\|^2\\
&\,+2\,\left(\,2\,|\delta(t)|
+\epsilon\,|\varepsilon(t)|\,\right)
\,\|\theta_h(t,\cdot)\|\,|\theta_h(t,\cdot)|_1
\quad\forall\,t\in[0,T],\ \ \forall\,\epsilon>0,\\
\end{split}
\end{equation*}
which, along with \eqref{TNT55}, yields
\begin{equation}\label{TNT61}
\epsilon\,\tfrac{d}{dt}\|\theta_h(t,\cdot)\|^2\leq\,
|\varepsilon(t)|\,|\partial_t\theta_h(t,1)|^2
+C\,(\epsilon+\epsilon^3)\,\|\theta_h(t,\cdot)\|^2
+\epsilon\,h^{2(r+1)}\,\Gamma_{\ssy D}(t) \quad\forall\,t\in[0,T],
\ \ \forall\,\epsilon>0.
\end{equation}
\par
Set $\chi=\partial_t\theta_h$ in \eqref{TNT5} to obtain
\begin{equation*}\label{TNT7}
\begin{split}
\tfrac{d}{dt}\big[a(t)\,|\theta_h(t,\cdot)|_1^2
-\delta(t)\,|\theta_h(t,1)|^2\,\big]&=
-2\,|\varepsilon(t)|
\,|\partial_t\theta_h(t,1)|^2
-{\dot\delta}(t)\,|\theta_h(t,1)|^2\\
&\quad+{\dot a}(t)\,|\theta_h(t,\cdot)|_1^2
-2\,\|\partial_t\theta_h(t,\cdot)\|^2\\
&\quad +2\,(\beta(t,\cdot)\,\theta_h(t,\cdot),\partial_t\theta_h(t,\cdot))
+2\,(\Phi_{\star}(t,\cdot),
\partial_t\theta_h(t,\cdot))\quad\forall\,t\in[0,T],\\
\end{split}
\end{equation*}
which, in view of \eqref{TraceIneq} and \eqref{TNT55}, yields that
\begin{equation}\label{TNT8}
\tfrac{d}{dt}\big[a(t)\,|\theta_h(t,\cdot)|_1^2
-\delta(t)\,|\theta_h(t,1)|^2\,\big]
\leq\,-2\,|\varepsilon(t)|\,|\partial_t\theta_h(t,1)|^2
+C\,\|\theta_h(t,\cdot)\|_1^2
+h^{2(r+1)}\,\Gamma_{\ssy D}(t)
\quad\forall\,t\in[0,T].
\end{equation}
\par
For positive $\epsilon$ we define
\begin{equation}\label{TNT70}
\nu_{\epsilon}(t):=\epsilon\,\|\theta_h(t,\cdot)\|^2
+a(t)\,|\theta_h(t,\cdot)|_1^2
-\delta(t)\,|\theta_h(t,1)|^2
\quad\forall\,t\in[0,T].
\end{equation}
Then, applying the trace inequality \eqref{TraceIneq}, we have
\begin{equation}\label{TNT72}
\begin{split}
\nu_{\epsilon}(t)\ge&\,\epsilon\,\|\theta_h(t,\cdot)\|^2
+a_{\star}\,|\theta_h(t,\cdot)|_1^2
-2|\delta(t)|\,\|\theta_h(t,\cdot)\|\,|\theta_h(t,\cdot)|_1 \\
\ge&\,\tfrac{a_{\star}}{2}\,|\theta_h(t,\cdot)|_1^2
+\left(\epsilon-\tfrac{2\,|\delta(t)|^2}{a_{\star}}\right)
\,\|\theta_h(t,\cdot)\|^2\quad\forall\,t\in[0,T].\\
\end{split}
\end{equation}
If $\epsilon_0:=\tfrac{a_{\star}}{2}+\tfrac{2}{a_{\star}}
\,\max_{\ssy[0,T]}|\delta|^2$, \eqref{TNT72} yields that
\begin{equation}\label{TNT73}
\nu_{\epsilon_0}(t)\ge\,\tfrac{a_{\star}}{2}
\,\|\theta_h(t,\cdot)\|_1^2
\quad\forall\,t\in[0,T].
\end{equation}
\par
Now, setting $\epsilon=\epsilon_0$ in \eqref{TNT61} and then
adding the resulting equation with \eqref{TNT8}, we obtain
\begin{equation}\label{TNT74}
\tfrac{d}{dt}\nu_{\epsilon_0}(t)\leq\,C\,\nu_{\epsilon_0}(t)
+(\epsilon_0+1)\,h^{2(r+1)}\,\Gamma_{\ssy D}(t)
\quad\forall\,t\in[0,T].
\end{equation}
Since $\theta_h(0,\cdot)=0$, the bound \eqref{TNT41} follows from
\eqref{TNT74} via Gr{\"o}nwall's lemma and \eqref{TNT73}. $\Box$
\end{proof}
%
%
%
\par
A simple consequence of \eqref{4.39} and \eqref{TNT41} is the
following optimal-order error estimate.
%
%
%
\begin{thm}\label{TNTConv1}
Let $u$ be the solution of \eqref{TNT1} and $u_h$ its
semidiscrete approximation defined by \eqref{TNT2}-\eqref{TNT3}.
If $\varepsilon(t)\leq0$ for $t\in[0,T]$, then
\begin{equation}\label{TNTErrorEstim}
\|u_h(t,\cdot)-u(t,\cdot)\| +h\,\|u_h(t,\cdot)-u(t,\cdot)\|_1
\leq\,C\,h^{r+1}\,\left[\|u(t,\cdot)\|_{r+1}+
\left(\int_0^t\Gamma_{\ssy D}(\tau)\;d\tau
\,\right)^{\frac{1}{2}}\right]\quad
\forall\,t\in[0,T],
\end{equation}
where $\Gamma_{\ssy D}$ is the function defined in the statement of
Proposition~\ref{P_TNT1}. $\Box$
\end{thm}
%
%
%
%
%
%
\subsubsection{Crank-Nicolson fully discrete approximations}
%
%
We use the notation of paragraph 2.2.2. For $n=0,\dots,N$, the
Crank-Nicolson method for the problem \eqref{TNT1} yields an
approximation $U_h^n\in S_h$ of $u(t^n,\cdot)$ as follows:
\par
{\tt Step 1}: Set
\begin{equation}\label{CNPRB1}
U^0_h:=R_hu_0.
\end{equation}
\par
{\tt Step 2}: For $n=1,\dots,N$, find $U^n_h\in S_h$ such that
\begin{equation}\label{CNPRB2}
\begin{split}
(\partial U_h^n,\chi)=&\,\left[\,\varepsilon^{n-\frac{1}{2}}\,
\partial U_h^n(1)+\delta^{n-\frac{1}{2}}\,{\mathcal A}U_h^n(1)
+g^{n-\frac{1}{2}}\,\right]\,\chi(1)\\
&-a^{n-\frac{1}{2}}\,{\mathcal B}\big({\mathcal A}U_h^{n},\chi\big)
+\big(\beta^{n-\frac{1}{2}}\,{\mathcal A}U_h^n,\chi\big)
+\big(f^{n-\frac{1}{2}},\chi\big)
\quad\forall\,\chi\in S_h,\\
\end{split}
\end{equation}
where $a^{n-\frac{1}{2}}:=a(t^{n-\frac{1}{2}})$,
$\delta^{n-\frac{1}{2}}:=\delta(t^{n-\frac{1}{2}})$,
$\varepsilon^{n-\frac{1}{2}}:=\varepsilon(t^{n-\frac{1}{2}})$
$g^{n-\frac{1}{2}}:=g(t^{n-\frac{1}{2}})$,
$f^{n-\frac{1}{2}}:=f(t^{n-\frac{1}{2}},\cdot)$ and
$\beta^{n-\frac{1}{2}}:=\beta(t^{n-\frac{1}{2}},\cdot)$.
%
%
%
%
\begin{proposition}\label{CNPRB_Existence}
Let $n\in\{1,\dots,N\}$ and suppose that $U_h^{n-1}\in S_h$ is
well defined. If $\varepsilon^{n-\frac{1}{2}}\leq0$, then, there
exists a constant $C_n$ such that if $k_n<C_n$, then $U_h^n$ is
well defined by \eqref{CNPRB2}.
\end{proposition}
%
%
%
%
%
%
%
\begin{proof}
It is enough to show that if there is a $V\in S_h$ such that
\begin{equation}\label{CNPRB3}
\tfrac{1}{k_n}\,(V,\phi)=\left[\,
\tfrac{\varepsilon^{n-\frac{1}{2}}}{k_n}\,V(1)
+\tfrac{\delta^{n-\frac{1}{2}}}{2}\,V(1)\,\right]
\,\phi(1)
-\tfrac{a^{n-\frac{1}{2}}}{2}
\,{\mathcal B}\big(V,\phi\big)
+\tfrac{1}{2}\,(\beta^{n-\frac{1}{2}}\,V,\phi\big)
\quad\forall\,\phi\in S_h,
\end{equation}
then $V=0$. To arrive at the desired conclusion, first set
$\phi=V$ in \eqref{CNPRB3} and use \eqref{TraceIneq} to obtain
\begin{equation*}
\|V\|^2+|\varepsilon^{n-\frac{1}{2}}|\,|V(1)|^2\leq
\,\tfrac{k_n}{2}\,\left(\,2\,|\delta^{n-\frac{1}{2}}|
\,\|V\|\,|V|_1-a^{n-\frac{1}{2}}\,|V|_1^2
+|\beta^{n-\frac{1}{2}}|_{\infty}\,\|V\|^2\,\right).
\end{equation*}
Then use the arithmetic-geometric mean
inequality, to get
\begin{equation*}
\|V\|^2\,(1-k_n\,\gamma_n)\leq0,
\end{equation*}
where $\gamma_n=
\tfrac{1}{2}\,\left(|\beta^{n-\frac{1}{2}}|_{\infty}
+\tfrac{|\delta^{n-\frac{1}{2}}|^2}{a^{n-\frac{1}{2}}}\right)$.
This yields $V=0$ if we require, for example,
$k_n<\tfrac{1}{1+\gamma_n}$. $\Box$
\end{proof}
%
%
%
%
%
\par
The following consistency result is analogous to that of
Proposition~\ref{2.8}.
%
%
\begin{proposition}
Let $u$ be the solution of \eqref{TNT1}. For $n=1,\dots,N$,
define $\sigma^n:\overline{D}\rightarrow\mathbb{R}$ by
\begin{equation}\label{CNPRB_ConsistencyEquations}
\tfrac{u^{n}-u^{n-1}}{k_n}=a^{n-\frac{1}{2}}\,
u_{xx}(t^{n-\frac{1}{2}},\cdot)
+\beta^{n-\frac{1}{2}}\,{\mathcal A}u^n
+f^{n-\frac{1}{2}}+\sigma^n.
\end{equation}
Then
\begin{equation}\label{CNPRB_Synepeia1}
\|\sigma^n\|\leq\,C\,(k_n)^2\,\,
\left(\,\max_{\ssy[t^{n-1},t^n]}\|\partial_t^2u\|
+\max_{\ssy[t^{n-1},t^n]}\|\partial_t^3u\| \,\right),\quad
n=1,\dots,N.\,\,\Box
\end{equation}
\end{proposition}
%
%
%
%
\begin{proposition}\label{CNPRBMagnaProposition}
Let $u$ be the solution of \eqref{TNT1} and
$(U_h^n)_{n=0}^{\ssy N}$ be the fully discrete approximations
that the method \eqref{CNPRB1}-\eqref{CNPRB2} produces.
Assume that $\varepsilon(t)\leq0$ for $t\in[0,T]$, and that
\eqref{MeshCondition} holds.
Then, there exists a constant $C_{\ssy D}\ge0$ such that: if
$\displaystyle{\max_{1\leq{n}\leq{\ssy
N}}}(k_n\,C_{\ssy D})\leq\tfrac{1}{3}$,
there exists a constant $C>0$ such that
\begin{equation}\label{CNPRB_MagnaEst}
\max_{1\leq{n}\leq{\ssy N}}\|U_h^n-R_hu^n\|_1\leq\,C
\,(k^2+h^{r+1})\,\,\Upsilon_{\ssy D}(u)
\quad\forall\,h\in(0,h_{\star}],
\end{equation}
where
$\Upsilon_{\ssy D}(u):=\sum_{\ell=0}^1\max_{\ssy[0,T]}\|\partial_t^{\ell}u\|_{r+1}
+\sum_{\ell=2}^3\max_{\ssy[0,T]}\|\partial_t^{\ell}u\|_1
+\sum_{\ell=2}^4\max_{\ssy t\in[0,T]}|\partial_t^{\ell}u(t,1)|$.
\end{proposition}
%
%
\begin{proof}
Let $h\in(0,h_{\star}]$, $\delta^n:=\delta(t^n)$,
$a^n:=a(t^n)$ and
$\theta_h^n:=U_h^n-R_hu^n$ for $n=0,\dots,N$.
We use \eqref{CNPRB2}, \eqref{CNPRB_ConsistencyEquations},
\eqref{R_def}, and \eqref{4.40}, to obtain
\begin{equation}\label{CNPRBTOGOO1}
\begin{split}
(\partial\theta_h^n,\chi)=&\,
\Big[\varepsilon^{n-\frac{1}{2}}\,\partial\theta_h^n(1)
+\delta^{n-\frac{1}{2}}\,{\mathcal A}\theta_h^{n}(1)
-{\mathcal Z}_3^n\,\Big]\,\chi(1)
-a^{n-\frac{1}{2}}\,{\mathcal B}({\mathcal A}\theta_h^n,\chi)
+(\beta^{n-\frac{1}{2}}\,{\mathcal A}\theta_h^n,\chi)\\
&+({\mathcal
Z}_1^n-\sigma^n,\chi)+a^{n-\frac{1}{2}}\,{\mathcal
B}({\mathcal Z}_2^n,\chi)\quad\forall\,\chi\in S_h,
\quad n=1,\dots,N,\\
\end{split}
\end{equation}
where
${\mathcal Z}_1^n:=\partial u^n-R_h(\partial
u^n)-P_h[\,\beta^{n-\frac{1}{2}}({\mathcal
A}u^n-R_h({\mathcal A}u^n))\,\big]$,
${\mathcal Z}_2^n:=u(t^{n-\frac{1}{2}})-{\mathcal A}u^n$
and
${\mathcal Z}_3^n:=\varepsilon^{n-\frac{1}{2}}
\,\big[\,\partial_tu(t^{n-\frac{1}{2}},1)-\partial
u^n(1)\big]+\delta^{n-\frac{1}{2}}\,\big[\,u(t^{n-\frac{1}{2}},1)
-{\mathcal A}u^n(1)\,\big]$.
Using Taylor's formula and \eqref{4.39}, we deduce the following
estimates:
\begin{gather}
\|{\mathcal Z}_1^n\|
\leq\,C\,h^{r+1}\,\left[\,\max_{\ssy[t^{n-1},t^n]}\|u\|_{r+1}
+\max_{\ssy[t^{n-1},t^n]}\|\partial_tu\|_{r+1}
\,\right],\label{CNPRBFragmaE1}\\
|{\mathcal Z}_2^n|_1\leq\,C\,k_n^2
\,\max_{\ssy[t^{n-1},t^n]}|\partial_t^2u|_1,
\label{CNPRBFragmaE3}\\
|{\mathcal Z}_3^n|\leq\,C\,k_n^2\,
\left[\,\max_{\ssy t\in[t^{n-1},t^n]}|\partial_t^3u(t,1)|
+\max_{\ssy t\in[t^{n-1},t^n]}|\partial_t^2u(t,1)|\,\right],
\label{CNPRBFragmaE4}
\end{gather}
for $n=1,\dots,N$. The proof now proceeds in four steps.
%
%
\par\noindent\smallskip
\par
{\tt Step I}: Set $\chi={\mathcal A}\theta_h^n$ in \eqref{CNPRBTOGOO1}
and use \eqref{TraceIneq} and the
Cauchy-Schwarz inequality, to obtain
\begin{equation*}
\begin{split}
\|\theta_h^n\|^2+|\varepsilon^n|\,|\theta_h^n(1)|^2
\leq&\,\|\theta_h^{n-1}\|^2
+|\varepsilon^{n-1}|\,|\theta_h^{n-1}(1)|^2\\
&+(\varepsilon^{n-\frac{1}{2}}-\varepsilon^n)
\,|\theta_h^n(1)|^2
+(\varepsilon^{n-1}-\varepsilon^{n-\frac{1}{2}})
\,|\theta_h^{n-1}(1)|^2\\
&
+2\,k_n\,\Big[\,
2\,|\delta^{n-\frac{1}{2}}|\,\|{\mathcal A}\theta_h^n\|
\,|{\mathcal A}\theta_h^n|_1-a_{\star}\,|{\mathcal A}\theta_h^n|_1^2
+|\beta^{n-\frac{1}{2}}|_{\infty}\,\|{\mathcal A}\theta_h^n\|^2\\
&\hskip1.5truecm
+\sqrt{2}\,|{\mathcal Z}_3^n|\,\|{\mathcal A}\theta_h^n\|^{\frac{1}{2}}
\,|{\mathcal A}\theta_h^n|_1^{\frac{1}{2}}\\
&\hskip1.5truecm+\left(\,\|\sigma^n\|+\|{\mathcal Z}_1^n\|\right)
\,\|{\mathcal A}\theta_h^n\|+a^{n-\frac{1}{2}}\,|{\mathcal Z}_2^n|_1
\,|{\mathcal A}\theta_h^n|_1\,\Big],\quad n=1,\dots,N,\\
\end{split}
\end{equation*}
which, after the use of the arithmetic-geometric mean
inequality, yields
\begin{equation}\label{CNPRBMoon1}
\begin{split}
\|\theta_h^n\|^2+|\varepsilon^n|\,|\theta_h^n(1)|^2
\leq&\,\|\theta_h^{n-1}\|^2
+|\varepsilon^{n-1}|\,|\theta_h^{n-1}(1)|^2
+C\,k_n\,\left(\,|\theta_h^n(1)|^2
+|\theta_h^{n-1}(1)|^2\,\right)\\
&+C\,k_n\,\big(\,\|\theta_h^n\|^2+\|\theta_h^{n-1}\|^2\\
&\hskip1.5truecm
+\|\sigma^n\|^2+\|{\mathcal Z}_1^n\|^2
+|{\mathcal Z}_2^n|_1^2
+|{\mathcal Z}_3^n|^2\,\big),
\quad n=1,\dots,N.\\
\end{split}
\end{equation}
%
%
%
\par\noindent\smallskip
\par
{\tt Step II}: Now set $\chi=\partial\theta_h^n$ in
\eqref{CNPRBTOGOO1} to get
\begin{equation*}
\begin{split}
|\theta_h^n|_1^2-\tfrac{\delta^n}{a^n}\,|\theta_h^n(1)|^2
\leq&\,|\theta_h^{n-1}|_1^2-\tfrac{\delta^{n-1}}{a^{n-1}}
\,|\theta_h^{n-1}(1)|^2
-2\,k_n\,\tfrac{|\varepsilon^{n-\frac{1}{2}}|}{a^{n-\frac{1}{2}}}
\,|\partial\theta_h^n(1)|^2\\
&+2\,\tfrac{k_n}{a^{n-\frac{1}{2}}}\,\Big[-\|\partial\theta_h^n\|^2
+|\beta^{n-\frac{1}{2}}|_{\infty}\,\|{\mathcal
A}\theta_h^n\|\,\|\partial\theta_h^n\|
+a^{n-\frac{1}{2}}\,{\mathcal B}({\mathcal Z}_2^n,\partial\theta_h^n)\\
&\hskip2.0truecm
+\left(\|{\mathcal Z}_1^n\|+\|\sigma^n\|\right)\,\|\partial\theta_h^n\|
-{\mathcal Z}_3^n\,\partial\theta_h^n(1)\,\Big]\\
&+\left(\tfrac{\delta^{n-\frac{1}{2}}}{a^{n-\frac{1}{2}}}
-\tfrac{\delta^n}{a^n}\right)
\,|\theta_h^n(1)|^2
+\left(-\tfrac{\delta^{n-\frac{1}{2}}}{a^{n-\frac{1}{2}}}
+\tfrac{\delta^{n-1}}{a^{n-1}}\right)
\,|\theta_h^{n-1}(1)|^2,\quad n=1,\dots,N,\\
\end{split}
\end{equation*}
which, after the use of the arithmetic-geometric mean
inequality, yields
\begin{equation}\label{CNPRBMoon2}
\begin{split}
|\theta_h^n|_1^2-\tfrac{\delta^n}{a^n}\,|\theta_h^n(1)|^2
\leq&\,|\theta_h^{n-1}|_1^2-\tfrac{\delta^{n-1}}{a^{n-1}}
\,|\theta_h^{n-1}(1)|^2
+C\,k_n\,\left(\,|\theta_h^n(1)|^2
+|\theta_h^{n-1}(1)|^2\,\right)\\
&+2\,k_n\,\left[
\,{\mathcal B}({\mathcal Z}_2^n,\partial\theta_h^n)
-\tfrac{1}{a^{n-\frac{1}{2}}}
\,{\mathcal Z}_3^n\,\partial\theta_h^n(1)\right]\\
&+C\,k_n\,\left(\,\|\theta_h^n\|^2+\|\theta_h^{n-1}\|^2
+\|{\mathcal Z}_1^n\|^2+\|\sigma^n\|^2\,\right),
\quad n=1,\dots,N.\\
\end{split}
\end{equation}
%
%
\par\noindent\smallskip
\par
{\tt Step III}: For $\rho>0$, we introduce the quantities
\begin{equation}\label{Vouvalos}
{\mathcal V}_{\rho}^m:=\rho\left(\,\|\theta_h^m\|^2
+|\varepsilon^m|\,|\theta_h^m(1)|^2\,\right)
+|\theta_h^m|_1^2-\tfrac{\delta^m}{a^m}\,|\theta_h^m(1)|^2,
\quad m=0,\dots,N.
\end{equation}
Now, using \eqref{Vouvalos} and \eqref{TraceIneq} we have
\begin{equation*}
\begin{split}
{\mathcal V}_{\rho}^m\ge&\,\rho\,\|\theta_h^m\|^2+|\theta_h^m|_1^2
-2\,\tfrac{|\delta^m|}{a^m}\,|\theta_h^m|_1\,\|\theta_h^m\|\\
\ge&\,\tfrac{1}{2}\,|\theta_h^m|_1^2+\|\theta_h^m\|^2
\,\left(\,\rho-2\,\tfrac{(\delta^m)^2}{(a^m)^2}\right),
\quad m=0,\dots,N.\\
\end{split}
\end{equation*}
Thus, choosing $\rho=\rho_0:=\tfrac{1}{2}
+2\max_{\ssy[0,T]}\tfrac{\delta^2}{a^2}$, we obtain
\begin{equation}\label{Vouvalos12}
{\mathcal V}_{\rho_0}^m\ge\,\tfrac{1}{2}
\,\|\theta_h^m\|_1^2,\quad m=0,\dots,N.
\end{equation}
%

%
%
\par\noindent\smallskip
\par
{\tt Step IV}: Combining \eqref{CNPRBMoon1}, \eqref{CNPRBMoon2},
\eqref{TraceIneq}, \eqref{CNPRBFragmaE1}, \eqref{CNPRBFragmaE3},
\eqref{CNPRBFragmaE4}, and \eqref{CNPRB_Synepeia1} we obtain
\begin{equation}\label{Vouvalos13}
\begin{split}
{\mathcal V}^n_{\rho_0}\leq&\,{\mathcal V}_{\rho_0}^{n-1}
+C\,k_n\,\left(\,\|\theta_h^n\|_1^2+\|\theta_h^{n-1}\|_1^2\,\right)
+C\,k_n\,\left[\,(k_n)^2+(h^{r+1})^2\,\right]\,(\Upsilon_1(u))^2\\
&+2\,k_n\,\left[
\,{\mathcal B}({\mathcal Z}_2^n,\partial\theta_h^n)
-\tfrac{1}{a^{n-\frac{1}{2}}}
\,{\mathcal Z}_3^n\,\partial\theta_h^n(1)\right],\quad n=1,\dots,N,\\
\end{split}
\end{equation}
where
$\Upsilon_1(u):=
\max_{\ssy[0,T]}\|u\|_{r+1}
+\max_{\ssy[0,T]}\|\partial_tu\|_{r+1}
+\max_{\ssy[0,T]}\|\partial_t^2u\|_1
+\max_{\ssy[0,T]}\|\partial_t^3u\|
+\max_{\ssy t\in[0,T]}|\partial_t^2u(t,1)|
+\max_{\ssy t\in[0,T]}|\partial_t^3u(t,1)|$.
Using \eqref{Vouvalos13} and \eqref{Vouvalos12} we conclude that
there exist constants $C_1\ge0$ and $C_2\ge 0$, such that
\begin{equation}\label{Vouvalos14}
\begin{split}
(1-C_1\,k_n)\,{\mathcal V}_{\rho_0}^n
\leq&\,(1+C_1\,k_n)\,{\mathcal V}^{n-1}_{\rho_0}
+C_2\,k_n\,(h^{r+1}+k_n^2)^2\,(\Upsilon_1(u))^2\\
&+2\,k_n\,\left[
\,{\mathcal B}({\mathcal Z}_2^n,\partial\theta_h^n)
-\tfrac{1}{a^{n-\frac{1}{2}}}
\,{\mathcal Z}_3^n\,\partial\theta_h^n(1)\right],
\quad n=1,\dots,N.\\
\end{split}
\end{equation}
To continue, we assume that
$\displaystyle{\max_{1\leq{n}\leq{\ssy N}}}(C_1\,k_n)
\leq\,\tfrac{1}{3}$, which allows us to conclude that
$\tfrac{1+C_1\,k_n}{1-C_1\,k_n}\leq\,e^{3C_1k_n}$ for
$n=1,\dots,N$. Hence, \eqref{Vouvalos14} yields
\begin{equation*}
\begin{split}
{\mathcal V}_{\rho_0}^n\leq&\,e^{3C_1k_n}
\,{\mathcal V}_{\rho_0}^{n-1}
+\tfrac{C_2\,k_n}{1-C_1\,k_n}
\,(h^{r+1}+k_n^2)^2\,(\Upsilon_1(u))^2\\
&+\tfrac{2\,k_n}{1-C_1\,k_n}\,\Big[
\,{\mathcal B}({\mathcal Z}_2^n,\partial\theta_h^n)
-\tfrac{1}{a^{n-\frac{1}{2}}}\,{\mathcal Z}_3^n
\,\partial\theta_h^n(1)\,\Big],\quad n=1,\dots,N.\\
\end{split}
\end{equation*}
\par
Letting
$\lambda_j^n:=\tfrac{\exp\left(3C_1
\sum_{\ell=j+1}^{n}k_{\ell}\right)}{1-C_1\,k_j}$
and using a simple induction argument we arrive at
\begin{equation*}
\begin{split}
{\mathcal V}_{\rho_0}^n\leq&\,C_2\,(\Upsilon_1(u))^2
\,\sum_{j=1}^{n}k_j
\,\lambda_j^n\,(h^{r+1}+k_j^2)^2\\
&+2\,\sum_{j=1}^nk_j \,\lambda_j^n\,\Big[
\,{\mathcal B}({\mathcal Z}_2^j,\partial\theta_h^j)]
-\tfrac{1}{a^{j-\frac{1}{2}}}\,{\mathcal Z}_3^j
\,\partial\theta_h^j(1)\,\Big],\quad n=1,\dots,N,\\
\end{split}
\end{equation*}
which yields
\begin{equation}\label{Vouvalos15}
{\mathcal V}_{\rho_0}^n\leq\,C\,(h^{r+1}+k^2)^2
\,(\Upsilon_1(u))^2
+{\mathcal T}_1^n+{\mathcal T}_2^n,\quad n=1,\dots,N,
\end{equation}
where
${\mathcal T}_1^n:=2\,\sum_{j=1}^n\lambda_j^n\,
{\mathcal B}({\mathcal Z}_2^j,\theta_h^j-\theta_h^{j-1})$
and
${\mathcal T}_2^n:=-2\,\sum_{j=1}^n
\tfrac{\lambda_j^n}{a^{j-\frac{1}{2}}}\,{\mathcal Z}_3^j
\,(\theta_h^j(1)-\theta_h^{j-1}(1))$.
First, we proceed as in bounding the quantity $T_{\ssy A}^n$
in the proof of Proposition~\ref{MagnaProposition} to get
\begin{equation}\label{Vouvalos16}
|{\mathcal T}_1^n|\leq\,C\,k^2\,\Upsilon_2(u)
\,\max_{1\leq{m}\leq{n}}|\theta_h^m|_1,
\quad n=1,\dots,N,
\end{equation}
where
$\Upsilon_2(u):=\max_{\ssy [0,T]}|\partial_t^2u|_1
+\max_{\ssy [0,T]}|\partial_t^3u|_1$.
In addition, we have
\begin{equation}\label{Vouvalos17}
\begin{split}
-{\mathcal T}_2^n=&\tfrac{2}{1-C_1\,k_n}\,\tfrac{1}{a^{n-\frac{1}{2}}}
\,{\mathcal Z}_3^n\,\theta_h^n(1)
+2\,\sum_{j=1}^{n-1}\tfrac{\lambda_j^n}{a^{j-\frac{1}{2}}}
\,\left({\mathcal Z}_3^j-{\mathcal Z}_3^{j+1}\right)\,\theta_h^j(1)\\
&+\,2\,\sum_{j=1}^{n-1}
\tfrac{\exp\left(3\,C_1\sum_{\ell=j+2}^{n}k_{\ell}\right)}
{a^{j-\frac{1}{2}}}
\,\left[\tfrac{\exp\left(3\,C_1\,k_{j+1}\right)-1+C_1\,k_j}{1-C_1\,k_j}
-\tfrac{C_1\,k_{j+1}}{1-C_1\,k_{j+1}}\right]
\,{\mathcal Z}_3^{j+1}\,\theta_h^j(1)\\
&+\,2\,\sum_{j=1}^{n-1}\left(\,\tfrac{1}{a^{j-\frac{1}{2}}}
-\tfrac{1}{a^{j+\frac{1}{2}}}\right)\,
\lambda_{j+1}^n
\,{\mathcal Z}_3^{j+1}\,\theta_h^j(1),\quad n=1,\dots,N.\\
\end{split}
\end{equation}
Observing that
\begin{equation*}
|{\mathcal Z}_3^j-{\mathcal Z}_3^{j+1}|
\leq\,C\,(k_j+k_{j+1})\,\left[\,(k_j)^2+|k_j-k_{j+1}|\,\right]
\,\Upsilon_3(u),\quad j=1,\dots,N-1,
\end{equation*}
with
$\Upsilon_3(u):=\sum_{\ell=2}^4
\max_{\ssy t\in[0,T]}|\partial_t^{\ell}u(t,1)|$,
we see that \eqref{Vouvalos17}, \eqref{MeshCondition},
\eqref{CNPRBFragmaE4} and \eqref{TraceIneq} yield

\begin{equation}\label{Vouvalos18}
|{\mathcal T}_2^n|\leq\,C\,k^2\,\Upsilon_3(u)
\,\max_{1\leq{m}\leq{n}}\|\theta_h^m\|_1,
\quad n=1,\dots,N.
\end{equation}
Now, from \eqref{Vouvalos15}, \eqref{Vouvalos16} and \eqref{Vouvalos18}
there follows that
\begin{equation*}
{\mathcal V}_{\rho_0}^n\leq\,C\,(h^{r+1}+k^2)^2\,(\Upsilon_1(u))^2
+C\,k^2\,\big(\,\Upsilon_2(u)+\Upsilon_3(u)\,\big)
\,\max_{1\leq{m}\leq{n}}\|\theta_h^m\|_1,
\quad n=1,\dots,N.
\end{equation*}
Use then \eqref{Vouvalos12} to arrive at
\begin{equation*}
\max_{0,\leq{n}\leq{\ssy N}}\|\theta_h^n\|_1^2
\leq\,C\,(h^{r+1}+k^2)^2\,(\Upsilon_1(u)
+\Upsilon_2(u)+\Upsilon_3(u))^2,
\end{equation*}
which is the desired estimate \eqref{CNPRB_MagnaEst}. $\Box$
\end{proof}
%
%
\par
As a simple consequence of \eqref{CNPRB_MagnaEst} and \eqref{4.39}
we obtain the following optimal-order error estimates in $L^2$ and $H^1$
norms.
%
%
\begin{thm}
Let $u$ be the solution of \eqref{TNT1} and $(U_h^n)_{n=0}^{\ssy
N}$ be the fully discrete approximations that the method
\eqref{CNPRB1}-\eqref{CNPRB2} produces. Assume that
$\varepsilon(t)\leq0$ for $t\in[0,T]$, that \eqref{MeshCondition}
holds and $\displaystyle{\max_{1\leq{n}\leq{\ssy N}}}k_n\,C_{\ssy
D}\leq\tfrac{1}{3}$, where $C_{\ssy D}$ is the constant specified
in Proposition~\ref{CNPRBMagnaProposition}. Then
\begin{gather*}
\max_{0\leq{n}\leq{\ssy N}}\|U_h^n-u^n\|_{\ell}\leq\,C\,(k^2+h^{r+1-\ell})
\,\,\Upsilon_{\ssy D}(u),\quad\forall\,h\in(0,h_{\star}],
\end{gather*}
for $\ell=0,1$, where $\Upsilon_{\ssy D}(u)$ was specified in
Proposition~\ref{CNPRBMagnaProposition}. $\Box$
\end{thm}
%
%
%
%
%
\subsection{The reactive case}
In this paragraph, we propose finite element approximations when
the dynamical boundary condition in \eqref{TNT1} is of reactive
type, i.e. $\varepsilon(t)>0$ for $t\in[0,T]$. According to
\cite{VazVit2008}, \cite{BBR}, the problem is well posed only in
the one-dimensional case. To construct a finite element method
for this problem we follow the idea (cf. paragraph 3.3) to
replace the term $u_t$ in the dynamical boundary condition using
the partial differential equation in \eqref{TNT1}. Hence we
obtain: $a(t)\,u_{xx}(t,1)=\tfrac{a(t)}{\varepsilon(t)}u_x(1,t)
-\left[\tfrac{\delta(t)}{\varepsilon(t)}+\beta(t,1)\right]\,u(t,1)
-\left[\tfrac{g(t)}{\varepsilon(t)}+f(t,1)\right]$ for
$t\in[0,T]$. Then, to use this as a boundary condition, we
formulate a variational formulation using ${\mathcal
B}(\cdot,\cdot)$ instead of the $L^2(D)$ inner product
$(\cdot,\cdot)$. Of course this approach works also if
$\varepsilon(t)<0$ for $t\in[0,T]$.
%
%
\subsubsection{Preliminaries}
Let $r\in{\mathbb N}$ with $r\ge3$, and ${\mathcal S}_h$ be a
finite-dimensional subspace of ${\mathbb H}^2(D)$ consisting of
$C^1$ functions that are polynomials of degree less or equal to
$r$ in each interval of a non-uniform partition of $D$ with
maximum length $h\in(0,h_{\star}]$. It is well-known,
\cite{BrHilbert1970}, that the following approximation property
holds:
\begin{equation}\label{H2AppProp}
\inf_{\chi\in{\mathcal S}_h}\|v-\chi\|_2\leq\,
C\,h^{s-1}\,\|v\|_{s+1},
\quad\forall\,v\in{\mathbb H}^{s+1}(D),
\,\,\,\forall\,h\in(0,h_{\star}],\quad s=1,\dots,r.
\end{equation}
\par
We introduce bilinear forms
${\mathcal B}^{\star}$, $\gamma^{\star}:
H^2(D)\times H^2(D)\rightarrow{\mathbb R}$
given by ${\mathcal B}^{\star}(v,w):=(v'',w'')$
and $\gamma^{\star}(v,w):=(v'',w'')+(v',w')$ for $v$ and
$w\in H^2(D)$, and set $|v|_2:=\|v''\|$ for $v\in H^2(D)$.
Also, we define a new elliptic projection
$R_h^{\star}:H^2(D)\rightarrow{\mathcal S}_h$ by
\begin{equation}\label{RDEF2}
\gamma^{\star}(R_h^{\star}v,w)
=\gamma^{\star}(v,\chi)
\quad\forall\,\chi\in{\mathcal S}_h.
\end{equation}
%
%
\begin{lem}
The elliptic projection $R_h^{\star}$ has the following
property
\begin{equation}\label{Reactive102}
(R_h^{\star}v)'(1)=v'(1)+(R_h^{\star}v-v)(1)
-\tfrac{1}{6}\,
{\mathcal B}(R_h^{\star}v-v,\omega)
\quad\forall\,v\in{\mathbb H}^2(D),
\end{equation}
where $\omega(x)=x^3$.
\end{lem}
%
%
%
%
%
\begin{proof}
Let $v\in{\mathbb H}^2(D)$ and $\rho=R_h^{\star}v-v$. Since
$\omega\in{\mathcal S}_h$, setting $\chi=\omega$ in \eqref{RDEF2}
we obtain $\int_{\ssy
D}\rho''(x)\,x\;dx=-\tfrac{1}{6}\,(\rho',\omega')$. Then,
integrating by parts we get
$\rho'(1)=\rho(1)-\rho(0)-\tfrac{1}{6}\,(\rho',\omega')$, which
is the desired equality, since $\rho(0)=0$. $\Box$
\end{proof}
%
%
%
%
\begin{proposition}\label{Le_Propo_chemin}
The elliptic projection $R_h^{\star}$ has the following
approximation properties:
\begin{equation}\label{StarEll1}
\sum_{\ell=1}^2h^{\ell}\,\|R_h^{\star}v-v\|_{\ell}
\leq\,C\,h^{s+1}\,\|v\|_{s+1}
\end{equation}
and
\begin{equation}\label{StarEll4}
|(R_h^{\star}v-v)'(1)|+|(R_h^{\star}v-v)(1)|
\leq\,C\,h^s\,\|v\|_{s+1}
\end{equation}
for $s=1,\dots,r$, $v\in{\mathbb H}^{s+1}(D)$
and $h\in(0,h_{\star}]$.
\end{proposition}
%
%
%
%
%
%
\begin{proof}
Let $h\in(0,h_{\star}]$, $s\in\{1,\dots,r\}$,
$v\in{\mathbb H}^{s+1}(D)$ and $e=R_h^{\star}v-v$.
Using \eqref{RDEF2} we have $\gamma^{\star}(e,e)=
\gamma^{\star}(e,\chi-v)$ for $\chi\in{\mathcal S}_h$, which
along with \eqref{H2AppProp} yields
\begin{equation}\label{Reactive104}
|e|_2+|e|_1\leq\,C\,h^{s-1}\,\|v\|_{s+1}.
\end{equation}
Now, let $w\in H^3(D)$
such that
\begin{equation}\label{ReactiveCorr1}
\gathered
-w'''+w'=e'\quad\text{\rm in}\ \  D,\\
w(0)=w''(1)=w''(0)=0.\\
\endgathered
\end{equation}
It is easily seen that \eqref{ReactiveCorr1} conceals a standard
two-point boundary-value problem with respect to $w'$ and thus
existence and uniqueness of its solution follows in a
straightforward way; in addition we have that
\begin{equation}\label{ReactiveCorr2}
\|w\|_3\leq\,C\,|e|_1.
\end{equation}
Thus, we obtain $\|e'\|^2=\gamma^{\star}(e,w-\chi)$ for
$\chi\in{\mathcal S}_h$. Then, we use \eqref{Reactive104},
\eqref{H2AppProp} and \eqref{ReactiveCorr2} to get
\begin{equation*}
\begin{split}
|e|_1^2\leq&\,C\,\left(|e|_2+|e|_1\right)\,h^1\,\|w\|_3\\
\leq&\,C\,h^s\,\|v\|_{s+1}\,|e|_1,
\end{split}
\end{equation*}
which yields
\begin{equation}\label{Reactive105}
|e|_1\leq C\,h^{s}\,\|v\|_{s+1}.
\end{equation}
Hence, \eqref{StarEll1} follows as a simple consequence
of \eqref{Reactive104} and \eqref{Reactive105}.
\par
Using \eqref{Reactive102}, \eqref{TraceIneq}, \eqref{PoincareF},
and \eqref{StarEll1} we have
\begin{equation*}
\begin{split}
|e'(1)|^2+|e(1)|^2\leq&\,C\,\left(\,|e(1)|^2+\|e'\|^2\,\right)\\
\leq&\,C\,|e|_1^2\\
\leq&\,C\,h^{2s}\,\|v\|_{s+1}^2,\\
\end{split}
\end{equation*}
which obviously yields \eqref{StarEll4}. $\Box$
\end{proof}
%
%
%
\par
For later use, we close this section by extending
\eqref{TraceIneq} as follows:
\begin{lem}
For $v\in H^2(D)$ it holds that
\begin{equation}\label{TraceIneq2}
|v'(1)|^2\leq\,|v|_1^2+2\,|v|_1\,|v|_2.
\end{equation}
\end{lem}
%
%
\begin{proof}
Let $v\in H^2(D)$. Observing that
$|v'(1)|^2=\int_{\ssy D}[(v'(x))^2\,x]'\,dx$,
we obtain
$|v'(1)|^2=\|v'\|^2+2\int_{\ssy D}x\,v'(x)\,v''(x)\;dx$,
which yields \eqref{TraceIneq2} via the Cauchy-Schwarz inequality.
$\Box$
\end{proof}
%
%
%
%
%
\subsubsection{Semidiscrete approximation}
%
%
We define $u_h:[0,T]\rightarrow{\mathcal S}_h$, a space-discrete
approximation of $u$, requiring
\begin{equation}\label{Reactive1}
\begin{split}
{\mathcal B}(\partial_tu_h(t,\cdot),\chi)=&\,
\left\{\,\tfrac{a(t)}{\varepsilon(t)}
\,\partial_xu_h(t,1)
-\left[\,\tfrac{\delta(t)}{\varepsilon(t)}
+\beta(t,1)\,\right]\,u_h(t,1)
-\left[\,\tfrac{g(t)}{\varepsilon(t)}
+f(t,1)\right]\,\right\}\,\chi'(1)\\
&\quad+f(t,0)\,\chi'(0)
-a(t)\,{\mathcal B}^{\star}(u_h(t,\cdot),\chi)\\
&\quad+{\mathcal B}(\beta(t,\cdot)\,u_h(t,\cdot),\chi)
+{\mathcal B}(f(t,\cdot),\chi)
\quad\forall\,\chi\in{\mathcal S}_h,
\quad\forall\,t\in[0,T],\\
\end{split}
\end{equation}
and
\begin{equation}\label{Reactive2}
u_h(0,\cdot)=R^{\star}_hu_0(\cdot).
\end{equation}
%
%
%
%
\begin{proposition}
If $\varepsilon(t)>0$ for $t\in[0,T]$, then the problem
\eqref{Reactive1}-\eqref{Reactive2} admits a unique
solution $u_h\in C^1([0,T];{\mathcal S}_h)$.
\end{proposition}
%
%
%
\begin{proof}
The result follows if we argue along the lines of the proof of
Proposition~\ref{2.3.3**}. $\Box$
\end{proof}
%
%
%
In the sequel, we assume that the solution of the ibvp
\eqref{TNT1} in the reactive case is sufficiently smooth.
\begin{thm}\label{P_REACT_1}
Let $u$ be the solution of \eqref{TNT1}, $u_h$ its semidiscrete
approximation defined by \eqref{Reactive1}-\eqref{Reactive2}, and
$\Gamma_{\ssy D}$ be the function specified in
Proposition~\ref{P_TNT1}. If $\varepsilon(t)>0$ for $t\in[0,T]$,
then
\begin{equation}\label{React_ErrorEstim}
\|u_h(t,\cdot)-u(t,\cdot)\|_1
\leq\,C\,h^r\,\left[\|u(t,\cdot)\|_{r+1}+
\left(\int_0^t\Gamma_{\ssy D}(\tau)\;d\tau
\,\right)^{\frac{1}{2}}\right]
\quad\forall\,t\in[0,T],\ \ \forall\,h\in(0,h_{\star}].
\end{equation}
\end{thm}
%
%
%
%
%
%
%
\begin{proof}
Let $h\in(0,h_{\star}]$, $\theta_h:=u_h-R^{\star}_hu$,
and $\eta=R_h^{\star}u-u$. Using \eqref{Reactive1},
\eqref{RDEF2} and \eqref{TNT1}, we obtain
\begin{equation}\label{Reactive4}
\begin{split}
{\mathcal B}(\partial_t\theta_h(t,\cdot),\chi)&=\left\{\,
\tfrac{a(t)}{\varepsilon(t)}\,\partial_x\theta_h(t,1)
-\left[\tfrac{\delta(t)}{\varepsilon(t)}
+\beta(t,1)\right]\,\theta_h(t,1)+{\mathcal E}_{{\ssy R},2}(t)
\,\right\}\,\chi'(1)\\
&\quad\quad-a(t)\,{\mathcal B}^{\star}(\theta_h(t,\cdot),\chi)
+{\mathcal B}(\beta(t,\cdot)\,\theta_h(t,\cdot),\chi)\\
&\quad\quad
+{\mathcal B}\big({\mathcal E}_{{\ssy R},1}(t,\cdot),\chi\big)
+a(t)\,{\mathcal B}(R_h^{\star}u(t,\cdot)-u(t,\cdot),\chi)
\quad\forall\,\chi\in{\mathcal S}_h,\quad\forall\,t\in[0,T],\\
\end{split}
\end{equation}
where
${\mathcal E}_{{\ssy R},1}:=[\partial_tu-R^{\star}_h(\partial_tu)]
-\beta\,(u-R^{\star}_hu)$ and ${\mathcal E}_{{\ssy R},2}(t)
:=\tfrac{a(t)}{\varepsilon(t)}\,\partial_x\eta(t,1)
-\left[\tfrac{\delta(t)}{\varepsilon(t)}
+\beta(t,1)\right]\,\eta(t,1)$.
First observe that using \eqref{StarEll1}, \eqref{StarEll4}
and \eqref{PoincareF}, it follows that
\begin{equation}\label{Reactive5}
\left|\,
{\mathcal B}({\mathcal E}_{{\ssy R},1}(t,\cdot),\chi)
+a(t)\,{\mathcal B}(\eta(t,\cdot),\chi)\,\right|\leq\,C\,h^r
\,\left(\,\|u(t,\cdot)\|_{r+1}
+\|\partial_tu(t,\cdot)\|_{r+1}\,\right)
\,|\chi|_1
\end{equation}
and
\begin{equation}\label{Reactive6}
|{\mathcal E}_{{\ssy R},2}(t)\,\chi'(1)|\leq\,C\,h^{r}
\,\|u(t,\cdot)\|_{r+1}\,|\chi'(1)|
\end{equation}
for $\chi\in{\mathcal S}_h$ and $t\in[0,T]$.
Then, set $\chi=\theta_h$ in \eqref{Reactive4} and use the
Cauchy-Schwarz inequality, \eqref{PoincareF}, \eqref{Reactive5},
\eqref{Reactive6}, \eqref{TraceIneq2}, and \eqref{TraceIneq}, to
get
\begin{equation*}
\tfrac{1}{2}\,\tfrac{d}{dt}|\theta_h(t,\cdot)|_1^2\leq
\,-a_{\star}\,|\theta_h(t,\cdot)|_2^2
+C\,\Big[\,|\theta_h(t,\cdot)|_1^2
+h^{2r}\,\Gamma_{\ssy D}(t)
+|\theta_h(t,\cdot)|_1\,|\theta_h(t,\cdot)|_2\,\Big]
\quad\forall\,t\in[0,T],
\end{equation*}
which, along the arithmetic-geometric
mean inequality, yields
\begin{equation}\label{Reactive7}
\tfrac{d}{dt}|\theta_h(t,\cdot)|_1^2\leq\,
C\,\left[\,|\theta_h(t,\cdot)|_1^2
+h^{2r}\,\Gamma_{\ssy D}(t)\,\right]
\quad\forall\,t\in[0,T].
\end{equation}
Since $\theta_h(0,\cdot)=0$, using Gr{\"o}nwall's lemma from
\eqref{Reactive7} we see that
\begin{equation}\label{Reactive3}
|\theta_h(t,\cdot)|_1^2\leq\,C\,h^{2r}
\,\left(\int_0^t\Gamma_{\ssy D}(\tau)\;d\tau\right)
\quad\forall\,t\in[0,T].
\end{equation}
Finally, we combine \eqref{PoincareF}, \eqref{Reactive3} and \eqref{StarEll1}
to arrive at the error estimate \eqref{React_ErrorEstim}. $\Box$
\end{proof}
%
%
%
%
%
\subsubsection{Crank-Nicolson fully discrete approximations}
%
%
For $n=0,\dots,N$, the Crank-Nicolson method for the problem
\eqref{TNT1} yields an approximation $U_h^n\in{\mathcal S}_h$ of
$u(t^n,\cdot)$ as follows:
\par
{\tt Step 1}: Set
\begin{equation}\label{Reactive100}
U^0_h:=R_h^{\star}u_0.
\end{equation}
\par
{\tt Step 2}: For $n=1,\dots,N$, find $U^n_h\in{\mathcal S}_h$ such that
\begin{equation}\label{Reactive101}
\begin{split}
{\mathcal B}(\partial U_h^n,\chi)=&\,\left\{\,
\tfrac{a^{n-\frac{1}{2}}}{\varepsilon^{n-\frac{1}{2}}}\,
({\mathcal A}U_h^n)'(1)-\left[\,\tfrac{\delta^{n-\frac{1}{2}}}
{\varepsilon^{n-\frac{1}{2}}}+\beta^{n-\frac{1}{2}}(1)
\,\right]\,{\mathcal A}U_h^n(1)
-\left[\,\tfrac{g^{n-\frac{1}{2}}}{\varepsilon^{n-\frac{1}{2}}}
+f^{n-\frac{1}{2}}(1)\right]\,\right\}\,\chi'(1)\\
&+f^{n-\frac{1}{2}}(0)\,\chi'(0)
-a^{n-\frac{1}{2}}\,{\mathcal B}^{\star}
\big({\mathcal A}U_h^{n},\chi\big)
+{\mathcal B}\big(\beta^{n-\frac{1}{2}}\,
{\mathcal A}U_h^n,\chi\big)
+{\mathcal B}\big(f^{n-\frac{1}{2}},\chi\big)
\quad\forall\,\chi\in{\mathcal S}_h.\\
\end{split}
\end{equation}
%
%
%
%
\begin{proposition}
Let $n\in\{1,\dots,N\}$ and suppose that $U_h^{n-1}\in S_h$ is
well defined. If $\varepsilon^{n-\frac{1}{2}}>0$, then, there
exists a constant $C_n$ such that if $k_n<C_n$, then $U_h^n$ is
well defined by \eqref{Reactive101}.
\end{proposition}
%
%
%
\begin{proof}
It is enough to show that if there is a $V\in{\mathcal S}_h$ such that
\begin{equation}\label{Reactive200}
\tfrac{1}{k_n}\,{\mathcal B}(V,\phi)=\tfrac{1}{2}\,\left\{\,
\tfrac{a^{n-\frac{1}{2}}}{\varepsilon^{n-\frac{1}{2}}}\,V'(1)
-\left[\,\tfrac{\delta^{n-\frac{1}{2}}}{\varepsilon^{n-\frac{1}{2}}}
+\beta^{n-\frac{1}{2}}(1)\,\right]\,V(1)\,\right\}
\,\phi'(1)
-\tfrac{a^{n-\frac{1}{2}}}{2}
\,{\mathcal B}^{\star}\big(V,\phi\big)
+\tfrac{1}{2}\,{\mathcal B}(\beta^{n-\frac{1}{2}}\,V,\phi\big)
\end{equation}
for all $\phi\in{\mathcal S}_h$, then $V=0$. To arrive at the
desired conclusion, first set $\phi=V$ in \eqref{Reactive200} and
use \eqref{PoincareF} to obtain
\begin{equation*}
|V|_1^2\leq\tfrac{k_n}{2}\,\left[\,
\left(1+\tfrac{a^{n-\frac{1}{2}}}{\varepsilon^{n-\frac{1}{2}}}\right)\,|V'(1)|^2
+\Big|\tfrac{\delta^{n-\frac{1}{2}}}{\varepsilon^{n-\frac{1}{2}}}
+\beta^{n-\frac{1}{2}}(1)\Big|^2\,|V(1)|^2
-a^{n-\frac{1}{2}}\,|V|_2^2
+|\beta^{n-\frac{1}{2}}|_{1,\infty}\,(1+C_{\ssy P\!F})
|V|_1^2\right].
\end{equation*}
Then, use \eqref{TraceIneq}, \eqref{PoincareF}, and
\eqref{TraceIneq2}, to get
$\|V\|^2\,\big(1-\tfrac{k_n}{2}\,c_n\big)\leq0$,
where $c_{n}:=|\beta^{n-\frac{1}{2}}|_{1,\infty}\,(1+C_{\ssy
P\!F})
+2\,\Big(\,\tfrac{\delta^{n-\frac{1}{2}}}{\varepsilon^{n-\frac{1}{2}}}
+\beta^{n-\frac{1}{2}}(1)\,\Big)^2\,C_{\ssy P\!F}
+1+\tfrac{a^{n-\frac{1}{2}}}{\varepsilon^{n-\frac{1}{2}}}
+\tfrac{1}{a^{n-\frac{1}{2}}}\,\Big(
1+\tfrac{a^{n-\frac{1}{2}}}{\varepsilon^{n-\frac{1}{2}}}\Big)^2$.
Thus, assuming that $k_n<\tfrac{2}{2+c_n}$, we easily conclude
that $V=0$. $\Box$
\end{proof}
%
%
%
\begin{thm}
Let $u$ be the solution of \eqref{TNT1} and
$(U_h^n)_{n=0}^{\ssy N}$ be the fully discrete approximations
that the method \eqref{Reactive100}-\eqref{Reactive101} produces.
If $\varepsilon(t)>0$ for $t\in[0,T]$,
then, there exists a constant $C_{\ssy R}\ge0$ such that: if
$\displaystyle{\max_{1\leq{n}\leq{\ssy
N}}}(k_n\,C_{\ssy R})\leq\tfrac{1}{3}$, there exists a constant $C>0$
such that
\begin{equation}\label{Lestimation100}
\max_{1\leq{n}\leq{\ssy N}}\|U_h^n-u^n\|_1\leq\,C
\,(k^2+h^r)\,\,\Upsilon_{\ssy R}(u),
\quad\forall\,h\in(0,h_{\star}],
\end{equation}
where
$$\Upsilon_{\ssy
R}(u):=\sum_{\ell=0}^1\max_{\ssy[0,T]}\|\partial_t^{\ell}u\|_{r+1}
+\sum_{\ell=2}^3\max_{\ssy[0,T]}|\partial_t^{\ell}u|_1
+\sum_{m=0}^1\max_{\ssy t\in[0,T]}|\partial_t^2\partial_x^mu(t,1)|
+\max_{\ssy [0,T]}|\partial_t^2u|_2.$$
\end{thm}
%
%
\begin{proof}
Let $h\in(0,h_{\star}]$,
$\theta_h^n:=U_h^n-R_h^{\star}u^n$ and
$\eta^n:=R_h^{\star}u^n-u^n$ for $n=0,\dots,N$.
Use \eqref{Reactive101}, \eqref{CNPRB_ConsistencyEquations} and
\eqref{RDEF2}, to obtain
\begin{equation}\label{MagnaRea1}
\begin{split}
{\mathcal B}(\partial\theta_h^n,\chi)=&\,
\left\{\,\tfrac{a^{n-\frac{1}{2}}}{\varepsilon^{n-\frac{1}{2}}}
\,({\mathcal A}(\partial_x\theta_h^n)(1)
-\left[\,\tfrac{\delta^{n-\frac{1}{2}}}{\varepsilon^{n-\frac{1}{2}}}
+\beta(t,1)\right]\,({\mathcal A}\theta_h^{n})(1)
+E^n_3+E_4^n\,\right\}\,\chi'(1)\\
&-a^{n-\frac{1}{2}}\,{\mathcal B}^{\star}({\mathcal A}\theta_h^n,\chi)
+{\mathcal B}(\beta^{n-\frac{1}{2}}
\,{\mathcal A}\theta_h^n,\chi)\\
&+{\mathcal B}(E_1^n-\sigma^n,\chi)
+a^{n-\frac{1}{2}}\,{\mathcal B}^{\star}(E_2^n,\chi)
+a^{n-\frac{1}{2}}\,{\mathcal B}({\mathcal A}\eta^n,\chi)
\quad\forall\,\chi\in S_h,
\quad n=1,\dots,N,\\
\end{split}
\end{equation}
where $\sigma^n:\overline{D}\rightarrow\mathbb{R}$ is defined by
\eqref{CNPRB_ConsistencyEquations} and for $n=1,\dots,N$
\begin{equation*}
\begin{split}
E_1^n&:=\,\partial u^n-R^{\star}_h(\partial
u^n)-\beta^{n-\frac{1}{2}}\,
{\mathcal A}(u^n-R_h^{\star}u^n),\\
E_2^n&:=\,u(t^{n-\frac{1}{2}})-{\mathcal A}u^n,\\
E_3^n&:=\,\tfrac{a^{n-\frac{1}{2}}}{\varepsilon^{n-\frac{1}{2}}}
\,{\mathcal A}(\partial_x\eta^n(1))
-\left(\tfrac{\delta^{n-\frac{1}{2}}}{\varepsilon^{n-\frac{1}{2}}}
+\beta^{n-\frac{1}{2}}(1)\right)\,{\mathcal A}(\eta^n(1)),\\
E_4^n&:=\,\tfrac{a^{n-\frac{1}{2}}}{\varepsilon^{n-\frac{1}{2}}}
\,\left[\,{\mathcal A}(u_x(t^n,1))
-u_x(t^{n-\frac{1}{2}},1)\,\right]
-\left(\tfrac{\delta^{n-\frac{1}{2}}}{\varepsilon^{n-\frac{1}{2}}}
+\beta^{n-\frac{1}{2}}(1)\right)\,\left[{\mathcal A}(u(t^n,1))
-u(t^{n-\frac{1}{2}},1)\right].\\
\end{split}
\end{equation*}
\par
Using Taylor's formula, \eqref{StarEll1}
and \eqref{StarEll4}, we derive the following bounds:
\begin{equation}\label{MagnaRea2}
\left|{\mathcal B}(E_1^n-\sigma^n,\chi)
+a^{n-\frac{1}{2}}\,{\mathcal B}({\mathcal A}\eta^n,\chi)\right|
\leq\,C\,\left[\,|\sigma^n|_1+
h^{r}\,\left(\,\max_{\ssy[t^{n-1},t^n]}\|u\|_{r+1}
+\max_{\ssy[t^{n-1},t^n]}\|\partial_tu\|_{r+1}
\,\right)\,\right]\,|\chi|_1,
\end{equation}
\begin{equation}\label{ReactiveFD2}
|\sigma^n|_1\leq\,C\,(k_n)^2\,\,
\left(\,\max_{\ssy[t^{n-1},t^n]}|\partial_t^2u|_1
+\max_{\ssy[t^{n-1},t^n]}|\partial_t^3u|_1\,\right),
\end{equation}
\begin{equation}\label{MagnaRea3}
\left|a^{n-\frac{1}{2}}\,{\mathcal B}^{\star}(E_2^n,\chi)\right|
\leq\,C\,k_n^2\,\,\max_{\ssy[t^{n-1},t^n]}|\partial_t^2u|_2\,|\chi|_2,
\end{equation}
\begin{equation}\label{MagnaRea4}
\big|\,E_3^n\,\chi'(1)\,\big|\leq\,C\,h^{r+\frac{1}{2}}
\,\max_{\ssy[t^{n-1},t^n]}\|u\|_{r+1}\,|\chi'(1)|,
\end{equation}
and
\begin{equation}\label{MagnaRea5}
|E_4^n\,\chi'(1)|\leq\,C\,k_n^2\,\left(\,
\max_{\ssy t\in[t^{n-1},t^n]}|\partial_t^2u(t,1)|
+\max_{\ssy t\in[t^{n-1},t^n]}|\partial_t^2\partial_xu(t,1)|\,\right)
\,|\chi'(1)|
\end{equation}
for $n=1,\dots,N$ and $\chi\in{\mathcal S}_h$.
\par
Now, set $\chi={\mathcal A}\theta_h^n$ in \eqref{MagnaRea1}
and use \eqref{PoincareF}, the Cauchy-Schwarz inequality
and the estimates \eqref{MagnaRea2},
\eqref{MagnaRea3}, \eqref{MagnaRea4} and \eqref{MagnaRea5},
to obtain
\begin{equation}\label{MagnaRea6}
\begin{split}
\tfrac{1}{2k_n}\,\left(\,|\theta_h^n|_1^2-|\theta_h^{n-1}|_1^2\,\right)
\leq&\,-a_{\star}\,|{\mathcal A}\theta_h^n|_2^2
+C\,\Big[\,|{\mathcal A}\theta_h^n|_1^2
+k^2\,\max_{\ssy[0,T]}|\partial_t^2u|_2\,|{\mathcal A}\theta_h^n|_2\\
&\hskip2.9truecm
+(k^2+h^r)^2\,(\Upsilon_{\ssy R}(u))^2\\
&\hskip2.9truecm
+|\partial_x({\mathcal A}\theta_h^n)(1)|^2
+|({\mathcal A}\theta_h^n)(1)|^2\,\Big],
\quad n=1,\dots,N.\\
\end{split}
\end{equation}
After use of the trace inequalities \eqref{TraceIneq} and
\eqref{TraceIneq2}, of the inequality \eqref{PoincareF} and of
the arithmetic-geometric mean inequality, \eqref{MagnaRea6}
yields the existence of a constant $C_{\ssy R}>0$ such that
\begin{equation}\label{MagnaRea7}
(1-C_{\ssy R}\,k_n)\,|\theta_h^n|_1^2
\leq\,(1+C_{\ssy R}\,k_n)\,|\theta_h^{n-1}|_1^2
+C\,k_n\,(k^2+h^r)^2\,(\Upsilon_{\ssy R}(u))^2,
\quad n=1,\dots,N.
\end{equation}
Assuming that
$\displaystyle{\max_{1\leq{n}\leq{\ssy N}}}(C_{\ssy R}\,k_n)
\leq\,\tfrac{1}{3}$, and following a discrete Gr{\"o}nwall
argument similar to that of Proposition~\ref{MagnaProposition}
we arrive at
\begin{equation}\label{FinalRea}
\max_{0,\leq{n}\leq{\ssy N}}|\theta_h^n|_1
\leq\,C\,(k^2+h^r)\,\Upsilon_{\ssy R}(u).
\end{equation}
Thus the desired estimate \eqref{Lestimation100} follows easily
if we combine \eqref{FinalRea}, \eqref{StarEll1} and
\eqref{PoincareF}. $\Box$
\end{proof}
%
%
%
%
\subsection*{\bf Acknowledgments} This work was supported by
a {\tt Pythagoras} grant to the University of Athens,
co-funded by the E.U. European Social Fund and the Greek
Ministry of Education. The authors would like to thank
Ms. Evangelia Flouri for her help with the numerical
experiments using IFD.
%
%
%
%
\bibliographystyle{plain}
\end{document}